\documentclass[a4paper, 10pt]{article}
\usepackage[toc,page]{appendix}
\usepackage[bottom]{footmisc}
\usepackage{microtype}
\usepackage{tocbibind}
\usepackage{stmaryrd}   
\usepackage{amsmath}
\usepackage{verbatim}
\usepackage{mathrsfs}
\usepackage{mathtools}
\usepackage{amsthm} 
\usepackage{amssymb} 
\usepackage{mdframed}

\usepackage{faktor}
\usepackage[colorlinks={true},linkcolor={blue},citecolor={blue}, urlcolor={blue} ]{hyperref}
\usepackage{braket}
\usepackage{transparent}
\usepackage{graphicx}
\usepackage{color}
\usepackage{rotating}
\usepackage{tikz}
\usepackage{tikz-cd}
\usepackage[margin=2.5 cm]{geometry}

\usepackage{tikz}
\usepackage{tikz-cd}
\usetikzlibrary{intersections}
\usetikzlibrary{spy} 
\usetikzlibrary{decorations.pathmorphing}
\usetikzlibrary{decorations.pathreplacing}

\DeclareMathOperator{\spec}{Spec}

\DeclareMathOperator{\fr}{Frac}

\DeclareMathOperator{\res}{res}
\DeclareMathOperator{\tr}{Tr}
\DeclareMathOperator{\Pic}{Pic}
\DeclareMathOperator{\Div}{Div}

\DeclareMathOperator{\Cres}{Cres}
\DeclareMathOperator{\sep}{sep}
\DeclareMathOperator{\Hom}{Hom}

\DeclareMathOperator{\ar}{Ar}
\DeclareMathOperator{\ara}{Ar}

\DeclareMathOperator{\ideg}{ideg}

\DeclareMathOperator{\ord}{ord}
\DeclareMathOperator{\divi}{div}
\DeclareMathOperator{\Princ}{Princ}
\DeclareMathOperator{\CH}{CH}
\DeclareMathOperator{\cts}{cts}

\DeclareMathOperator{\MyProd}{\scalebox{1.4}{$\mathrm{I\kern-0.2ex I}$}}

\theoremstyle{plain}

\newtheorem{theorem}{Theorem}[section]
\newtheorem{lemma}[theorem]{Lemma}
\newtheorem{proposition}[theorem]{Proposition}

\newtheorem{corollary}[theorem]{Corollary}
\theoremstyle{definition}
\newtheorem{definition}[theorem]{Definition}

\theoremstyle{remark}

\newtheorem{remark}[theorem]{Remark}

\newcommand{\sbt}{\,\begin{picture}(-1,1)(-1,-3)\circle*{3}\end{picture}\ }

\newcommand{\catname}[1]{\mathbf{#1}}

\makeatletter

\title{Adelic geometry on arithmetic surfaces II: completed adeles and idelic Arakelov intersection theory}

\author{
  Weronika Czerniawska
  \and
  Paolo Dolce
}
\date{}

\newcommand{\Addresses}{{
  \bigskip
  \footnotesize

  W.~Czerniawska, \textsc{University of Nottingham}\par\nopagebreak
  \textit{E-mail address}: 
  \texttt{Weronika.Czerniawska1@nottingham.ac.uk}

  \medskip

  P.~Dolce, \textsc{University of Nottingham}\par\nopagebreak
  \textit{E-mail address}: \texttt{Paolo.Dolce1@nottingham.ac.uk}

}}

\begin{document}

\maketitle

\begin{abstract}
We work with completed adelic structures on an arithmetic surface and justify that the construction under consideration is compatible with Arakelov geometry. The ring of completed adeles is algebraically and topologically self-dual and fundamental adelic subspaces are self orthogonal with respect to a natural differential pairing. We show that the Arakelov intersection pairing can be lifted to an idelic intersection pairing. 
\end{abstract}

\makeatletter
\@starttoc{toc}
\makeatother

\setcounter{section}{-1}
\section{Introduction}
\subsection{Background}
 Adelic theory for global fields was introduced for the first time by Chevalley in the 1930's as a tool for studying the completions of a number field with respect to all possible absolute values at the same time. This is a great expression of ``local-to-global'' principles as well as an example of geometric approaches to number theory which have proven to be very powerful. One of  the principal applications of adelic theory for number fields was published in John Tate's  thesis \cite{tatethesis} which presented a proof of meromorphic continuation and functional equation of $\zeta$ functions of number fields in  clearer and more compact way than the proof given before by Hecke. When $C$ is a curve over a perfect field, one can define the adelic ring $\mathbf A_{C}$ associated to $C$ as the restricted product of the complete discrete valuation fields $K_c$ for any closed point $c\in C$ with respect to their valuation rings $\mathcal O_c$. It is possible to obtain a very elegant proof of the Riemann-Roch theorem for curves by using adeles (see for example \cite[0.]{fe0}).

 Adelic approach  has been generalized for higher dimensions by Beilinson in \cite{bei} where he defined adelic structures as functors on the category of quasi-coherent sheaves. An explicit  theory of $2$-dimensional adelic cohomology and dualities for algebraic surfaces was outlined in \cite{pa1}, where hope for proving adelic Riemann-Roch theorem for a surface over a finite field was expressed. However, the explicit adelic structures introduced in \cite{pa1} are not equivalent to Beilinson's, since \cite{pa1} worked with objects that now are called rational adeles. The gap on the definitions was partially fixed in \cite{pa2}, but a complete account of $2$-dimensional explicit adelic theory was given by Fesenko in \cite{fe0}, where he also proved  an adelic  Riemann-Roch theorem for an algebraic  surface over a perfect field  by using properties of adelic cohomology. In particular, Fesenko showed that the function field of an algebraic surface $X$ can be seen as a discrete subspace inside the ring of $2$-dimensional adeles attached to $X$. Such a result generalizes the classical result of \cite{tatethesis} which shows that a global field is a discrete object inside the ring of adeles. 

The non-cohomological part of (explicit) adelic theory for algebraic surfaces can be summarized in the following way: fix a nonsingular algebraic surface $(X,\mathscr O_X)$ over a perfect field $k$, then to each ``flag''  $x\in y$ made of a closed point $x$ inside an integral curve $y\subset X$ we can associate the ring $K_{x,y}$ which will be a $2$-dimensional local field if $y$ is nonsingular at $x$, or a finite product of $2$-dimensional local fields if we have a singularity. Note how the geometric dimension of $X$ matches  the ``dimension'' of the ring $K_{x,y}$, and this happens roughly speaking because for a flag $x\in y$ (assuming that $x$ is a nonsingular point of $y$) we have two distinct levels of discrete valuations: there is the discrete valuation associated to the containment $x\in y$ and the discrete valuation associated to $y\subset X$. $K_{x,y}$ is obtained through a process of successive localisations and completions starting with $\mathscr O_{X,x}$ and by the symbol $\mathcal O_{x,y}$ we denote the product of valuation rings inside $K_{x,y}$. The step to the global theory is obtained  by performing a ``double restricted product'' of the rings $K_{x,y}$: first over all points ranging on a fixed curve  and then over all curves in $X$, in order to obtain the $2$-dimensional adelic ring:
$$\mathbf A_X:=\sideset{}{''}\prod_{\substack {x\in y\\ y\subset X}} K_{x,y}\subset\prod_{\substack {x\in y\\ y\subset X}} K_{x,y} \,.$$
The topology on $K_{x,y}$ can be defined canonically thanks to the construction by completions and localisations, and by starting with the standard $\mathfrak m_x$-adic topology on $\mathscr O_{X,x}$. The topology on $\mathbf A_X$ is obtained after a process of several inductive an projective limits by starting from the local topologies on all $K_{x,y}$. In \cite{fe0} it is shown that  $\mathbf A_X$ is self-dual as $k$-vector space. For $2$-dimensional local fields with the same structure of $K_{x,y}$ there is a well known theory of differential forms and residues (e.g. \cite{yek1}); one can globalize the constructions in order to obtain a $k$-character  $\xi^{\omega}:\mathbf A_X\to k$ associated to a rational differential form $\omega\in\Omega^1_{k(X)|k}$ and the differential pairing:
\begin{eqnarray*}
d_\omega:\mathbf A_X\times\mathbf A_X &\to & k\\
(\alpha,\beta) &\mapsto& \xi^\omega(\alpha\beta)\,.
\end{eqnarray*} 
Fesenko in \cite{fe0} proves that the subspace $\mathbf A_X/k(X)^{\perp}$ is a linearly compact $k$-vector space (orthogonal spaces are calculated with respect to $d_\omega$) and the function field  $k(X)$ is discrete in $\mathbf A_X$. It is possible to define some important subspaces of $\mathbf A_X$ denoted as:
$k(X)=A_0$, $A_{1}$, $A_2$, $A_{01}$, $A_{02}$, $A_{12}$, $A_{012}=\mathbf A_X$ which generate an idelic complex assuming the following form:
\begin{equation*}
\begin{tikzcd}
\mathcal A^\times_X: & A^\times_0\oplus A^\times_1\oplus A^\times_2 \arrow[r, "d_\times^0"]& A^\times_{01}\oplus A^\times_{02}\oplus A^\times_{12}\arrow[r, "d_\times^1"]& A^\times_{012}\\
\end{tikzcd}
\end{equation*}
It can be shown that the space $\ker(d^1_\times)$ is a generalization of the group $\Div(X)$ since there is a surjective map  $\ker(d^1_\times)\to \Div(X)$ and the intersection pairing on $\Div(X)$ can be extended to a pairing on $\ker(d^1_\times)$ (cf. \cite[3.]{dolce_phd}).

The main aim of our work is to obtain a two-dimensional adelic theory, for arithmetic surfaces i.e. objects of the form $\varphi:X\to\spec O_K$ where $K$ is a number field. The problem is motivated by Fesenko's ``analysis on arithmetic schemes programme''. The programme develops a two-dimensional generalization of Tate's thesis, i.e. two-dimensional measure, integration and Fourier analysis. Fesenko's work reveals relationships between geometry and analysis not visible without adelic tools (see also \cite{wero_phd} for an alternative presentation). 

In \cite{MM3} and \cite{MM2} Morrow,  develops an explicit approach to residues and dualizing sheaves of arithmetic surfaces. In particular he defines the residue map for $2$-dimensional local fields arising from an arithmetic surface and he formulates and proves reciprocity laws around a point and along a curve of an arithmetic surface. To have a reciprocity law along a horizontal curve, he completes horizontal curves with points at infinity, i.e. real or complex embeddings of the function field of the horizontal curve.

\subsection{Results in this paper}\label{results}
At the center of our considerations there is an adelic object for an arithmetic surface $\varphi:X\to B=\spec O_K$. One expects that one has to take into account (archimedean)   ``data at infinity'' of the arithmetic surface. Such an adelic space completed by data at infinity was proposed for the first time in \cite{fe1}. In section \ref{comp_adeles} we present a simpler and slightly different version of it. Already at the level of local theory, adelic geometry for arithmetic surfaces is quite interesting, in fact the rings $K_{x,y}$ can be equal characteristic or mixed characteristic $2$-dimensional local fields depending whether $y$ is horizontal or vertical. Over each point at infinity $\sigma\in B_\infty$, i.e. an embdedding $\sigma: K\to\mathbb C$, we obtain, by a base change, a Riemann surface $X_\sigma$ that can be thought as a fibre at infinity. The completed  adelic ring $\mathbf A_{\widehat X}$ will then contain the one dimensional adelic rings $\mathbf A_{X_\sigma}$ relative to the fibres at infinity $X_\sigma$, but counted twice:
$$
\mathbf A_{\widehat X}= \mathbf A_{X}\oplus\prod_{\sigma\in B_\infty}(\mathbf A_{X_\sigma}\oplus \mathbf A_{X_\sigma})\,.
$$
The arithmetic counterparts $A_{\widehat \ast}$ of the fundamental subspaces $A_{\ast}$ are also defined.
There is a specific geometric reason that suggests why we should count adeles at infinity twice, and it involves the interpretation horizontal curves on $\widehat X$ in terms of Arakelov geometry i.e. we have to consider their ``intersection'' with fibres at infinity.

By slightly generalizing the local theory of residues for two dimensional local fields developed in \cite{MM3}, in section \ref{res_sect} we define a global adelic residue 
$$\xi^\omega:\mathbf A_{\widehat X}\to \mathbb T$$
($\omega$ is a fixed nonzero rational differential form and $\mathbb T$ is the unit complex circle) and we show that $\xi^\omega$ is sequentially continuous. 

Section \ref{self} is entirely dedicated to the proof of the self-duality of $\mathbf A_{\widehat X}$ as topological additive group. In particular we show that $\mathbf A_{\widehat X}\cong \widehat{\mathbf A_{\widehat X}}$ as topological groups and moreover that there is a character $\psi:\mathbf A_{\widehat X}\to\mathbb T$ such that any other character of $\mathbf A_{\widehat X}$ is of the form $\psi(a\cdot)$ for $a\in \mathbf A_{\widehat X}$. 

 In section \ref{prop_diff} we define the arithmetic differential pairing \begin{eqnarray*}
d_{\omega}\colon \mathbf A_{\widehat X}\times\mathbf A_{\widehat X}&\to&\mathbb T\\
(\alpha,\beta)&\mapsto& \xi^{\omega}(\alpha\beta)\,.
\end{eqnarray*}
We improve the reciprocity laws proved in \cite{MM2} by  giving a set of ``completed'' reciprocity laws, i.e. taking into account all flags coming from points at infinity. We show that both $A_{\widehat{01}}$ and $A_{\widehat{02}}$ (adelic subspaces corresponding to curves and points respectively) are self-orthogonal with respect to $d_\omega$ i.e. $A_{\widehat{01}}=A_{\widehat{01}}^\perp$ and  $A_{\widehat{02}}=A_{\widehat{02}}^\perp$. The inclusions $A_{\widehat{01}}\subseteq A_{\widehat{01}}^\perp$ and  $A_{\widehat{02}}\subseteq A_{\widehat{02}}^\perp$ are a direct consequence of the completed reciprocity laws, thus  the self-orthogonality of $A_{\widehat{01}}$ and $A_{\widehat{02}}$ can be interpreted as ``strong reciprocity laws'' for arithmetic surfaces. The ``strong reciprocity laws'' for surfaces over a perfect field were proved in \cite{fe0}.
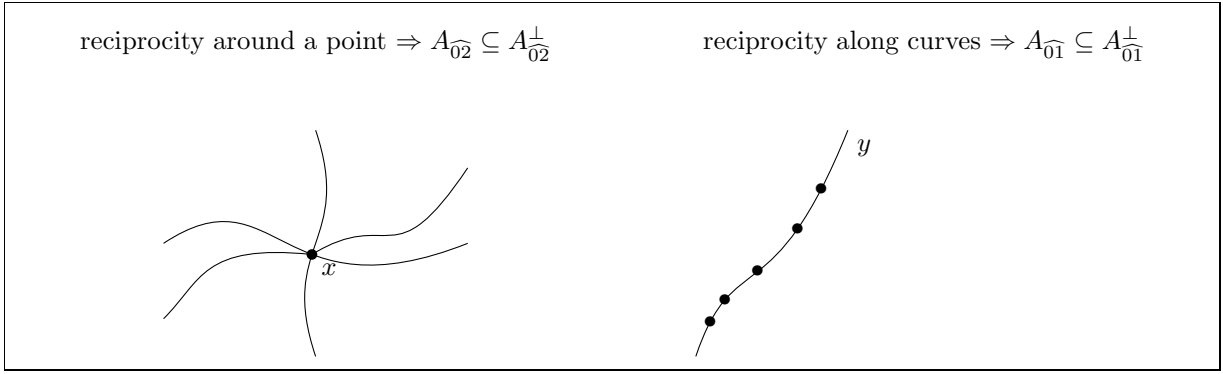
\begin{figure}[htp!]
\begin{mdframed}
\centering
\begin{tikzpicture}

\draw [name path=c1](2, -1.5)..controls (1.5,0) and (2.5,0).. (2,1.5);
\draw [name path=c2](0, 0)..controls (1.5,1) and (1.5,-1).. (4,0);
\fill [name intersections={
of=c1 and c2,  by={b}}]
(b) circle (2pt)node[anchor=north west]{$x$};

\draw[name path=c3] (0,-1)..controls (.5,-.5) and (.5,0)..(b);

\draw[name path=c4] (b)..controls (3,0.5) and (3,-.5)..(4,1);

\draw (7, -1.5)..controls (7.5,0) and (8,-1).. (9,1.5)
node[near end] {$\bullet$}
node[very near end] {$\bullet$}
node[near start] {$\bullet$}
node[very near start] {$\bullet$}
node[midway] {$\bullet$};

\draw (9,1.5) node[anchor=north west]{$y$};

\draw (10,3) node[anchor=north]{reciprocity along curves $\Rightarrow A_{\widehat{01}}\subseteq A_{\widehat{01}}^\perp$ };

\draw (2,3)node[anchor=north]{reciprocity around a point $\Rightarrow A_{\widehat{02}}\subseteq A_{\widehat{02}}^\perp$ };
\end{tikzpicture}
\end{mdframed}
\caption{\footnotesize{The sum of local two-dimensional residues when a point is fixed and curves passing through it vary is zero. The sum of local two-dimensional residues when a curve is fixed and the points sitting on in vary is zero.}}
\end{figure}
\noindent

The problems of finding proofs of the discreteness of the function field $K(X)$ inside $\mathbf A_{\widehat X}$ and of the compactness of the quotient $\mathbf A_{\widehat X}/K(X)^\perp$ are still open, but we plan to publish a solution in a forthcoming paper. Finally, in analogy with the case of algebraic surfaces we show that the Arakelov intersection pairing can be lifted to the idelic group $\mathbf A^\times_X$. The schematic part of the lifting was already proved in \cite{adI}, so here we solve the problem of the data carried by Green functions on fibres at infinity. It is worth remembering that Arakelov theory is the only known theory that provides consistent intersection theory on arithmetic surfaces, therefore we would expect that a theory of adeles on arithmetic surfaces should resonate with Arakelov geometry. 

The text contains also two appendices which are indispensable for the understanding  and moreover prerequisites for this paper are \cite{adI} and a basic knowledge of the theory of higher local fields (e.g. \cite{MM1}).

\paragraph{Acknowledgements.} The authors want to express their gratitude to \emph{I. Fesenko} and  \emph{E. Lepage} for their suggestions and their kindness in answering our several questions on the topic. Special thanks go also to \emph{O. Braunling}, \emph{R. Waller}, \emph{W. Porowski} and \emph{A. Minamide}.

This work is supported by the EPSRC programme grant EP/M024830/1 (Symmetries and correspondences: intra-disciplinary developments and applications).

\section{Preliminaries}
\subsection{Basic notions}

\paragraph{General notations}
All rings are considered commutative and unitary. Let $(A,\mathfrak m)$ be a Noetherian local ring and let $M$ be an $A$-module, then we put $M^{\sep}:=M\big/\cap_{j\ge 1} \mathfrak m^jM$. When we pick a point $x$ in a scheme $X$ we generally mean a \emph{closed point} if not otherwise specified, also all sums $\sum_{x\in X}$ are meant to be ``over all closed points of $X$''. The cardinality of a set $T$ is denoted as $\#(T)$. If $F$ is a discrete valuation  field, then $\overline{F}$ doesn't denote the algebraic closure but its residue field. In particular if $a\in F$ then $\overline a$ is the image of $a$ in $\overline F$. For a morphisms of schemes $f:X\to S$, the schematic preimage of $s\in S$ is $X_s$. Sheaves are denoted with the ``mathscr'' late$\chi$ font;in particular the structure sheaf of a scheme $X$ is $\mathscr O_X$ (note the difference with the font $\mathcal O$). With the symbol $\mathbb T$ we denote the unit circle in the complex plane. The superscript $\,\widehat{\;}\,$ is used several times in this paper to denote completely different objects: the dual of a topological group, the completion of a local ring or a ``completed structure'' in the framework of Arakelov geometry. This superposition of notation is harmless because the specific meaning of $\,\widehat{\;}\,$ will be clear from the context.

\paragraph{Topological groups.} If not otherwise specified we assume that any topological group is abelian and Hausdorff. The dual of a topological group $G$ is the group of (unitary) characters:
$$\widehat{G}:=\Hom_{\operatorname{cont}}(G,\mathbb T).$$
It is a topological group endowed with the compact-to-open topology. 
Moreover for a compact subset $C\subset G$ and an open  $U\subset \mathbb T$ neighborhood of $1$ we denote 
$$\mathcal W(C,U)=\left\{\chi\in\widehat G\colon \chi(C)\subset U\right\}\subset\widehat G\,.$$
The sets of the type  $\mathcal W(C,U)$ form an open base at $1$ for the compact-to-open topology in $\widehat G$.

If $G$ is algebraically and topologically isomorphic to $\widehat G$, then we say that \emph{$G$ is self-dual}. If  $G$ is also a ST ring and $\xi\in \widehat{G}$ is a nontrivial character, then for any $a\in G$ the map 
\begin{eqnarray*}
\xi_a: G &\to& \mathbb T\\
x&\mapsto& \xi(ax)
\end{eqnarray*}
is a character. If the map
 \begin{eqnarray*}
\Theta_\xi\colon G&\to &  \widehat G\\
a &\mapsto &\xi_a
\end{eqnarray*}
is an algebraic and topological isomorphism for any $a\in G$, we say that $\xi$ is a \emph{standard character}. For any subsets $S\subseteq G$  and $R\subseteq \widehat G$ we put:
$$S^{\perp}:=\{\chi\in\widehat{G}\colon \chi(S)=1\}\subseteq \widehat{G}\,,$$
$$R^{\perp}:=\{g\in G \colon \chi(g)=1,\;\forall \chi\in R\}\subseteq G\,.$$
If $H$ is a subgroup of $G$, we say that $H$ is \emph{dually closed} if  for every element $ g\in G\setminus H$, there is a character  $\psi\in H^\perp$ such that  $\psi(g)\neq 1$.

We will often use the following simple general result:

\begin{proposition}\label{easycomp}
Let $G$ be a topological group such that $G=\varinjlim_{r\in\mathbb Z} H_i$ where $H_i\subset G$ is a subgroup and $H_i\supset H_{i+1}$  for any $i\in \mathbb Z$. Then any compact subset $C\subset G$ is contained in some $H_i$. 
\end{proposition} 
\proof
For $i\le 0$ put $n=-i$ and consider the sequence of subgroups $\{H_n\}_{n\in \mathbb N}$ with $H_n\subset H_{n+1}$ and $G=\bigcup_n H_n$. Assume that the claim is false and choose a sequence of points $x_n\in C\cap (H_n\setminus H_{n-1})$. Put $A=\{x_n\}_{n\ge0}$. Then if $B\subseteq A$, then $B\cap H_n$ is finite for each $n$, so since points are closed in $H_n$, $B\cap H_n$ is closed in $H_n$. This means that $B$ is closed in $G$. In particular, $A$ is a closed subset of $G$, and every subset of $A$ is closed so it has the discrete topology. But a closed subset of a compact space is compact, and a compact discrete space must be finite. This is a contradiction with the construction of $A$.
\endproof
\subsection{Geometric setting}\label{setting}
Let's fix the main objects and notations that we will use throughout the whole paper. Some of the material contained in this section can be found with more details in \cite{adI}.  In particular we assume that the reader is familiar with the notion of $2$-dimensional local field. Moreover, topological aspects of this section rely on appendix \ref{ap_ST}.

Let $K$ be a number field with ring of integers $O_K$. Fix the arithmetic surface $\varphi: X\to B=\spec O_K $ which is a $B$-scheme with the following properties:

\begin{enumerate}
\item[\sbt] $X$ is  two dimensional, integral, and regular. The generic point of $X$ is $\eta$ and the function field of $X$ is denoted by $K(X)$.
\item[\sbt] $\varphi$ is proper and flat.
\item[\sbt] The generic fibre, denoted by $X_K$, is a geometrically integral, smooth, projective curve over $K$. The generic point of $B$ is denoted by $\xi$.
\end{enumerate}
It is well known that $\varphi$ is a projective morphism, so in particular also $X$ is projective (see \cite[Theorem 8.3.16]{Liu}). Let's recall a useful result which characterizes all points of dimension $1$ on $X$:
\begin{proposition}\label{hori_chara}
If $x$ is a closed point of the curve $X_K$, then $\overline{\{x\}}$ is a horizontal (prime) divisor in $X$. Vice versa if $D$ is a prime divisor on $X$, then either $D\subseteq X_b$ for a closed point $b\in B$ or $D=\overline{\{x\}}$ where $x$ is a closed point of $X_K$.
\end{proposition}
\proof
See for example \cite[Proposition 8.3.4]{Liu}.
\endproof
Let $B_\infty$ be the set of field embeddings $\sigma\colon K\hookrightarrow \mathbb C$ up to conjugation, its  cardinality is $[K:\mathbb Q]$ and the \emph{completion of $B$} is the set $\widehat{B}:=B\cup B_\infty$. 
For any point (i.e. nonzero prime ideal) $b=\mathfrak p\in B$ we put:
\begin{itemize}
\item[\sbt] $\mathcal O_b:=\widehat{\mathscr O_{B,b}}$. It is a complete DVR.
\item[\sbt] $K_b:=\fr \mathcal O_b$. It is a local field with finite residue field. The valuation is denoted by $v_b$.
\end{itemize}
From now on, we \emph{always} fix a set of representatives in $B_\infty$. Therefore  $B_\infty$ is simply a finite set of embeddings viewed as points at infinity of $B$.  For the non-archimedean place associated to $b=\mathfrak p\in B$, on $K$ we choose the absolute value
$$|\cdot|_{b}:=\mathfrak N(\mathfrak p)^{-v_{b}(\cdot)}$$
where $\mathfrak N(\mathfrak p)$ is the cardinality of $O_K/\mathfrak p$.
Moreover:
\begin{itemize}
\item[\sbt] For any real embedding $\tau:K\to \mathbb R$ we consider the absolute value:
$$|\cdot|_{\tau}:=|\tau(\cdot)|$$
where on the right hand side we mean the usual absolute value on $\mathbb R$. In this case we define the real valuation associated to $\tau$ as 
$$v_\tau(\cdot):=-\log|\cdot|_{\tau}$$
\item[\sbt] For any couple of conjugate embeddings $\sigma,\overline{\sigma}:K\to \mathbb C$ we choose:
$$|\cdot|_{\sigma}:=|\sigma(\cdot)|$$
where on the right hand side we have the usual absolute value on $\mathbb C$\footnote{Many authors in this case take the square of the complex absolute value  to keep track of the fact that point at infinity induced by $|\cdot|_\sigma$ is ``complex'', so roughly speaking ``of order two". We will fix this by using the coefficient $2$ when necessary.}. Note that $|\cdot|_\sigma$ doesn't depend on the choice between $\sigma$ and $\overline \sigma$, since they give the same absolute value. The associated real valuation is 
$$v_\sigma(\cdot):=-\log|\cdot|_{\sigma}\,.$$
\end{itemize}
For $\sigma\in B_\infty$, $K_\sigma$ is the completion of $K$ with respect to $|\cdot|_{\sigma}$, thus $K_\sigma=\mathbb C$ or $K_\sigma=\mathbb R$. Furthermore, let's introduce a  constant, associated to each $\sigma\in B_{\infty}$:
$$\epsilon_{\sigma}:=\begin{cases} 1 & \mbox{if $\sigma$ is real}   \\ 2 & \mbox{if $\sigma$ is complex\;.}  \end{cases}$$
The adelic ring of $\widehat B$ (or equivalently of the number field $K$) is denoted by $\mathbf A_{\widehat B}$, whereas $\mathbf A_{B}:= \mathbf A_{\widehat B}\cap\prod_{b\in B} K_b$. For any $\sigma\in B_\infty$ consider the base change diagram:
\begin{equation}\label{diag1}
\begin{tikzcd}
X_\sigma:=X\times_{B}\spec \mathbb C\arrow{r}\arrow{d}{\varphi_\sigma}&\spec \mathbb C\arrow{d}{\spec\sigma} \\
X\arrow{r}{\varphi}& B\rlap{\ .}
\end{tikzcd}
\end{equation}
By the properties of the fibred product, it turns out that $X_\sigma\to \spec \mathbb C$ is a complex integral (integrality is a consequence of the geometrical integrality of $X_K$), regular projective curve. We denote the function field of $X_\sigma$ by the symbol $\mathbb C(X_\sigma)$.
\begin{remark}
Diagram (\ref{diag1}) arises from the following rather obvious commutative diagram:
$$
\begin{tikzcd}
X_\sigma\arrow["\beta"]{d} \arrow["\varphi_\sigma"',bend right]{dd}\arrow{r}& \spec \mathbb C \arrow["\spec \sigma"]{d}\\
X_K\arrow{r}\arrow{d} & \spec K\arrow["\spec \iota"]{d}\\
X\arrow{r} & B
\end{tikzcd}
$$
where $\iota:O_K\hookrightarrow K$ is the natural embedding and the map $\beta$ is surjective. In other words $\varphi_\sigma$ maps surjectively $X_\sigma$ onto the curve $X_K$. Since the morphisms $\iota$ and $\sigma$ are both flat and flatness is preserved after base change, we can conclude that $\varphi_\sigma$ is flat.
\end{remark}
With the notation $\widehat {X}$, we define the ``completed surface'' 
$$\widehat {X}:=X\cup\bigcup_{\sigma\in B_\sigma} X_\sigma\,.$$
A curve $Y$ on $X$ will always be an integral curve and its unique generic point will be denoted with the letter $y$. For simplicity we will often identify  $Y$ with its generic point $y$, which means that by an abuse of language  and notation we will use sentences like ``let $y\subset X$ be a curve on $X$...". A \emph{flag} on $X$ is a couple $(x,y)$ where $x$ is a closed  point sitting on a curve $y\subset X$, it will be denoted simply as $x\in y$. 
\begin{definition}
Fix a closed point $x\in X$, then:
\begin{itemize}
\item[\sbt] $\mathcal O_x:=\widehat{\mathscr O_{X,x}}$. It is a Noetherian, complete, regular, local, domain of dimension $2$ with maximal ideal $\widehat{\mathfrak m_x}$.
\item[\sbt] $K'_x:=\fr\mathcal O_x$.
\item[\sbt] $K_x:=K(X)\mathcal O_x\subseteq K'_x$. 
\end{itemize}
For a curve $y\subset X$ we put:
\begin{itemize}
\item[\sbt] $\mathcal O_y:=\widehat{\mathscr O_{X,y}}$. It is a complete DVR  with maximal ideal $\widehat{\mathfrak m_y}$.
\item[\sbt] $K_y:=\fr\mathcal O_y$. It is a complete discrete valuation field with valuation ring $\mathcal O_y$. The valuation is denoted by $v_y$.
\end{itemize}
\end{definition}

For a flag $x\in y\subset X$, we have a surjective local homomorphism $\mathscr O_{X,x}\to\mathscr O_{y,x}$, with kernel $\mathfrak p_{y,x}$, induced by the closed embedding $y\subset X$ (note that $\mathfrak p_{y,x}$ is a prime ideal of height $1$). The inclusion $\mathscr O_{X,x}\subset \mathcal O_x$ induces a morphism of schemes $\varphi\colon\spec\mathcal O_x\to\spec\mathscr O_{X,x}$ and we define the \emph{local branches of $y$ at $x$} as the elements of the set
$$y(x):=\varphi^{-1}(\mathfrak p_{y,x})=\set{\mathfrak z\in\spec\mathcal O_x\colon \mathfrak z\cap \mathscr O_{X,x}=\mathfrak p_{y,x}}\,.$$
If $y(x)$ contains only an element, we say that $y$ is unbranched at $x$. Fix a flag $x\in y\subset X$ with $\mathfrak z\in y(x)$, then we have the $2$-dimensional local field $$K_{x,\mathfrak z}:=\fr\left(\widehat{\left(\mathcal O_x\right)_\mathfrak z}\right)$$
explicitly obtained in the following way: we localise $\mathcal O_x$ at the prime ideal $\mathfrak z$,  complete it at its maximal ideal and finally we take the fraction field. The ring of integers of $K_{x,\mathfrak z}$ is denoted by $\mathcal O_{x,\mathfrak z}:=\mathcal O_{K_{x,\mathfrak z}}=\widehat{\left(\mathcal O_x\right)_\mathfrak z}$. All the needed material about higher local fields is contained in \cite[1.1]{adI}, whereas for a deeper study the reader can consult  \cite{fesbook2}; see also a more recent introduction in \cite{MM1}.
\begin{definition}
Let $x\in y\subset X$ be a flag and let $\mathfrak z\in y(x)$, then the first residue field of $K_{x,\mathfrak z}$ is $E_{x,\mathfrak z}:=K^{(1)}_{x,\mathfrak z}$ and the second residue field is  $k_\mathfrak z (x):=K^{(2)}_{x,\mathfrak z}$. The valuation on $K_{x,\mathfrak z}$ is $v_{x,\mathfrak z}$ and the valuation on $E_{x,\mathfrak z}$ is $v^{(1)}_{x,\mathfrak z}$; whereas $\mathcal O_{K_{x,\mathfrak z}}^{(2)}:=\{a\in\mathcal O_{x,\mathfrak z}\colon \overline a\in\mathcal O_{E_{x,\mathfrak z}}\}$. 
$$
\begin{tikzcd}
K_{x,\mathfrak z}\arrow[r, phantom, "\supset"]\arrow[dr, bend right,dashed, no head
] & \mathcal O_{x,\mathfrak z}:=\mathcal O_{K_{x,\mathfrak z}}\arrow[d]\arrow[r, phantom, "\supset"]& \mathcal O_{K_{x,\mathfrak z}}^{(2)}\arrow[d]\\
& E_{x,\mathfrak z}:=K^{(1)}_{x,\mathfrak z}\arrow[r, phantom, "\supset"]\arrow[dr, bend right,dashed, no head
] & \mathcal O_{E_{x,\mathfrak z}}\arrow[d]\\
& & k_\mathfrak z(x):=K^{(2)}_{x,\mathfrak z}
\end{tikzcd}
$$
Moreover we put:
$$K_{x,y}:=\prod_{\mathfrak z\in y(x)}K_{x,\mathfrak z}\,,\quad\mathcal O_{x,y}:=\prod_{\mathfrak z\in y(x)}\mathcal O_{x,\mathfrak z}\,,$$
$$E_{x,y}:=\prod_{\mathfrak z\in y(x)}E_{x,\mathfrak z}\,,\quad k_y(x):=\prod_{\mathfrak z\in y(x)}k_\mathfrak z (x)\,.$$
\end{definition}
Let's endow  $\mathscr O_{X,x}$ with the $\mathfrak m_x$-adic topology with respect to its maximal ideal, then $K_{x,\mathfrak z}$ can be endowed with a canonical topology by using the following steps explained in appendix \ref{indpro}:

\begin{equation}\label{top_geo_2loc}
\begin{tikzcd}
\mathscr O_{X,x}\arrow[r,squiggly,"(C)"]& \mathcal O_x=\widehat{\mathscr O_{X,x}} \arrow[r,squiggly,"(L)"]& (\mathcal O_x)_{\mathfrak z} \arrow[r,squiggly,"(C)"]& \widehat{(\mathcal O_x)_{\mathfrak z}} \arrow[r,squiggly,"(L)"]& K_{x,\mathfrak z}=\fr\left(\widehat{\left(\mathcal O_x\right)_\mathfrak z}\right)\,.
\end{tikzcd}
\end{equation}
Then $K_{x,y}$ is endowed with the product topology and it is a ST ring (see appendix \ref{ap_ST} for an introduction to semi-topological structures). Here it is very important to point out that $K_{x,y}$ is not a topological ring, since it turns out that the multiplication is not continuous as function of two variables.

\begin{remark}
This is one of the several ways to topologise $K_{x,y}$; see for example \cite[1.]{Brau} for a survey. It is not the most explicit topology for $K_{x,y}$, but it is independent from the choice of the uniformizing parameter since it is obtained by a general process of localizations and completions.
\end{remark}

If $y$ is a horizontal curve then $K_{x,\mathfrak z}$ is of equal characteristic and isomorphic to $E_{x,\mathfrak z}((t))$ where  $E_{x,\mathfrak z}$ is a finite extension of $\mathbb Q_p$ and $t$ is  (the image of) a uniformizing parameter. If $y$ is a vertical curve then  $K_{x,\mathfrak z}$ is of mixed characteristic and isomorphic to a finite extension of $K_p\{\!\{t\}\!\}$ where $K_p$ is a finite extension of $\mathbb Q_p$ (see  \cite[example 1.7]{adI} for the definition of $K_p\{\!\{t\}\!\}$). In this case $t$ it is not (the image of) a uniformizing parameter, but it is (the image of) a uniformizing parameter for  $E_{x,\mathfrak z}\cong\overline{K_p}((t))$. It is always possible to choose a uniformizing parameter of $K_y$ to be also the uniformizing parameter of  $K_{x,\mathfrak z}$ for all $x\in y$.

If $\varphi(x)=b$ we have an embedding $K_b\hookrightarrow K_{x,\mathfrak z}$, and we say that $K_{x,\mathfrak z}$ is  an arithmetic $2$-dimensional local field over $K_b$. The module of differential forms relative to $x$ and $\mathfrak z\in y(x)$ is  the  $K_{x,\mathfrak z}$-vector space:
$$\Omega^1_{x,\mathfrak z}:=\left(\Omega^{1}_{\mathcal O_{x,\mathfrak z}|\mathcal O_b}\right)^{\sep} \otimes_{\mathcal O_{x,\mathfrak z}} K_{x,\mathfrak z}\,,$$
where $\Omega^{1}_{\mathcal O_{x,\mathfrak z}|\mathcal O_b}$ is the usual module of K\"ahler differential forms and the operator ``$\sep$'' was defined  at the end of section \ref{results} in the ``General notations'' paragraph.
Then, $\Omega^1_{x,\mathfrak z}$ is endowed with the  vector space topology over $K_{x,\mathfrak z}$.  In \cite{MM3} and \cite{MM2} it is defined the residue map:
$$\res_{x,\mathfrak z}: \Omega^1_{x,\mathfrak z}\to K_b$$
with the following properties:
\begin{itemize}
    \item[\sbt] It is $K_b$-linear.
    \item[\sbt] It is continuous (this is shown in \cite[Lemma 2.8, Remark 2.9]{MM2}).
\end{itemize}
A more detailed description of $\Omega^1_{x,\mathfrak z}$ and $\res_{x,\mathfrak z}$ will be given in section \ref{res_sect}.

The global adelic theory for the projective scheme $X$ is described in \cite[1.2]{adI}. We obtain the adelic ring $\mathbf A_X$ as a ``double restricted product'' of the rings $K_{x,y}$ performed first over closed points ranging on curves, and then over all curves in $X$. Fix any curve $y\subset X$ and  denote by $\mathfrak J_{x,y}$ the Jacobson radical of $\mathcal O_{x,y}$; we put

\begin{equation*}
\mathbb A^{(0)}_y := \left\{ 
  \begin{aligned}
  & (\alpha_{x,y})_{x\in y}\in\prod_{x\in y}\mathcal O_{x,y} \colon  \forall s >0,\; \alpha_{x,y}\in\mathcal O_x+\mathfrak J_{x,y}^s \\ 
  & \text{for all but finitely many $x\in y$.}
  \end{aligned}
\right\}\subset\prod_{x\in y} \mathcal O_{x,y}\,.
\end{equation*}
Then for any $r\in\mathbb Z$ and for any choice of uniformizing parameter $t_y$
$$\mathbb A^{(r)}_y:= \widehat{\mathfrak m}^{r}_{y}\mathbb A_y^{(0)}=t_y^r\mathbb A_y^{(0)}\subset \prod_{x\in y} K_{x,y}\,.$$
Clearly  $\mathbb A^{(r)}_y\supseteq\mathbb A^{(r+1)}_y$ and $\bigcap_{r\in\mathbb Z} \mathbb A^{(r)}_y=0$; moreover we  define 
$$\mathbb A_y:=\bigcup_{r\in\mathbb Z} \mathbb A^{(r)}_y\,.$$
\begin{definition}
The ring of adeles of $X$ is
$$\mathbf A_X:=\left\{(\beta_y)_{y\subset X}\in\prod_{y\subset X}\mathbb A_y\colon \beta_y\in\mathbb A^{(0)}_y\; \text{for all but finitely many $y$}\right\}\subset\prod_{\substack {x\in y,\\ y\subset X}} K_{x,y}\,,$$
\end{definition}

Finally we recall the definitions of  some important subspaces of $\mathbf A_X$. Consider the following diagonal embeddings:$$K_x\subset\prod_{y\ni x} K_{x,y}, \quad K_y\subset\prod_{x\in y} K_{x,y}\,,$$
so  we can put:
$$\prod_{x\in X}K_x\subset\prod_{\substack {x\in y\\ y\subset X}} K_{x,y}, \quad \prod_{y\subset X}K_y\subset\prod_{\substack {x\in y\\ y\subset X}} K_{x,y}\,.$$
then we define
$$A_{012}:=\mathbf A_X\,;\quad A_{12}:=\mathbf A_X\cap \prod_{\substack{x\in y\\y\subset X}}\mathcal O_{x,y}=\prod_{y\subset X} \mathbb A^{(0)}_y; $$
$$A_{02}:=\mathbf A_X\cap\prod_{x\in X} K_x\,;\quad A_{2}:=\mathbf A_X\cap\prod_{x\in X} \mathcal O_x\,;\quad A_{01}:=\mathbf A_X\cap\prod_{y\subset X} K_y\,;$$
$$\quad A_{1}:=\mathbf A_X\cap\prod_{y\subset X} \mathcal O_y\,;\quad A_0:=K(X)\,.$$
The subspaces satisfy a series of inclusion relations depicted in the following diagram:
$$
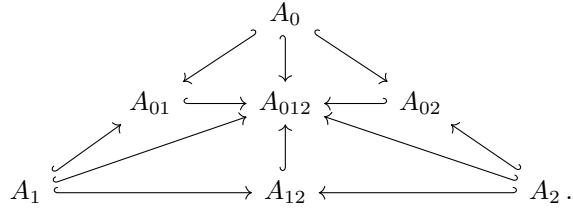
\begin{tikzcd}
&&A_0\arrow[d,hook']\arrow[dl,hook']\arrow[dr,hook]&&\\
&A_{01}\arrow[r,hook]& A_{012} &A_{02}\arrow[l,hook']&\\
A_1\arrow[rr,hook]\arrow[ur,hook]\arrow[urr,hook]&&A_{12}\arrow[u,hook]&&A_2\arrow[ll,hook']\arrow[ul,hook']\arrow[ull,hook']\,.\\
\end{tikzcd}
$$
When $X$ is an algebraic surface over a perfect field $k$, the algebraic and topological properties of the subspaces $A_\ast$ were studied in \cite{fe0}.

\subsection{Topology on adelic structures.}\label{thetop} In this crucial subsection we explain how to put a topology on all adelic structures introduced so far. We point out that all categorical limit considered here are in the category linear topological groups (so linear direct/inverse limits). For more details see appendix \ref{ap_ST}. 

\begin{itemize}
\item[\sbt] For any $s>0$ let's put:
 $$\mathbb A^{(0)}_y\{s\}:=\{(a_{x,y})_{x\in y}\in\prod_{x\in y} \mathcal O_{x,y}\colon a_{x,y}\in\mathcal O_x+\mathfrak J^s_{x,y} \text{ for all but finitely many} x\in y\}\,.$$
 Endow $\mathbb A^{(0)}_y\{s\}$ with the  restricted product topology (i.e. linear direct limit).
\item[\sbt] $\mathbb A^{(0)}_y=\bigcap_{s\ge 0} A^{(0)}_y\{s\}$, so we put on  $\mathbb A^{(0)}_y$ the linear inverse limit topology.
\item[\sbt] The topology is transferred from  $\mathbb A^{(0)}_y$  to $\mathbb A^{(r)}_y$ for any $r\in \mathbb Z$, by the multiplication by $t^r_y$. 
\item[\sbt] Each $\mathbb A^{(r)}_y/\mathbb A^{(r+j)}_y$, for $j>0$, is endowed with the quotient topology.
\item[\sbt] We endow $\mathbb A_y=\varinjlim \mathbb A^{(r)}_y=\bigcup_r \mathbb A^{(r)}_y$ with the linear direct limit topology. 
\item[\sbt] $\mathbf A_X$ is the restricted product (seen as linear direct limit) of the topological groups $\mathbb A_y$ with respect to $\mathbb A^{(0)}_y$. 
\end{itemize}
\noindent
Since $\mathcal O_x+\mathfrak J_{x,y}$ surjects onto $E_{x,y}$, it is easy to see that the natural projection (which is continuous and open)
\begin{eqnarray*}
p_y:\mathbb A^{(0)}_y &\to& \mathbf A^f_{k(y)}\\
(a_{x,y})_{x\in y} &\mapsto&  (\overline{a_{x,y}})_{x\in y}
\end{eqnarray*}
induces an algebraic and topologic isomorphism betweent $\mathbb A^{(0)}_y/\mathbb A^{(1)}_y$ and the ring of the one dimensional finite adeles $\mathbf A^f_{k(y)}$.  Consider the exact sequence:
$$0\to \mathbb A_y^{(1)}/\mathbb A_y^{(2)}\to \mathbb A_y^{(0)}/\mathbb A_y^{(2)}\to \mathbb A_y^{(0)}/\mathbb A_y^{(1)}\to 0$$
Since $\mathbb A_y^{(1)}/\mathbb A_y^{(2)}$ and $\mathbb A_y^{(0)}/\mathbb A_y^{(1)}$ are locally compact and self-dual, then   $\mathbb A_y^{(0)}/\mathbb A_y^{(2)}$ is locally compact. We conclude that for any $j>0$ the quotient $\mathbb A_y^{(r)}/\mathbb A_y^{(r+j)}$ is a locally compact topological group (hence complete).
\begin{proposition}\label{comple/neigh}
The following two fundamental topological properties hold:
\begin{itemize}
\item[$(i)$] $\mathbb A^{(r)}_{y}$ is complete for any $r\in\mathbb Z$ (but in general is not locally compact).
\item[$(ii)$] For each open neighborhood $U\subset\mathbb A^{(r)}_y$  of $0$ there is $s>r$ such that $\mathbb A^{(s)}_y\subset U$.
\end{itemize}
\end{proposition}
\proof $(i)$ is true  since $\mathcal O_{x,y}$ is complete and $\mathcal O_x+\mathfrak J^s_{x,y}$ is closed in $\mathcal O_{x,y}$. $(ii)$ can be checked directly from the above definition of the topology. 
\endproof

\begin{proposition}
There is an algebraic and topological isomorphism
$$\mathbb A^{(r)}_y\cong\varprojlim_{j>0}\mathbb A^{(r)}_y/\mathbb A^{(r+j)}_y$$
\end{proposition}
\proof
Thanks to proposition \ref{comple/neigh}, we can apply directly  \cite[III \S7.3, Corollary 1]{bou}.
\endproof
In particular $\mathbb A^{(r)}_y\cong t^r_y\mathbb A^{(0)}_{y}[[t_y]]$ and  any Laurent power series in $t^r_y\mathbb A^{(0)}_{y}[[t_y]]$ is a truly convergent series.  The open subgroups of $\mathbb A_y$ that form a local basis at $0$ can be described in the following way: fix  a sequence $\{U_i\}_{i\in\mathbb Z}$ of open sets in $\mathbb A^{(0)}_y$ with the property that there exists $k\in\mathbb Z$ such that $U_{i}=\mathbb A^{(0)}_y$ for $i\ge k$.  Then we consider the open set 
$$\sideset{}{'}\sum U_it_y^i:=\left\{\text{Laurent series  $\sum a_jt_y^j$ such that $a_j\in U_j$}\right \}\,.$$
Each open neighborhood $U\subset\mathbb A_y$ of $0$ contains some $\mathbb A^{(r)}_y$.
 
\section{The ring of completed adeles \texorpdfstring{$\mathbf A_{\widehat X}$}{} and its subspaces}\label{comp_adeles}

We want to define adeles for arithmetic surfaces in a way that preserves the most fundamental properties of the adelic theory and is compatible with Arakelov geometry. In particular, we have to consider points at infinity of the base and, corresponding to them, infinite fibres.  When we add a fibre at infinity $X_\sigma$ to the picture, we have to take in account \emph{all} possible flags on the completed surface $\widehat X$: a point $p$ on a fibre at infinity $X_\sigma$ originates a flag  $p\in X_\sigma$, but it can be seen also as an ``intersection point'' between a completed horizontal curve $\overline y$ and $X_\sigma$. 

Let $y$ be a curve on $X$, if $y$ is vertical then we put $\overline y=y$, if $y$ is horizontal, then by $\overline y$ we mean:
$$\overline y=y\cup \!\bigcup_{\sigma\in B_\infty} \!y_\sigma$$
where 
\begin{equation*}
y_\sigma=\varphi_\sigma^\ast(y)\in \Div(X_\sigma)\,.
\end{equation*}
By simplicity we also put $y_\infty:=\cup_{\sigma\in B_\infty} y_\sigma $, so we have the decomposition $\overline y= y\cup y_\infty$. Any point $p\in X_\sigma$ lies on a  completed horizontal curve $\overline y$ because we have the map $\varphi_\sigma: X_\sigma\to X_K\subset X$  and points of the generic fibre $X_K$ are in bijective correspondence with horizontal curves. From now on, a curve on $\widehat X$ will be always a completed curve $\overline y$, and a point $x\in \overline y$ can be also a point lying on some ``part at infinity'' $y_\sigma$  (when $y$ is horizontal), if not explicitly said otherwise.

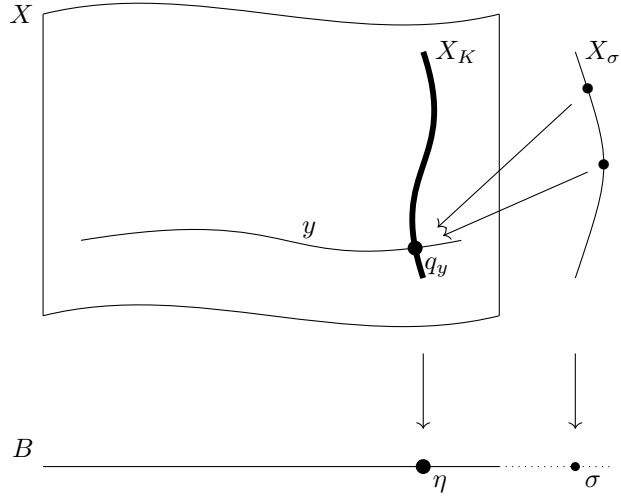
\begin{figure}[htp!]
\centering
\begin{tikzpicture}

\draw (-3,2)node[anchor=east]{$X$}..controls (-1,2.5) and (1,1.5)..(3,2);
\draw (-3,2)--(-3,-2);
\draw (-3,-2)..controls (-1,-1.5) and (1,-2.5)..(3,-2);
\draw (3,2)--(3,-2);

\draw[name path=c2](-2.5,-1)..controls (0.5,-0.5) and (0,-1.5)..(2.5,-1);

\draw (0.5,-1.1)node[anchor=south]{$y$};

\draw [name path=c1,line width=2.2pt](2, -1.5)..controls (1.5,0) and (2.5,0).. (2,1.5)node[anchor=west]{$X_K$};

\draw(4, -1.5)..controls (4.5,0) and (4.5,0).. (4,1.5)node[anchor=west]{$X_\sigma$}
node(a1)[very near end] {$\bullet$}
node(a2)[midway] {$\bullet$};

\fill [name intersections={
of=c1 and c2,  by={b}}]
(b) circle (2.8pt)node[anchor=north west]{$q_y$};

\draw (-3,-4)node[anchor= south east]{$B$}--(3,-4);
\draw[dotted] (3,-4)--(4.5,-4);
\draw (2,-4) circle [radius=2.2pt];
\fill (2,-4) circle [radius=2.8pt]node[anchor=north west]{$\eta$};

\draw [->] (2,-2.5)--(2,-3.5);
\draw [->] (4,-2.5)--(4,-3.5);
\draw [->,  shorten >=0.4cm] (a1)--(b);
\draw [->,  shorten >=0.4cm] (a2)--(b);

\draw (4,-4) circle [radius=1.5pt];
\fill (4,-4) circle [radius=1.5pt]node[anchor=north west]{$\sigma$};

\end{tikzpicture}
\caption{\footnotesize{A visual example where $ y_\sigma$ is made of two points (marked on the curve $X_\sigma$). The generic point of the curve $y$ here is denoted by $q_y$.}}
\end{figure}
\noindent
The local data of the completed adelic ring will be the following ones:

\begin{itemize}
    \item[\sbt] For any flag at infinity $p \in  X_\sigma$ we put
$$K_{p,\sigma}:=\fr \left (\widehat{\mathscr O_{X_\sigma,p}}\right )\,.$$
In other words $K_{p,\sigma}$ is a local field isomorphic to $\mathbb C((t))$. The valuation ring of $K_{p,\sigma}$ is $\mathcal O_{p,\sigma}\cong \mathbb C[[t]]$ and  $E_{p,\sigma}\cong \mathbb C$ is the residue field. 
    \item[\sbt] If $p\in \overline y$ and $p\in y_\sigma$ for some $\sigma\in B_{\infty}$, we put 
$$K_{p,\overline y}:= K_{p,\sigma},\quad\mathcal O_{p,\overline y}:= \mathcal O_{p,\sigma},\quad E_{p,\overline y}:= E_{p,\sigma}\,;$$
    \item[\sbt] For any other point $x\in y$ we have:
$$K_{x,\overline y}:= K_{x,y},\quad \mathcal O_{x,\overline y}:= \mathcal O_{x,y},\quad E_{x,\overline y}:= E_{x,y},\quad k_{\overline y}(x)=k_y(x)\,.$$
\end{itemize}
\noindent
When $p$ is a point at infinity we want to consider the fields  $K_{p,\sigma}$ and $K_{p,\overline y}$ as $2$-dimensional local fields, but if we use a completion/localization topology as described in equation (\ref{top_geo_2loc}), we obtain the usual one dimensional valuation topology. Therefore we fix some isomorphisms $K_{p,\sigma}\cong K_{p,\overline y}\cong \mathbb C((t))$ (parametrizations), we consider $\mathbb C$ with the standard topology given by its archimedean norm, and we endow $\mathbb C((t))$ with the ind/pro-topology (see appendix \ref{indpro}). Then we carry such a topology on  $K_{p,\sigma}$ and $K_{p,\overline y}$ through the parametrizations. The ind/pro-topology on $\mathbb C((t))$ is coarser than the 1-dimensional valuation topology. 

\begin{remark}\label{almostadeles}
In \cite{fe1} the construction of local fields at infinity is slightly different, indeed $K_{p,\sigma}$ is $\mathbb R((t))$ or $\mathbb C((t))$, depending whether $\sigma$ is a real or complex embedding. This might seem a very natural choice, but in the framework of Arakelov geometry $X_\sigma$ is always  a Riemann surface, even if $\sigma$ is real. We want to build deep link between Arakelov geometry and adelic geometry, therefore we prefer to put $K_{p,\sigma}\cong\mathbb C((t))$.
\end{remark}

\begin{remark}
In the product $\prod_{\substack {x\in \overline y,\\ \overline y\subset \widehat X}}K_{x,\overline y}$ we find three different types of 2-dimensional local fields: $K_p((t))$, finite extensions of $K_p\{\!\{t\}\!\}$ and  $\mathbb C((t))$.
\end{remark}

We are going to define a new ring $\overline{\mathbf A}_{X}$ which will be a subspace of the big product $\prod_{\substack {x\in \overline y,\\ \overline y\subset \widehat X}}K_{x,\overline y}$. Let's first extend the spaces $\mathbb A^{(r)}_y$ for completed curves: 
\begin{definition}
For any completed curve $\overline y$ let's put:
$$\mathbb{A}_{\overline y}:=\mathbb{A}_y\oplus\prod_{p\in y_\infty} K_{p,\overline y}\,,$$
$$\mathbb{A}^{(0)}_{\overline y}:=\mathbb{A}^{(0)}_y\oplus\prod_{p\in y_\infty} \mathcal O_{p,\overline y}\,,$$
$$\mathbb{A}^{(r)}_{\overline y}:=\mathbb{A}^{(r)}_y\oplus\prod_{p\in y_\infty} \mathfrak p^r_{K_{p,\overline y}}\mathcal O_{p,\overline y}\,,$$
and endow them with the finite product topology.
\end{definition}
Again each  $\mathbb{A}^{(r)}_{\overline y}$ is closed in $\mathbb A_{\overline y}$ and the latter can be thought as a first restricted product performed on the completed curve $\overline y$. We can use the formal notation:
$$\mathbb A_{\overline y}=\sideset{}{'}\prod_{x\in \overline y} K_{x,\overline y}\,.$$
Let's assume by simplicity that $y$ is a regular horizontal curve, then $K_{x,y}\cong E_{x,y}((t))$ where $E_{x,y}$ is a finite extension of $\mathbb Q_p$ and  it is the completion of the field $k(y)$ with respect to the valuation induced by the inclusion $x\in y$. Moreover  $y=\spec O_L$ where $L$ is a finite extension of $K$. In general if $y$ is any horizontal curve admitting singular points, then $y=\spec R$ where $R$ is an order of $L$.
\noindent
For any curve $y$ we put 
$$\mathbf A_{\overline y}:=\sideset{}{'}\prod_{x\in y}k(y)_x\oplus\prod_{q\in y_{\infty}}\mathbb C\,,$$
where the restricted product is with respect to the complete discrete valuation rings corresponding to the points $x\in y$. In other words $\mathbf A_{\overline y}$ is in general slightly bigger than the classical $1$-dimensional adelic ring of $\overline y$. If the point $q\in y_{\infty}$ is present (recall that in the case of vertical curves there is no archimedean data) and corresponds to a real embedding $\sigma$, then the ``$q$-component'' of $\mathbf A_{\overline y}$ is $\mathbb C$ and not $\mathbb R$, i.e. we take $\mathbb{C}$ for all archimedean places. The finite part of $\mathbf{A}_{\overline y}$, denoted by  $\mathbf{A}^f_{\overline y}$, and the finite part of classical $1$-dimensional adeles coincide.  This of course descends from our choice of data at infinity (see remark \ref{almostadeles}), but all adelic properties of $\mathbf A_{\overline y}$ are clearly the same of the one dimensional adeles. In particular all results of \cite{tatethesis} hold for $\mathbf{A}_{\overline y}$. 

\begin{lemma}\label{explemma}
Let $y$ be a regular horizontal curve. For any $r\in\mathbb Z$,  $\mathbb A_y^{(r)}$ is equal to the following ring:
$$
\Xi^{(r)}_y:=\left\{(\alpha_{x,y})_{x\in y}\in\prod_{x\in y} K_{x,y}\colon \alpha_{x,y} \text{ satisfies the following conditions $(\ast)$ and $(\ast\ast)$}\right\}
$$
\begin{itemize}
\item[$(\ast)$] $\alpha_{x,y}\in t^rE_{x,y}[[t]]$.
\item[$(\ast\ast)$] Assume that: 
$$\alpha_{x,y}= t^r\sum_{i\ge 0}\Gamma_{x,i} t^i\qquad \text{with }\; \Gamma_{x,i}\in E_{x,y}\,,$$
then for any fixed index $i$ the sequence $(\Gamma_{x,i})_{x\in y}\in \mathbf A^{f}_{\overline y}$. In other words for all but finitely many $x\in y$ we have that $\Gamma_{x,i}\in\mathcal O_{E_{x,y}}$.
\end{itemize}
\end{lemma}
\proof
\emph{Inclusion $\mathbb A^{(r)}_y\subseteq \Xi^{(r)}_y$}. Let's start with $r=0$, the general case will follow trivially. Consider an element $(\alpha_{x,y})_{x\in y}$, then clearly $(\ast)$ is true because $\mathcal O_{x,y}=E_{x,y}[[t]]$. Suppose that
$\alpha_{x,y}=\sum_{i\ge 0}\Gamma_{x,i} t^i$, then there exists a  decomposition:
$$\alpha_{x,y}=\sum_{i\ge 0}\Theta_{x,i} t^i+\sum_{i\ge 0}\Lambda_{x,i} t^i \in\mathcal O_x+\mathcal O_{x,y}$$
where $\Theta_{x,i}\in\mathcal O_{E_{x,y}}$, $\Lambda_{x,i}\in E_{x,y}\setminus\mathcal O_{E_{x,y}}$, and $\Gamma_{x,i}=\Theta_{x,i}+\Lambda_{x,i}$. Now fix an index $h\ge 0$, then the set
$$S_h:=\left\{x\in y\colon \Lambda_{x,h}\neq 0\right\}$$
is finite, indeed note that $\mathcal O_x+\mathfrak J^s_{x,y}=\mathcal O_{E_{x,y}}[[t]]+t^s E_{x,y}[[t]]$, thus if $\Lambda_{x,h}\neq 0$, then $\alpha_{x,y}\not\in\mathcal O_x+\mathfrak J^{h+1}_{x,y}$. In other words if for infinitely many $x\in y$ we had that  $\Lambda_{x,h}\neq 0$, then for the same points $\alpha_{x,y}\not\in\mathcal O_x+\mathfrak J_{x,y}^{h+1}$ against the definition of $\mathbb A^{(0)}_y$. We have shown that for all but finitely many $x\in y$, $\Gamma_{x,i}=\Theta_{x,i}\in  \mathcal O_{E_{x,y}}$ which is equivalent to say that $(\Gamma_{x,i})_{x\in y}\in \mathbf A^f_{\overline y}$.\\
The case when $r\neq 0$ follows easily from the fact that $\widehat{\mathfrak m}^r_y\Xi^{(0)}_y=\Xi^{(r)}_y$.\\\\
\emph{Inclusion $\Xi^{(r)}_y\subseteq \mathbb A^{(r)}_y$}. As above it is enough to write the proof for $r=0$. Let $(\alpha_{x,y})_{x\in y}\in\Xi^{(0)}_y $, then for any index $i\ge 0$ define:
$$T_i:=\left\{x\in y\colon \Gamma_{x,i}\not\in \mathcal O_{E_{x,y}}[[t]]\right\}\,;$$
by the property $(\ast\ast)$ $T_i$ is a finite set. Now fix an index $h> 0$ then for all $x\in y\setminus \cup_{i=1}^{h-1}T_i$, (i.e. for all but finitely many $x\in y$) it holds that $\Gamma_{x,i}=\Theta_{x,i}$ when $i< h $, which means that
$$\alpha_{x,y}=\sum_{i\ge 0}\Theta_{x,i} t^i+\sum_{i\ge h}\Lambda_{x,i} t^i \in\mathcal O_x+\mathfrak J^h_{x,y}\,.$$ 
\endproof

\begin{proposition}\label{adelepol}
Let $y$ be a regular horizontal curve. For any $r\in\mathbb Z$, $\mathbb A_{\overline y}^{(r)}\cong t^r\mathbf A_{\overline y}[[t]]$ as rings (here $t$ is simply a variable). In particular $\mathbb A_{\overline y}\cong \mathbf A_{\overline y}((t))$ and $\mathbb A^{(0)}_{\overline y}\cong \mathbf A_{\overline y}[[t]]$. 
\end{proposition}
\proof
By lemma \ref{explemma} we have the equality $\mathbb A_y^{(r)}=\Xi^{(r)}_y$ and the map $\Xi^r_y\to t^r\mathbf A^f_{\overline y}[[t]]$ is given in the following way  and it is well defined:
$$
(\alpha_{x,y})_{x\in y}=\left(t_2^r\sum_{i\ge 0}\Gamma_{x,i}t_2^i\right)_{x\in y}\mapsto  t^r\sum_{i\ge 0}(\Gamma_{x,i})_{x\in y}\, t^i\,.
$$
It is routine check to show that is is a ring isomorphism.
\endproof

\begin{remark}\label{remsing}
Proposition \ref{adelepol} is true also when $y$ is a singular curve. The proof is based on a slightly modified version of lemma \ref{explemma}; the only difference consists in the fact that if $x\in y$ is singular then  $K_{x,y}=\prod_{\mathfrak z\in y(x)} K_{x,\mathfrak z}$ is a sum of 2-dimensional valuation fields and $\mathfrak J_{x,y}$ is the sum of the maximal ideals of  $K_{x,\mathfrak z}$. Here we restricted the proof to the case of non-singular curves just by simplicity of notations.
\end{remark}

\begin{definition}
The modified version of $\mathbf A_X$ which takes in account the completed curves is:
$$\overline{\mathbf A}_X:=\left\{(\beta_{\overline y})_{\overline y\subset \widehat X}\in\prod_{\overline y\subset \widehat X}\mathbb A_{\overline y}\colon \beta_{\overline y}\in\mathbb A^{(0)}_{\overline y}\; \text{for all but finitely many $\overline y$}\right\}\subset\prod_{\substack {x\in \overline y,\\ \overline y\subset \widehat X}} K_{x,\overline y}\,.$$
We also  introduce the formal notation
$$\overline{\mathbf A}_X=\sideset{}{''}\prod_{\substack {x\in \overline y\\ \overline y\subset \widehat X}} K_{x,\overline y}\,.$$
The topology on $\overline{\mathbf A}_{X}$ is the restricted topology of the additive groups $\mathbb A_{\overline y}$ with respect to $\mathbb{A}^{(0)}_{\overline y}$.
\end{definition}
\begin{definition}
The \emph{completed adelic ring} attached to $\widehat X$ is
$$\mathbf A_{\widehat X}:=\overline{\mathbf A}_{X}\oplus\prod_{\sigma\in B_\infty}\mathbf A_{X_{\sigma}}$$
where  each $\mathbf A_{X_\sigma}$ is the adelic ring of the Riemann surface $X_\sigma$. The topology on $\mathbf A_{\widehat X}$ is the product topology. 
\end{definition}

Let $\Upsilon$ be the collection of all finite sets of completed curves of $\widehat X$, then for $S\in\Upsilon$ we define

$$\mathbf A_{\widehat X}(S):=\prod_{\overline y\in S} \mathbb A_{\overline y}\times \prod_{\overline y\notin S} \mathbb A^{(0)}_{\overline y}\times \prod_{\sigma\in B_{\infty}}\mathbf A_{X_\sigma}$$
then:

$$\bigcup_{S\in \Upsilon}\mathbf A_{\widehat X}(S)=\mathbf A_{\widehat X}\,,\quad \bigcap_{S\in \Upsilon}\mathbf A_{\widehat X}(S)=\prod_{\overline y\subset \widehat X}\mathbb A^{(0)}_{\overline y}\times \prod_{\sigma\in B_{\infty}}\mathbf A_{X_\sigma}\,.$$
The following proposition establishes a nice relationship between $\mathbf A_{\widehat X}$ and $\mathbf A_X$.
\begin{proposition}
The following equality holds:
$$\mathbf A_{\widehat X}= \mathbf A_{X}\oplus\prod_{\sigma\in B_\sigma}(\mathbf A_{X_\sigma}\oplus \mathbf A_{X_\sigma})\,.$$
\end{proposition}
\proof
Let $\alpha\in \overline{\mathbf A}_{X}$, then it can be decomposed in the following way:
$$\alpha=(a_y)_{y\subset X}\times(a_{p,\sigma})_{\substack {p\in X_{\sigma},\\ \sigma\in B_\infty}}$$
where:
\begin{itemize}
\item[\sbt] $a_y\in\mathbb A_y$ for all $y\subset X$ and $a_y\in\mathbb A^{(0)}_y$ for all but finitely many $y$.
\item[\sbt] For any fixed $\sigma$ we have $a_{p,\sigma}\in K_{p,\sigma}$ and $a_{p,\sigma}\in \mathcal O_{p,\sigma}$ for all but finitely many $p\in X_\sigma$.
\end{itemize}
This means that $\alpha\in \overline{\mathbf A}_{X}\subseteq \mathbf A_X\oplus \prod_{\sigma\in B\sigma}\mathbf A_{X_\sigma}$, so obviously 
$$\mathbf A_{\widehat X}\subseteq \mathbf A_{X}\oplus\prod_{\sigma\in B_\sigma}(\mathbf A_{X_\sigma}\oplus \mathbf A_{X_\sigma})\,.$$
Vice versa, let $\alpha \in\mathbf A_X\oplus \prod_{\sigma\in B\sigma}\mathbf A_{X_\sigma}$ then:
$$\alpha=(a_y)_{y\subset X}\times(a_{p,\sigma})_{\substack {p\in X_{\sigma},\\ \sigma\in B_\infty}}$$
where $a_y$ and $a_{p,\sigma}$ satisfy the conditions listed above. Since each $\varphi_\sigma: X_\sigma\to X_K$ is surjective and points of $X_K$ correspond to horizontal curves on $X$, we can write easily:
$$\alpha=(a_y)_{y\subset X}\times(a_{p,\sigma})_{\substack {p\in X_{\sigma},\\ \sigma\in B_\infty}}=(a_y)_{y\subset X}\times\left((a_{p,\sigma})_{p\in y_\infty}\right)_{y_\infty\subset X_\infty}=(a_{\overline y})_{\overline y\subset\widehat X}\in \overline{\mathbf A}_{X}\,.$$

\endproof

\begin{remark}
The above definition of $\mathbf A_{\widehat X}$ is new, but the object is very similar to the ring of  completed adeles given in \cite{fe1}. One difference was already emphasized in remark \ref{almostadeles}; moreover in \cite[25.]{fe1} the spaces $\mathbb A_y^{(r)}$ are obtained through some local lifting maps of $E_{x,y}$ to $\mathcal O_{x,y}$.
\end{remark}

\begin{remark}
At first glance, one might think that a reasonable definition of the adelic ring  $\mathbf A_{\widehat X}$  can be just $\mathbf A_X\oplus\prod_{\sigma\in B_\infty} \mathbf A_{X_{\sigma}}$. With such a definition of $\mathbf A_{\widehat X}$ we totally forget about the flags of the type $p\in\overline y\subset \widehat X$ where $y$ is horizontal and $p\in X_\infty$. So, we only add the flags of the type $p\in X_\sigma\subset \widehat X$ to the usual geometric picture. 
\end{remark}

Now we give the definitions of the completed spaces $A_{\widehat{\ast}}$: denote by $K_{\overline y}$  the diagonal embedding of $K_y$ inside $\prod_{x\in\overline y}K_{x,\overline y}$, then we put:
$$\overline{A}_{01}:=\overline{\mathbf A}_{X}\cap\prod_{\overline y\subset \widehat{X}} K_{\overline y}\,.$$
Moreover for any $\sigma$ let $A_{0}(\sigma)$ be the diagonal embedding  $\mathbb C(X_\sigma)\hookrightarrow\prod_{p\in X_\sigma} K_{p,\sigma}$, then:
$$A_{\widehat{01}}:=\overline{A}_{01}\oplus\prod_{\sigma\in B_\infty}A_{0}(\sigma)\,.$$
If $x\in X$ we have the natural embedding $K_{x}\hookrightarrow\prod_{\overline y\ni x} K_{x,\overline y}$; if $p\in X_\sigma$ then we consider the diagonal $\Delta_{p,\sigma}\subset K_{p,\overline y}\times K_{p,\sigma}$, where $\overline y$ is the unique horizontal curve containing $p$ (remember that  $K_{p,\overline y}= K_{p,\sigma}$). Thus we define:

$$\displaystyle A_{\widehat{02}}:=\mathbf A_{\widehat X}\cap \left(\prod_{x\in X} K_x \times \prod_{\substack{p\in X_\sigma, \\\sigma\in B_\infty}}\Delta_{p,\sigma}\right)\,.$$
$\overline{A}_0$ is the diagonal embedding of $K(X)$ in $\overline{\mathbf A}_{X}$ and: 
$$A_{\widehat 0}:=\overline{A}_0\oplus\prod_{\sigma\in B_\infty}A_{0}(\sigma)\,.$$
Note that: $A_{\widehat{01}}\cap A_{\widehat {02}}=A_{\widehat 0}$\,.
\begin{remark}
On a completed arithmetic surface we have a ``generalized version'' of the function field, it is not just $K(X)$ because we have fibres at infinity. It should be intended as $K(X)\oplus\prod_{\sigma}\mathbb C(X_\sigma)$ and note that this coincides with $A_{\widehat 0}$.
\end{remark}
The other adelic subspaces are the followings:

$$A_{\widehat {12}}:= \mathbf A_{\widehat X}\cap \left(\prod_{\substack{x\in \overline y,\\ \overline y\in \widehat X}} \mathcal O_{x,\overline y}\times\prod_{\substack{p\in X_\sigma,\\ \sigma\in B_\infty}} \mathcal O_{p,\sigma}\right)\,,\quad A_{\widehat{1}}:=A_{\widehat{01}}\cap A_{\widehat{12}}\,,\quad A_{\widehat{2}}:=A_{\widehat{02}}\cap A_{\widehat{12}}\,.$$
and the containment relations are the same as the geometric case:

$$
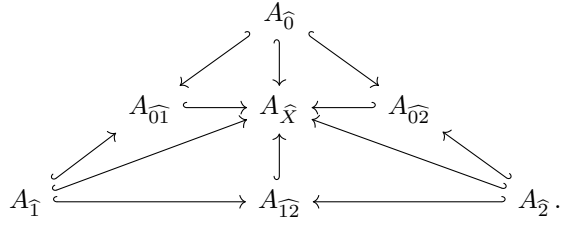
\begin{tikzcd}
&&A_{\widehat 0}\arrow[d,hook']\arrow[dl,hook']\arrow[dr,hook]&&\\
&A_{\widehat{01}}\arrow[r,hook]& A_{\widehat X} &A_{\widehat{02}}\arrow[l,hook']&\\
A_{\widehat 1}\arrow[rr,hook]\arrow[ur,hook]\arrow[urr,hook]&&A_{\widehat{12}}\arrow[u,hook]&&A_{\widehat 2}\arrow[ll,hook']\arrow[ul,hook']\arrow[ull,hook']\,.\\
\end{tikzcd}
$$

\section{Residue theory}\label{res_sect}
\subsection{Local multiplicative residues}

For any $b\in\widehat B$ we choose a standard character $\psi_b\colon K_b\to\mathbb T$ in the following way. By the classification theorem of locally compact local fields $K_b$ is a finite extension of $K_0$, where $K_0=\mathbb Q_p$, $K_0=\mathbb F_p((t))$, or $K_0=\mathbb R$, then we define $\psi^0_b\colon K_0\to\mathbb T$:
\begin{itemize}
\item[\sbt] If $K_0=\mathbb R$, put $ \psi_b^0(x):=e^{-2\pi ix}$.
\item[\sbt] If $K_0=\mathbb Q_p$ and $x=\sum_{j\ge m} a_j p^j\in\mathbb Q_p $, then put 
$$\psi_b^0(x):=e^{2\pi i\sum_{j=m}^{-1} a_j p^j}\,.$$
\item[\sbt] If $K_0=\mathbb F_p((t))$ and $x=\sum_{j\ge m} a_j t^j\in\mathbb F_p((t)) $, then for any lift $a'_{-1}\in\mathbb Z$ of $a_{-1}$, put 
$$\psi_b^0(x):=e^{2\pi i \frac {a'_{-1 }}{p}}\,.$$
\end{itemize}
Finally we define $\psi_b:=\psi_b^0\circ \tr_{K_b|K_0}$. Fix a completed curve $\overline y\subset \widehat X$, by considering all local branches in $y(x)$ we also define:
$$\Omega^1_{x, \overline y}:=\bigoplus_{\mathfrak z\in  y(x)}\Omega^1_{x,\mathfrak z}\,.$$
The structure of $\Omega^1_{x, \overline y}$  and the explicit expression of the $\res_{x,\mathfrak z}$ depend on the nature of $\overline y$:

\paragraph{$\overline y$ horizontal.} The local field $E_{x,\mathfrak z}$ is the constant field of  $K_{x,\mathfrak z}$ i.e. $K_{x,\mathfrak z}\cong E_{x,\mathfrak z}((t))$ and $[E_{x,\mathfrak z}: K_b]<\infty$. In \cite[2.2]{MM3} it is shown that there is an isomorphism  
\begin{equation}\label{iso_hori}
\Omega^{1}_{x,\mathfrak z}\cong E_{x,\mathfrak z}((t))dt
\end{equation}
where $t$ is a uniformizing parameter and moreover  the local residue assumes the following form independently from the choice of the isomorphism (\ref{iso_hori})
\begin{eqnarray*}
\res_{x,\mathfrak z}: \Omega^{1}_{x,\mathfrak z} &\to& K_b\\
adt &\mapsto& \tr_{E_{x,\mathfrak z}|K_b}(a_{-1})
\end{eqnarray*}
where $a=\sum_{j\ge m}a_jt^j\in E_{x,\mathfrak z}((t))$. Moreover we put: 

$$\res_{x, \overline y}:=\sum_{\mathfrak z\in  y(x)}\res_{x, \mathfrak z}:\Omega^1_{x,\overline y}\to K_b\,,$$

$$\Cres_{x, \overline y}:=\psi_b\circ\res_{x, y}:\Omega^1_{x,\overline y}\to \mathbb T$$
where $\psi_b:K_b\to\mathbb T$ is the standard character.

\paragraph{$\overline y=y$ vertical.}  $K_{x,\mathfrak z}$  is a finite extension of the standard field $L=K_p\{\!\{t\}\!\}$  where $[K_p: K_b]<\infty$ and $t$ is a uniformizing parameter for the residue field $\overline L=\overline K_p((t))$. Thanks to \cite[2.3]{MM3} we have an isomorphism 

\begin{equation}\label{iso_verti}
\Omega^{\cts}_{L|K_b}:=\left(\Omega^{1}_{\mathcal O_L|\mathcal O_b}\right)^{\sep}\otimes_{\mathcal O_L} L\cong K_p\{\!\{t\}\!\}dt
\end{equation}
 and a local residue independent from isomorphism (\ref{iso_verti}):
\begin{eqnarray*}
\res_{L}: \Omega^{\cts}_{L|K_b} &\to& K_b\\
adt &\mapsto& -\tr_{K_p|K_b}(a_{-1})
\end{eqnarray*}
where $a=\sum_{j\in\mathbb Z}a_jt^j\in K_p\{\!\{t\}\!\}$. By \cite[Remark 2.6]{MM3}, we know that $\Omega^1_{x,\mathfrak z}=\Omega^{\cts}_{L|K_b}\otimes_L K_{x,\mathfrak z}$, so we obtain a well defined trace map 
$$\tr_{K_{x,\mathfrak z}|L}\colon\Omega^1_{x,\mathfrak z}\to \Omega^{\cts}_{L|K_b} $$
At this point we define:

$$\res_{x,\mathfrak z}:=\res_{L}\circ\tr_{K_{x,\mathfrak z}|L}\colon\Omega^1_{x,\mathfrak z}\to K_b\,,$$

$$\res_{x, \overline y}:=\sum_{\mathfrak z\in  y(x)}\res_{x, \mathfrak z}:\Omega^1_{x,\overline y}\to K_b\,,$$
$$\Cres_{x, \overline y}:=\psi_b\circ\res_{x, y}:\Omega^1_{x,\overline y}\to \mathbb T\,,$$
where $\psi_b:K_b\to\mathbb T$ is the standard character.\\

When $\overline y$ is a completed horizontal curve and $x=p\in y_\sigma\subset y_\infty$ is a point at infinity,  then:
$$\Omega^1_{x,\overline{y}}:=\Omega^1_{p,\sigma}=K_{p,\sigma}dt;$$
$$\Cres_{x,\overline y}:=\psi_\sigma\circ \res_{p,\sigma}: \Omega^1_{p,\sigma}\to\mathbb T\,.$$
Where in the last line,  $\res_{p,\sigma}$ is the one dimensional residue on the Riemann surface $X_\sigma$ at the point $p$ and $\psi_\sigma:\mathbb C\to\mathbb T$ is the standard character of $\mathbb C$. 

Finally for a  flag at infinity $p\in X_\sigma$:
$$\Cres_{p,\sigma}:=\psi_\sigma\circ (-\res_{p,\sigma}): \Omega^1_{p,\sigma}\to\mathbb T\,.$$

\begin{remark}
The choice of the minus sign in the definition of $\Cres_{p,\sigma}$ is coherent with the main theory since $X_{\sigma}$ is vertical curve on $\widehat X$ in our geometric construction.
\end{remark}
The following proposition is the the extension of \cite[Lemma 3.3]{MM2} to the adelic case. It says that it makes sense to take the product of residues along vertical curves; moreover  by looking at its proof one immediately realizes that in the definition of two dimensional adeles, ``the first restricted product'' along a fixed curve is a crucial operation.

 \begin{proposition}\label{conv_vert}
 Let $\alpha\in \mathbf A_X$ and fix a vertical  curve $y\subseteq X_b$. Then the series
 $$\sum_{x\in y}\res_{x,y}(\alpha_{x,y}dt)$$
 converges in $K_b$. In particular $\res_{x,y}(\alpha_{x,y}dt)\in\mathcal O_b$ for all but finitely many $x\in y$.
  \end{proposition}
  \proof
 For simplicity let's assume that $y$ is nonsingular. We know that $(\alpha_{x,y})_{x\in y}\in \mathbb A^{(r)}_y$ for some $r\in\mathbb Z$, it means  that $(\alpha_{x,y})_{x\in y}=(t^r_y\beta_{x,y})_{x\in y}$ where  $(\beta_{x,y})_{x\in y}\in\mathbb A^{(0)}_y$. Now we use the definition of $\mathbb A_y^{(0)}$ to say that for any $s>0$ we have $\res_{x,y}(\beta_{x,y})\in\mathfrak p^{s+m}_{K_b}\mathcal O_b$ at almost all $x\in y$. It follows that for any $s>0$,  $\res_{x,y}(\alpha_{x,y})\in\mathfrak p^{s+m+r}_{K_b}\mathcal O_b$ at almost all $x\in y$. This shows that $\sum_{x\in y}\res_{x,y}(\alpha_{x,y}dt)$ converges in $K_b$.
  \endproof

By the universal property of the module of differential forms we have a canonical map $\Omega^{1}_{K(X)|K}\to \Omega^1_{x,\overline y}$, therefore by abuse of notation, we can consider an element  $\omega\in\Omega^{1}_{K(X)|K}$ as an element lying in $\Omega^1_{x,\overline y}$. Moreover, by base change we know that $\Omega^{1}_{\mathbb C(X_{\sigma})|\mathbb C}\cong\Omega^1_{K(X)|K}\otimes_{K(X)}\mathbb C(X_{\sigma})$, so we have a canonical composition map:
$$\Omega^1_{K(X)|K}\to \Omega^{1}_{\mathbb C(X_{\sigma})|\mathbb C}\to \Omega^1_{p,\sigma}$$
and  when clear from the context we can consider $\omega\in\Omega^{1}_{K(X)|K}$ as an element lying in $\Omega^1_{p,\sigma}$. In other words, it always makes sense to take a residue of a ``rational'' differential form $\omega\in\Omega^{1}_{K(X)|K}$ for flags in $X$ and in $\widehat X$.

\begin{theorem}[2D arithmetic reciprocity laws]\label{arith_rec_laws}
Let $\omega\in \Omega^1_{K(X)|K}$ then:
\begin{enumerate}
\item[$(1)$] Let  $x\in X$, then $\sum_{\overline y\ni x}\res_{x,y}(\omega)=0$ and $\res_{x,\overline y}(\omega)=0$ for all but finitely many curves $\overline y$ containing $x$. In particular  $\prod_{\overline y\ni x} \Cres_{x,\overline y}(\omega)=1$ and   $\Cres_{x,\overline y}(\omega)=1$ for all but finitely many $x\in y$.
\item[$(2)$] Let $p\in X_\sigma$, then $\Cres_{p,\sigma}(\omega)\cdot\prod_{\overline y\ni p} \Cres_{p,\overline y}(\omega)=1$.
\item[$(3)$] Let $\overline y\subset X$ be a vertical curve or $\overline y=X_\sigma$ for some $\sigma\in B_\infty$, then $\sum_{x\in \overline y}\res_{x,y}(\omega)=0$. In particular $\prod_{x\in \overline y} \Cres_{x,\overline y}(\omega)=1$ and   $\Cres_{x,\overline y}(\omega)=1$ for all but finitely many $x\in y$.
\item[$(4)$] Let $\overline y\in\widehat{X}$ be a horizontal curve, then $\prod_{x\in \overline y} \Cres_{x,\overline y}(\omega)=1$.
\end{enumerate}
\end{theorem}
\proof
See \cite[2.4]{MM2}, \cite[5]{MM2} and \cite[3]{MM2} for $(1)$, $(4)$ and the non-archimedean part of $(3)$  respectively. For the archimedean case of $(3)$ see \cite[Corollary of theorem 3]{tate_res}.\\
(2) For any fixed $p\in X_\sigma$ there is exactly one horizontal curve $\overline y\in\widehat X$ ``passing by'' $p$. Therefore
$$\Cres_{p,\sigma}(\omega)\cdot\prod_{\overline y\ni p} \Cres_{p,\overline y}(\omega)=
$$
$$=\Cres_{p,X_\sigma}(\omega)\Cres_{p,\overline y}(\omega)=\psi_\sigma(-\res_{p,\sigma}(\omega))\cdot \psi_\sigma(\res_{p,\sigma}(\omega))=1\,.$$
\endproof
\begin{remark}
Note that statements $(1)$ and $(2)$ of theorem \ref{arith_rec_laws} describe reciprocity laws around a point, whereas statements $(3)$ and $(4)$ describe reciprocity laws for a fixed curve. Archimedean data are taken in account without any special treatment: points on $X_\sigma$ are considered as points of $\widehat X$ and achimedean fibres are considered as vertical curves on $\widehat X$. We point out that statement $(2)$ is new and it hasn't been published anywhere before.
\end{remark}

\subsection{Adelic residue}
  In this subsection we globalise the local residues in order to get a residue at the level of completed adeles. Fix a nonzero rational $1$-form $\omega\in\Omega^1_{K(X)|K}$, then we define the adelic residue map:
  \begin{equation}\label{arith_ad_char}
\begin{tikzpicture}
\node at (4,4.05) {$\xi^\omega:$}; 
\node(1) at (4.6,4) {$\mathbf A_{\widehat X}$};
\node(2) at (13,4) {$\mathbb T$};
\node(3) at (4.7,3) {$(a_{x,\overline y})_{\substack {x\in \overline y,\\ \overline y\subset \widehat X}}\times (a_{p,\sigma})_{\substack {p\in X_{\sigma},\\ \sigma\in B_\infty}}$};
\node(4) at (13,3) {$\begin{displaystyle}\prod_{\substack {x\in \overline y,\\ \overline y\subset \widehat X}}\Cres_{x,\overline y}(\omega a_{x,\overline y})\prod_{\substack {p\in X_{\sigma},\\ \sigma\in B_\infty}}\Cres_{p,\sigma}(\omega a_{p,\sigma})\end{displaystyle}$};
\node at (4.6,3.6) {\begin{sideways}$\in$\end{sideways}};
\node at (13,3.7) {\begin{sideways}$\in$\end{sideways}};
\path[->] (1) edge (2);
\path[|->] (6.5,3.2) edge (10,3.2);
\end{tikzpicture}
  \end{equation}

Let's  explain why $\xi^\omega$ is well defined (i.e. the product (\ref{arith_ad_char}) is convergent): along all but finitely many curves $y\subset X$ the local residue is zero due to the restricted product with respect to the spaces $\mathbb A^{(0)}_y$. For the remaining curves we use the following arguments 
\begin{itemize}
\item[\sbt] If $y$ is horizontal it is enough to look at property $(\ast\ast)$ of lemma \ref{explemma}. It follows that the residue is $0$ at all but finitely many points of $y$.
\item[\sbt] If $y$ is vertical we use proposition \ref{conv_vert}.
\item[\sbt] For curves at infinity it is enough to appeal to the  $1$-dimensional adelic restricted product. 
\end{itemize}

\begin{figure}[htp!]
\centering
\begin{tikzpicture}

\draw (-3,2)node[anchor=east]{$X$}..controls (-1,2.5) and (1,1.5)..(3,2);
\draw (-3,2)--(-3,-2);
\draw (-3,-2)..controls (-1,-1.5) and (1,-2.5)..(3,-2);
\draw (3,2)--(3,-2);

\draw[name path=c2](-2.5,-1)..controls (0.5,-0.5) and (-.5,-1.5)..(2.5,-1);

\draw (-2.4,-1)node[anchor=south]{$y_1$};

\draw [name path=c1,](.6, -1.5)..controls (1.5,-0.5) and (2.5,-0.5).. (.6,1.5)node[anchor=west]{$y_2$};

\draw(4, -1.5)..controls (4.5,0) and (4.5,0).. (4,1.5)node[anchor=west]{$X_\sigma$}
node(a2)[midway] {$\bullet$}
node (a3)[midway,anchor=west]{$p$};

\fill [name intersections={
of=c1 and c2,  by={x}}]
(x) circle (2.2pt)node[anchor=north west]{$x$};

\draw (-3,-4)node[anchor= south east]{$B$}--(3,-4);
\draw[dotted] (3,-4)--(5.5,-4);
\draw (1,-4) circle [radius=2pt];
\fill (1,-4) circle [radius=2pt]node[anchor=north west]{$b$};

\draw [->] (1,-2.5)--(1,-3.5);
\draw [->] (4.3,-2.5)--(4.3,-3.5);

\draw (1,-3)node[anchor=west]{$\res_{x,y_1}+\res_{x,y_2}$};
\draw (4.3,-3)node[anchor=west]{$-\res_{p,\sigma}$};

\draw (4.3,-4) circle [radius=1.5pt];
\fill (4.3,-4) circle [radius=1.5pt]node[anchor=north west]{$\sigma$};

\draw [->] (6.3,-4)--(8,-4);
\draw (9.2,-4) circle [radius=20pt];
\draw (7,-4)node[anchor=south]{$\psi_b\cdot \psi_\sigma$};
\draw (9.2,-3.1)node[anchor=west]{$\mathbb T$};

\end{tikzpicture}
\caption{\footnotesize{A graphic representation of the action of the adelic residue on $3$ different flags: $x\in y_1$, $x\in y_2$ and $p\in X_\sigma$.}}
\end{figure}
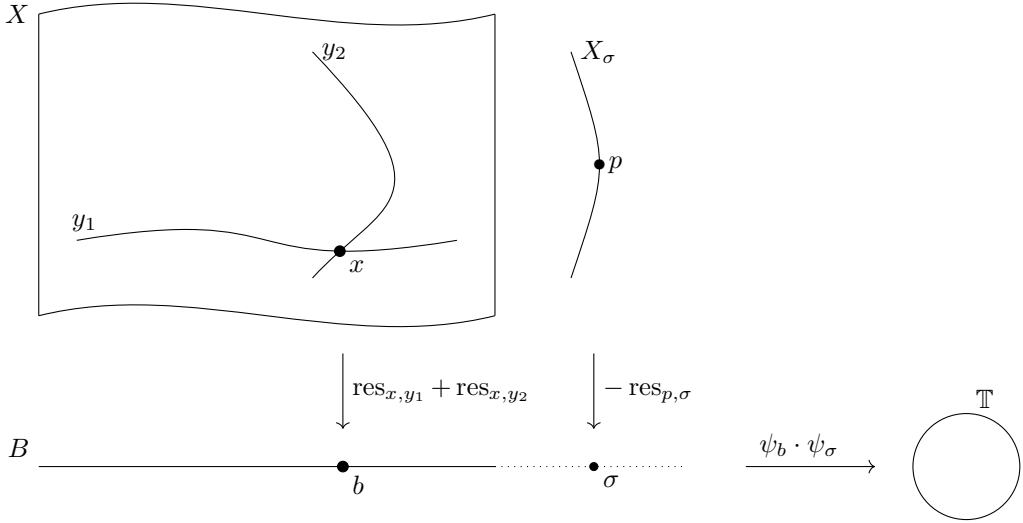

In \cite[Lemma 2.8, Remark 2.9]{MM2} it was proved that the local residues $\res_{x,y}:\Omega_{x,y}\to K_b$ are continuous, moreover it is clear that the local residues at infinity  $\res_{p,\sigma}:\Omega_{p,\sigma}\to \mathbb C$ are continuous (remember that here  $\mathbb C$ has the archimedean topology and $K_{p,\sigma}$ the $2$-dimensional topology).  We are interested in the global theory of residues and we will show that the adelic residue $\xi^{\omega}$ is sequentially continuous.

\begin{proposition}
The adelic residue $\xi^{\omega}$ is sequentially continuous.
\end{proposition}
\proof
To prove the sequential continuity of $\xi^{\omega}$ it is not necessary to consider the residues along curves at infinity because we have only a finite number of them and the  $1$-dimensional adelic residue is continuous. So, it is enough to discuss the schematic part of $\xi^{\omega}$ which will be denoted as  $\xi_S^{\omega}:\mathbf A_X\to\mathbb T$. Note that  we can write $\xi_S^\omega=\psi_S\circ \theta^{\omega}$ where $\psi_S:\mathbf A_{B}\to \mathbb T$ is the schematic part of the $1$-dimensional standard character and

$$\theta^{\omega}=(\theta^{\omega}_b)_{b\in B}:\mathbf A_X\to\mathbf A_B$$
with 

$$\theta^\omega_b:\prod_{\substack{x\in y,\\ y\subset X, \\ x\mapsto b}}K_{x,y}\to K_b$$
$$
\theta^{\omega}_b((\alpha_{x,y}))=\sum_{\substack {x\in X_b,\\ y\ni x}}\res_{x,y} (\omega\alpha_{x,y})=\overbrace{\sum_{\substack {y\subset X_b,\\ x\in y}}\res_{x,y} (\omega\alpha_{x,y})}^{(i)}+\overbrace{\sum_{\substack {x\in X_b,\\ y\ni x,\\ y \text{ horiz.}}}\res_{x,y}(\omega\alpha_{x,y})}^{(ii)}\in K_b\,. 
$$
 For any $n\in \mathbb N$, let   $\alpha^{(n)}:=(\alpha^{(n)}_{x,y})_{x,y}\in\mathbf A_X$ such that $\lim_{n\to\infty} \alpha^{(n)}_{x,y}=0$. Moreover put $\beta^{(n)}_{x,y}dt:=\omega\alpha^{(n)}_{x,y}$. Just for simplicity of notations we can assume that all curves involved in are nonsingular. We want to show that 
$$\lim_{n\to\infty}\xi^\omega_S\left(\alpha^{(n)}\right)=1\,.$$
Let's analyse carefully the summations $(i)$ and $(ii)$:\\

$(i)$ We know that $K_{x,y}$ is a finite extension of $L=K_p\{\!\{t\}\!\}$ and we can write 
$$\tr_{K_{x,y}|L}\left(\beta^{(n)}_{x,y}\right)=\sum_{i=-\infty}^{\infty}\beta^{(n)}_{x,y}(i)t^i\quad \text{for } \beta^{(n)}_{x,y}(i)\in K_p\,.$$
Then 
$$\res_{x,y}\left(\beta^{(n)}_{x,y}dt\right)= -\tr_{K_p|K_b}\left(\beta^{(n)}_{x,y}(-1)\right)\,.$$
Since $\lim_{n\to\infty} \beta^{(n)}_{x,y}=0$, there exists $n_0\in \mathbb N$ such that for $n\ge n_0$, we have $\beta^{(n)}_{x,y}(-1)\in\mathcal O_{K_p}$, i.e. $\res_{x,y}\left(\beta^{(n)}_{x,y}dt\right)\in\mathcal O_b$. This means that
$$\lim_{n\to\infty} \sum_{\substack {y\subset X_b,\\ x\in y}}\res_{x,y} \left(\beta^{(n)}_{x,y}dt\right)\in \mathcal O_b\,.$$

$(ii)$  
We know that $\beta^{(n)}_{x,y}=\sum_{i\ge m} \beta^{(n)}_{x,y}(i)t^i$, where $\beta^{(n)}_{x,y}(i)\in E_{x,y}$ and $E_{x,y}$ is a finite extension of $K_b$. Furthermore $\lim_{n\to\infty} \beta^{(n)}_{x,y}=0$.
We have:
\begin{equation}\label{thelim}
\lim_{n\to \infty}\sum_{\substack {x\in X_b,\\ y\ni x,\\ y \text{ horiz.}}}\res_{x,y}(\beta^{(n)}_{x,y}dt)=\lim_{n\to \infty}\sum_{\substack {x\in X_b,\\ y\ni x,\\ y \text{ horiz.}}}\tr_{E_{x,y}|K_b}(\beta^{(n)}_{x,y}(-1))\,.
\end{equation}
Due to the adelic restricted product, for all $n\ge n_0$ we have that $\res_{x,y}\left (\beta^{(n)}_{x,y}\right)=0$ along all but a fixed finite set of curves $y\subset X$, therefore we can exchange the summation and the limit in equation (\ref{thelim}). So we get:

$$
\lim_{n\to \infty}\sum_{\substack {x\in X_b,\\ y\ni x,\\ y \text{ horiz.}}}\res_{x,y}(\beta^{(n)}_{x,y}dt)=\sum_{\substack {x\in X_b,\\ y\ni x,\\ y \text{ horiz.}}}\lim_{n\to \infty}\tr_{E_{x,y}|K_b}(\beta^{(n)}_{x,y}(-1))=$$
$$
=\sum_{\substack {x\in X_b,\\ y\ni x,\\ y \text{ horiz.}}}\tr_{E_{x,y}|K_b}\left(\lim_{n\to \infty}\beta^{(n)}_{x,y}(-1)\right)=0\,.
$$
We can write:
$$\lim_{n\to\infty}\xi_S^\omega (\alpha^{(n)})=\lim_{n\to\infty}\psi_S(\theta^\omega(\alpha^{(n)}))\,.$$
For each $b\in B$ we have
$$\tau_b^{(n)}:=\theta^\omega_b(\alpha^{(n)})=\sum_{\substack {y\subset X_b,\\ x\in y}}\res_{x,y}(\beta^{(n)}_{x,y}dt)+\sum_{\substack {x\in X_b,\\ y\ni x,\\ y \text{ horiz.}}}\res_{x,y}(\beta^{(n)}_{x,y}dt)$$
and by $(i)$ and $(ii)$ we can conclude that:
$$\lim_{n\to\infty}\tau_b^{(n)}\in\mathcal O_b\,.$$
Finally:

$$\lim_{n\to\infty}\xi_S^\omega (\alpha^{(n)})=\lim_{n\to\infty}\psi_S\left(\left(\tau_b^{(n)}\right)_{b\in B}\right)=\psi_S\left(\left(\lim_{n\to\infty} \tau_b^{(n)}\right )_{b\in B}\right)=1\,.$$
\endproof

From the sequential continuity of the adelic residue we can deduce a stronger version of reciprocity laws:
\begin{proposition}\label{improved_rec}
Fix a rational differential form $\omega\in \Omega^1_{K(X)|K}$. Then the following statements hold:
\begin{enumerate}
\item[$(1)$] Fix a point $x\in X$. For any $\alpha\in K_{x}$ we have $\prod_{\overline y\ni x} \Cres_{x,\overline y} (\alpha\omega)=1$.
\item[$(2)$] Fix a  curve $\overline y\subset\widehat{X}$. For any $\alpha\in K_{\overline y}$ we have $\prod_{x\in\overline y} \Cres_{x,\overline y} (\alpha\omega)=1$.
\end{enumerate}
\end{proposition}
\proof
$(1)$ $K_x=K(X)\mathcal O_x$, but $\mathscr O_{X,x}$ is dense into its completion $\mathcal O_x$. Then the claim follows from the sequential continuity of the adelic residue and theorem \ref{arith_rec_laws}(1).

$(2)$ Again It follows from   the fact that $K(X)$ is dense in its completion (with respect to $y$)  $K_{y}$, the sequential continuity of the adelic residue,  and theorem \ref{arith_rec_laws}(3)-(4).
\endproof

\section{Self-duality of completed adeles}\label{self}
This section is entirely dedicated to the proof that the additive group $\mathbf A_{\widehat X}$ is self-dual. We will reduce the problem to show the self-duality of $\mathbb A_{y}$ and $\mathbf A_{X_{\sigma}}$.

The following two lemmas characterize the characters of $\mathbf A_{\widehat X}$:
\begin{lemma}\label{1stlemmachar}
Let $\chi\in \widehat{\mathbf A_{\widehat X}}$, then $\chi\left(\mathbb A^{(0)}_{\overline y}\right)=1$ for all but finitely many curves $\overline y\subset\widehat X$. In particular if $\beta:=\left(\beta_{\overline y}\right)_{\overline y}\times (\beta_{\sigma})_\sigma\in\mathbf A_{\widehat X}$ we have that 
$$\chi(\beta)=\prod_{\overline y\subset \widehat X}\chi(\beta_{\overline y})\prod_{\sigma\in B_\infty}\chi(\beta_\sigma)\,.$$
(In the above formula we clearly embedded each $\beta_{\overline y}$ and $\beta_\sigma$ naturally in $\mathbf A_{\widehat X}$).
\end{lemma}
\proof
Let $U\subset\mathbb T$ be an open neighborhood of $1$ which contains no subgroups of $\mathbb T$ other than $\{1\}$ and let $V\subset \mathbf A_{\widehat X}$ be an open subset such that $\chi(V)\subset U$. By the definition of restricted product as direct limit with the final topology, we know that for any finite set  $S$ of completed curves in $\widehat X$ the subset  $V\cap\mathbf A_{\widehat X}(S)$ is open in $\mathbf A_{\widehat X}(S)$. In particular by the definition of product topology, it contains an open subset of the following form:
$$W=\prod_{\overline y\notin S'} \mathbb A^{(0)}_{\overline y}\times\prod_{\overline y\in S'} W'_{\overline y}\times\prod_{\sigma\in B_\infty} W'_{\sigma}$$
where $S'$ is another finite set of completed curves in $\widehat X$ and $W'_{\overline y}\subset \mathbb A_{\overline y }$, $W'_{\sigma}\subset\mathbf A_{X_\sigma}$ are open. It follows that $H:=\chi\left(\prod_{\overline y\notin S'} \mathbb A^{(0)}_{\overline y}\right)\subset U$, but $H$ is a subgroup of $\mathbb T$, thus $H=\{1\}$ by the choice of $U$. In particular $\chi\left(\mathbb A^{(0)}_{\overline y}\right)=1$ for any $\overline y\notin S'$. The last assertion of the lemma is straightforward.
\endproof
\begin{lemma}\label{2ndlemmachar}
Let $\chi_{\overline y}\in\widehat{\mathbb A_{\overline y}}$ and let $\chi_{\sigma}\in\widehat{\mathbf A_{X_\sigma}}$. If $\chi_{\overline y}\left(\mathbb A^{(0)}_{\overline y}\right)=1$ for all but finitely many curves $\overline y\subset\widehat X$, then

 $$\chi:=\prod_{\overline y\subset \widehat X}\chi_{\overline y}\prod_{\sigma\in B_\infty}\chi_{\sigma}\in\widehat{\mathbf A_{\widehat X}}$$
\end{lemma}
\proof
The only thing that is not straightforward is the continuity of $\chi$, and there is no need to consider the fibres at infinity since they are finitely many. Let $U\subset\mathbb T$ be an open neighborhood of $1$ and choose $V\subset U$ such that\footnote{By $\MyProd$ we denote the actual complex multiplication of all elements in the open sets. In this particular case we are taking the ``$m$-th power of $V$''} $\MyProd_{m} V\subset U$. Now pick a finite set of completed curves  $S\subset\Upsilon$ of cardinality $m$, and for any $\overline y\in S$ take  $W_{\overline y}\subset \chi^{-1}_{\overline y}(V)$. Then $\prod_{\overline y\in S} W_{\overline y}\times \prod_{\overline y\not\in S}\mathbb A^{(0)}_{\overline y}$ is contained in the preimage of $\prod_{\overline y}\chi_{\overline y}$.
\endproof

The following proposition is basically the ``reduction argument'' that allows us to restrict our attention to $\mathbb A_y$ and $\mathbf A_{X_\sigma}$.

\begin{proposition}\label{red_arg}
The following isomorphism of topological groups holds:

$$\widehat{\mathbf A_{\widehat X}}\cong \sideset{}{'}\prod_{\overline y\subset \widehat X}\widehat{\mathbb A_{\overline y}}\times \!\!\!\prod_{\sigma\in B_\infty}\widehat{\mathbf A_{X_\sigma}}\,$$
where on the right hand side the restricted product is taken with respect to the subgroups $\left(\mathbb A^{(0)}_{\overline y}\right)^{\perp}\subset \widehat{\mathbb A_{\overline y}}$.

\end{proposition}

\proof
Consider the map:

\begin{eqnarray*}
\Psi:\sideset{}{'}\prod_{\overline y\subset \widehat X}\widehat{\mathbb A_{\overline y}}\times \!\!\!\prod_{\sigma\in B_\infty}\widehat{\mathbf A_{X_\sigma}}&\to&\widehat{\mathbf A_{\widehat X}}\\
(\chi_{\overline y})_{\overline y\subset \widehat X}\times(\chi_\sigma)_{\sigma\in B_\infty}&\mapsto&  \prod_{\overline y\subset \widehat X}\chi_{\overline y}\,\,\prod_{\sigma\in B_\infty}\chi_\sigma\\
\end{eqnarray*}
where clearly we naturally considered $\chi_{\overline y}, \chi_\sigma\in\widehat{\mathbf A_{\widehat X}} $. From lemmas \ref{1stlemmachar} and \ref{2ndlemmachar} it follows immediately that $\Psi$ is an isomorphism of groups, so we have to prove that it is continuous and open. Let $U$ be an open neighborhood of $1$ and consider the compact of $\mathbf A_{\widehat X}$:
$$C=\prod_{\overline y\in S} C_{\overline y}\times\prod_{\overline y\notin S} \mathbb A^{(0)}_{\overline y}\times\prod_{\sigma\in B_\infty} C_{\sigma}$$
where $C_\sigma, C_{\overline y}$ are compacts, and we assume that $S$ has cardinality $m$. Then $\mathcal W(C,U)$ is a basic open neighborhood of $\mathbf A_{\widehat X}$ around the identity character. Take now $V\subset U$ such that $\MyProd_{m+\#B_{\infty}} V\subset U$ and consider:
$$W=\prod_{\overline y\in S}\mathcal W(C_{\overline y},V)\times\prod_{\overline y\notin S}\left(\mathbb A^{(0)}_{\overline y}\right)^{\perp}\times\prod_{\sigma\in B_\infty} \mathcal W(C_{\sigma},V)\,. $$
Then clearly  $\Psi(W)\subseteq\mathcal W(C,U)$. The proof of openess is similar.
\endproof

\begin{remark}\label{restr_prod_rem}
So in order to show the self-duality of $\widehat{\mathbf A_{\widehat X}}$ we are reduced to show two things:
\begin{enumerate}
\item[\sbt] The self-duality of $\mathbf A_{X_\sigma}$.
\item[\sbt] There are topological and algebraic isomorphisms  $\theta_{\overline y}:\mathbb A_{\overline y}\to \widehat{\mathbb A_{\overline y}}$  mapping homeomorphically $\mathbb A^{(0)}_{\overline y}$ onto $\left(\mathbb A^{(0)}_{\overline y}\right)^{\perp}$ for all but finitely many completed curves.
\end{enumerate}
\end{remark}

 For the self-duality of $\mathbf A_{X_\sigma}$ we will use the following general results about Laurent power series over a self-dual group.
\begin{lemma}\label{standchar}
Let  $G$ be a ST ring and suppose that $(G,+)$ is endowed with a standard character. Then $G((t))$ has a standard character with conductor equal to $0$ (see appendix \ref{indpro} to see how $G((t))$ is topologised and for the definition of conductor). 
\end{lemma}
\proof
Let $\xi$ be a standard character of $G$. First of all let's find explicitly a nontrivial character  of $G((t))$ which has conductor equal to $0$. Consider:
\begin{eqnarray*}
\psi^0: G((t)) &\to& \mathbb T\\
\sum_{i\ge m}a_it^i&\mapsto& \xi(a_{-1})
\end{eqnarray*}
Let $\psi\in \widehat{G((t))}$, we want to show that that there exists a uniquely determined $\alpha\in G((t))$ such that $\psi=\psi^0_{\alpha}$. Assume that $c_{\psi}=i$, for any $b\in G$ the map $b\mapsto \psi(bt^{i-1})$ defines a character on $G$ that by hypothesis is equal to $\xi_{a_0}$ for a uniquely determined $a_0\in G$. So consider the character:
$$\psi^1(x):=\frac{\psi(x)}{\psi^0(xa_0t^{-i})}\quad \text{for } x\in G((t))\,,$$
it is easy to verify that $\psi^1\left(t^{i-1}G[[t]]\right)=1$. By iterating the above argument, for any $j\ge 1$ one finds a uniquely determined $a_j\in G$ such that
$$\psi^{j+1}(x):=\frac{\psi^j(x)}{\psi^0(xa_jt^{-i+j})}=\frac{\psi(x)}{\psi^0\left(x\sum_{h=0}^j a_h t^{-i+h} \right)}$$
is a character trivial on $t^{i-1-j}G[[t]]$. By taking the limit for $j\to\infty$ we obtain:
$$1=\lim_{j\to\infty} \psi^j(x)=\frac{\psi(x)}{\psi^0\left(x\sum_{h\ge 0} a_h t^{-i+h} \right)}\,.$$
So we put $\alpha:=\sum_{h\ge 0} a_h t^{-i+h}$ and it follows that $\psi(x)=\psi^0(x\alpha)$.\\
Now we show the continuity and the openess of the map $G((t))\to\widehat {G((t))} $ defined by $\alpha\mapsto \psi^{0}_\alpha$.  It is enough to prove the following simple things: 

\begin{itemize}
\item[$(a)$] Given a compact  $C\subset G((t))$ and an open $U\ni 1$ in $\mathbb T$,  there exist an open set $V\ni 0$ in $G((t))$ such that: $a\in V\Rightarrow aC\subseteq \psi^{-1}(U)$\,.

\item[$(b)$] Given an open $U\ni 0$ in $A$ there exist a compact $C\subset G((t))$ and an open $V\ni 1$ in $\mathbb T$,  such that: $ aC\subseteq \psi^{-1}(V)\Rightarrow a\in U$\,.
\end{itemize}

The explicit proofs of $(a)$ and $(b)$ are a respectively a very special case of the proofs of continuity and openess assertions of theorem \ref{bigselfdual}, so they are omitted here.

\endproof
\begin{proposition}\label{one_dim_dual}
The additive group $\mathbf A_{X_\sigma}$ is self-dual for every $\sigma\in B_\infty$
\end{proposition}
\proof
For any point $p\in X_\sigma$, we have $K_{p,\sigma}\cong\mathbb C((t))$, therefore we can apply lemma \ref{standchar} to conclude that $K_{p,\sigma}$ is self-dual and that a standard character with conductor equal to $0$ is given by $a\mapsto \Cres_{p,\sigma}(adt)$. At this point it is enough to follow line by line the argument in Tate's thesis that shows that adeles over a number field are self-dual (see for example \cite[5.1]{RamaVal}) to prove that $\mathbf A_{X_\sigma}$ is self-dual. Actually one needs the $1$-dimensional version of lemmas \ref{1stlemmachar}, \ref{2ndlemmachar} and \ref{red_arg}, but recall that we have a $2$-dimensional topological structure on $K_{p,\sigma}$ and $\mathbf A_{X_\sigma}$.
\endproof

 When $\overline y$ is horizontal one can apply lemma \ref{standchar} and the explicit expression of  $\mathbb A_{\overline y}$ to show that $\mathbb A_{\overline y}$ is self-dual, but when $\overline y$ is vertical, the proof is more problematic because we don't have any nice expression of  $\mathbb A_{\overline y}$ is terms of one dimensional adeles. 
A deeper analysis of the proof of lemma \ref{standchar} unravels that the only real advantage of having the expression $A=G((t))$, is  the ind-pro structure of $A$ over a self-dual group $G$ with a standard character. In general  also $\mathbb A_y$ has such property, and the following theorem is a generalisation of lemma \ref{standchar}, where $\mathbb A^{(-1)}_y/\mathbb A_y^{(0)}$ plays the role of $G$ and $\mathbb A_y^{(0)}$ plays the role of $G[[t]]$. We will heavily employ the topological properties described in subsection \ref{thetop}.

\begin{theorem}\label{bigselfdual}
 The additive group $\mathbb A_{\overline y}$  is self-dual with a standard character $\psi^0$. Moreover $\psi^0\in \left(\mathbb A^{(0)}_{\overline y}\right)^\perp$ and $\Theta_{\psi^0}\left(\mathbb A^{(0)}_{\overline y}\right)=\left(\mathbb A^{(0)}_{\overline y}\right)^\perp$.
\end{theorem}
\proof
It is enough to work with $\mathbb A_y$. For simplicity of notations let's put $A_r:= \mathbb A^{(r)}_y$, $A:=\mathbb A_y$ and $t:=t_y$. Let's summarize some properties (all categorical limits are in the category of linearly topologised groups):
\begin{itemize}
\item[(1)] $A_r$ is complete and $A_r=\varprojlim_{j\ge 1} A_r/A_{r+j}$
\item[(2)] $A_{r}/A_{r+1}$ is locally compact and self-dual with a standard character.
\item[(3)] $A_r/A_{r+j}$ is locally compact for every $j>0$.
\item[(4)] $A=\varinjlim_r A_r=\bigcup_r A_r$ and  $\bigcap_r A_r=\{0\}$.
\item[(5)] Any open neighborhood of $0$ in $A$ contains some $A_r$.

\end{itemize}
\noindent
Fix a standard character $\overline{\xi}\in\widehat{A_{-1}/{A_0}}$. Then consider the following commutative diagram of topological groups with exact short sequences:
$$
\begin{tikzcd}
&0            &0            & 0      &\\
 0\arrow{r}& A_{-1}/A_{0}\arrow{r}\arrow{u}& A_{-2}/A_{0}\arrow{r}\arrow{u}&A_{-2}/A_{-1}\arrow{r}\arrow{u}& 0\\
           &A_{-1}\arrow{u}        &A_{-2}\arrow{u}        & A_{-2}\arrow{u}      &\\ 
           &A_{0}\arrow{u}        &A_{0}\arrow{u}        & A_{-1}\arrow{u}      &\\
           &0\arrow{u}            &0\arrow{u}            & 0\arrow{u}      &\\
\end{tikzcd}
$$
Since the dual functor is exact on the category of LCA groups, we get the following diagram with exact short sequences:

$$
\begin{tikzcd}
&0 \arrow{d}           &0 \arrow{d}          & 0 \arrow{d}     &\\
 0& \widehat{A_{-1}/A_{0}}\arrow{d}\arrow{l}&\widehat{ A_{-2}/A_{0}}\arrow{d}\arrow{l}&\widehat{A_{-2}/A_{-1}}\arrow{l}\arrow{d}& 0\arrow{l}\\
&\widehat{A_{-1}}\arrow{d}        &\widehat{A_{-2}}\arrow{d}        & \widehat{A_{-2}}\arrow{d}      &\\ 
           &\widehat{A_{0}}       &\widehat{A_{0}}        & \widehat{A_{-1}}      &\\
\end{tikzcd}
$$

\noindent
In other words $\overline \xi$ lifts to a chacrater $\xi^1\in\widehat{A_{-1}}$ which is trivial on $A_0$, then we can   lift $\xi^1$ to a character $\xi ^2\in\widehat{A_{-2}}$ which extends $\xi^1$. By iterating this process we clearly construct a character  $\xi^n\in\widehat{A_{-n}}$ extending $\xi^1$. Now we can define a character $\psi:A \to\mathbb T$ in the following way:
$$\psi^0(a):=\xi^n(a)\quad \text{ if }\quad a\in A_{n}\setminus A_{n+1}\,.$$
By construction $\psi^0$ is trivial on $A_0$. A more explicit expression of $\psi^0$ can be given by using the identification $A_r=t^rA_0[[t]]$: if $a=\sum_{i\ge r}a_it^i\in A_r$, then $\psi^0(a)=\overline{\xi}(\overline{a_{-1}t^{-1}})$ where $\overline{a_{-1}t^{-1}}$ is the natural projection of $a_{-1}t^{-1}$ onto $A_{-1}/A_{0}$. We want to prove that $\psi^0$ is a standard character for $A$, so that the map:
\begin{eqnarray*}
\Theta_{\psi^0}:A &\to& \widehat{A}\\
a &\mapsto& \psi^0_a
\end{eqnarray*}
is an algebraic and topological isomorphism.

\emph{Surjectivity.} Since $\overline{\xi}$ is a standard character of $A_{-1}/A_{0}$, any other character in $\widehat{A^{-1}}$ which is trivial on $A_0$ is of the form $\xi^{1}(g\,\cdot)$ for $g\in A_{-1}$.  Consider any $\psi\in \widehat A$ and let $i=c_\psi$ the minimum integer $i\in \mathbb Z$ such that $\psi(A_i)=1$, note that this integer always exists thanks to $(5)$ and the fact that $\mathbb T$ has no small subgroups. Then:
 $$\psi_{|A_{i-1}}(\cdot\, t^{-i})=\xi^{1}(\cdot\, a_0t^{-i})\quad \text{for } a_0\in A_{-1}$$
 Let's define the following character
 $$\psi^{1}(\cdot)=\frac{\psi(\cdot)}{\psi^{0}(\cdot\,a_0t^{-i})}\,,$$
 then for any $t^{i-1}b\in A_{i-1}$ ($b\in A_0$):
 $$\psi^{1}_{|A_{i-1}}(t^{i-1}b)=\frac{\psi(t^{i-1}b)}{\psi^{0}(a_0t^{-i}t^{i-1}b)}=\frac{\xi^{1}(a_0t^{-1}b)}{\psi^{0}(a_0t^{-1}b)}=1\,.$$
 In other words $\psi^1$ is trivial on $A_{i-1}$. By iterating the above process for $j>1$, we find elements $a_h\in A_{-1}$ and characters:
 
$$\psi^j(\cdot)= \frac{\psi(\cdot)}{\psi^{0}(\cdot\, \sum^j_{h=0} a_ht^{-i+h})}$$ 
which are trivial on $A_{i-1-j}$. Now for $g\in A$  take the limit:
$$1=\lim_{j\to\infty}\psi^j(g)=\frac{\psi(g)}{\psi^{0}(a\sum_{g\ge 0} a_ht^{-i+h})}\,.$$
We conclude that $\psi(\cdot)=\psi^0(\cdot\, \alpha)$ for $\alpha:=\sum_{h\ge 0} a_ht^{-i+h}$. The partial sums defining $\alpha$ form a Cauchy sequence in $A_{-1}$, which is complete, so $\alpha$ is actually an element of $A_{-1}$.\\

\emph{Injectivity.} For every $a\in A\setminus 0$, there exists $r\in\mathbb Z$ such that $\ker\psi_a$ is trivial on $A_r$ but not on $A_{r-1}$. \\

\emph{Continuity.} 
We have to show that given a compact  $K\subset A$ and an open $U\ni 1$ in $\mathbb T$,  there exist an open set $V\ni 0$ in $A$ such that: $\psi(VK)\subseteq U$. Since $K$ is contained in some $A_m$, by simplicity we can ``shift'' $K$ thanks to the multiplication by $t^{m-1}$ and assume $K\subset A_{-1}$. Then $K=\varprojlim_j K_j$ with $K_j\in A_{-1}/A_j$. Now,  since $\overline \xi$ is a standard character for $A_{-1}/A_{0}$, it is not difficult to show by induction that the multiplication in $A$ and the character $\psi^0$ induce an algebraic and topological isomorphism $A_{-1}/A_j\cong \widehat{A_{-j}/A_1}$ for any $j\ge 0$. Thus we induce perfect pairing of LCA groups:
$$e_j:A_{-1}/A_j\times A_{-j}/A_1\to\mathbb T\,.$$ 
Consider the orthogonal complement $W_j=K_j^\perp:=\{a\in  A_{-j}/A_1\colon e_j(K_j,a)=1\}$, then $W_j$ is open in $A_{-j}/A_1$. Let $V_j\subset A_{-j}$ be the lift of $W_j$, it follows that the open set $V=\bigcup_j V_j$ is the open set we were looking for.\\

\emph{Openess.}  
We have to show that given an open $U\ni 0$ in $A$ there exist a compact $K\subset A$ and an open $V\ni 1$ in $\mathbb T$,  such that $ aK\subseteq \psi^{-1}(V)\Rightarrow a\in U$. The open set $U$ is contains  a basic open subgroup  $\sideset{}{'}\sum U_it^i$ where we assume that $U_i= A_0$ for $i\ge m$. Since $A_{-1}/A_{0}$ has a standard character, for any $i<m$ there exists a compact $C_i\subset A_{-1}/A_{0}$ and an open pen $V_i\ni 1$ in $\mathbb T$ such that:
$$\overline{\xi}(\overline a \overline{C_i})\subset V_i\Rightarrow \overline a\in\overline{U_it^{-1}}\subset A_{-1}/A_{0}\,.$$
Since $\mathbb T$ has no small subgroups, we can actually choose $\overline{V_i}$ in a way that 
$$\overline{\xi}(\overline a \overline{C_i})=1\Rightarrow \overline a\in\overline{U_it^{-1}}\,.$$
Now, since for any $r\ge 1$ we have a surjective homomorphisms of LCA groups $A_{-1}/A_r\to A_{-1}/A_0$,  we can lift $\overline {C_i}$ to $C_i^r\in A_{-1}/A_r$ which in turn gives $C_i=\varprojlim_{r} C_i^r$ compact in $A_{-1}$. We put $K_i=C_it\in A_0$. For $i\ge m$ we choose $K_i=0$, so we construct the compact set $K=\sum_i K_it^{-i}$ in $A$. It is easy to show that $K$ and a small enough $V\subset T$ containing $1$ satisfy the requirements nedeed to show openess.\\

Clearly $\Theta_{\psi^0}(A_0)\subseteq(A_0)^\perp$. Let $\psi^0_a\notin\left(A_0\right)^\perp $, then there exists $b\in  A_0$ such that $\psi^0(ab)\neq 1$, but this means that $a\notin A_0$ otherwise we would have $\psi^0(ab)=1$.
\endproof

\begin{corollary}
$\mathbf A_{\widehat X}$ is self-dual. 
\end{corollary}
\proof
The proof follows directly from propositions \ref{red_arg}, \ref{one_dim_dual} and theorem \ref{bigselfdual}. For more clarity, see also remark \ref{restr_prod_rem}.
\endproof

\section{Properties of the adelic differential pairing}\label{prop_diff}
Fix a nonzero rational differential form $\omega \in \Omega^1_{K(X)|X}$, then the \emph{adelic differential pairing (associated to $\omega$)} is defined as:

\begin{eqnarray*}
d_{\omega}\colon \mathbf A_{\widehat X}\times\mathbf A_{\widehat X}&\to&\mathbb T\\
(\alpha,\beta)&\mapsto& \xi^{\omega}(\alpha\beta)\,.
\end{eqnarray*}
For any subset $S\subseteq\mathbf A_{\widehat X}$ we define the \emph{orthogonal complement of $S$ with respect to $d_{\omega}$}:
\begin{equation}\label{orto}
S^\perp:=\{\beta\in\mathbf A_{\widehat X}\colon d_{\omega}(S,\beta)=1\}\,.
\end{equation}
The operator $\perp$ in this section shouldn't be confused with  the  one for topological groups and their duals.

\begin{proposition}
The map $d_{\omega}$ has the following properties:
\begin{enumerate}
\item[$(1)$]  It is symmetric and sequentially continuous.
\item[$(2)$] For any couple of subgroups $H_1,H_2\subseteq\mathbf A_{\widehat X}$ we have $H_1^\perp\cap H_2^\perp=(H_1+H_2)^\perp$.  

\end{enumerate}
\end{proposition}
\proof
$(1)$ Symmetry is obvious, sequential continuity follows easily from the fact that $\xi^{\omega}$ and the product are sequentially  continuous. 

$(2)$  If $h\in H_1^\perp\cap H_2^\perp $, then $d_\omega(h, H_1+H_2)=d(h, H_1)+d(h, H_2)=1$, so one inclusion is proved. Vice versa assume that $h\in (H_1+H_2)^\perp$, then $d(h, H_i)=0$ for $i=1,2$, so we have also the other inclusion.
\endproof
Now we show that the spaces $A_{\widehat{01}}$ and $A_{\widehat{02}}$ are equal to their orthogonal complements. Compare these results with the ``geometric counterpart'' in \cite{fe0}.
\begin{theorem}
$A_{\widehat{01}}^\perp= A_{\widehat{01}}$\,.
\end{theorem}
\proof
We show the equality by showing two inclusions.
\noindent
First we prove that $A_{\widehat{01}}\subseteq A_{\widehat{01}}^\perp$. It is essentially a consequence of our reciprocity laws for completed arithmetic surfaces. We have to show that for any $\alpha,\beta\in A_{\widehat{01}} $, $d_{\omega}(\alpha,\beta)=\xi^{\omega}(\alpha\beta)=1$. Let $a=\alpha\beta$, then 

$$\xi^{\omega}(a)=\prod_{\substack {x\in \overline y,\\ \overline y\subset \widehat X}}\Cres_{x,\overline y}(\omega a_{x,\overline y})\prod_{\substack {p\in X_{\sigma},\\ \sigma\in B_\infty}}\Cres_{p,\sigma}(\omega a_{p,\sigma})\,.$$
The first product is equal to $1$ thanks to proposition \ref{improved_rec}(2); the second product is $1$ thanks to the one dimensional reciprocity law.

Next we show the inclusion $A^\perp_{\widehat{01}}\subseteq A_{\widehat{01}}$. We take an element $a=(a_{\overline y})\times (a_\sigma) \in A^\perp_{\widehat{01}}$. We need to show that $a_{\overline y}\in K_{\overline y}$ and $a_\sigma\in \mathbb A_0(\sigma)$.  We consider $3$ cases.

\paragraph{Curves at infinity.} Pick any $g\in \mathbb A_0(\sigma)\subset A_{\widehat{01}}$, then since $a\in A^\perp_{\widehat{01}}$
$$d_{\omega}(a,g)=\prod_{p\in X_\sigma}\Cres_{p,\sigma}(a_{p,\sigma} g\omega)=1\,.$$
If $\psi$ is the standard character of $\mathbb C$, it follows that
\begin{equation}\label{therealperp}
\sum_{p\in X_\sigma}\res_{p,\sigma}(a_{p,\sigma} g\omega)\in\ker\psi=\frac{1}{2}\mathbb Z+\mathbb Ri,
\quad \forall g\in A_0(\sigma)\,. 
\end{equation}
By equation (\ref{therealperp}) for any $\lambda\in\mathbb R$ we have
$$
\sum_{p\in X_\sigma}\res_{p,\sigma}(a_{p,\sigma} \lambda\omega)=\lambda\sum_{p\in X_\sigma}\res_{p,\sigma}(a_{p,\sigma}\omega)\in \frac{1}{2}\mathbb Z+\mathbb Ri\,.
$$
It follows that $\sum_{p\in X_\sigma}\res_{p,\sigma}(a_{p,\sigma} \omega)=0$. We can replace $a_\sigma$ with $a_\sigma h$ for any $h\in  A_0(\sigma)$ to get $\sum_{p\in X_\sigma}\res_{p,\sigma}(a_{p,\sigma}h \omega)=0$. In other words $a_\sigma$ lies in the orthogonal complement of $A_{0}(\sigma)$ with respect to the pairing:
\begin{eqnarray*}
T_{\omega}:\mathbf A_{X_\sigma}\times\mathbf A_{X_\sigma}&\to& \mathbb C\\
\left((\alpha_{p,\sigma}), (\beta_{p,\sigma})\right) &\mapsto& \sum_{p\in X_\sigma}\res_{p,\sigma}(\alpha_{p,\sigma}\beta_{p,\sigma}\omega)
\end{eqnarray*}
But we know that $\mathbb A_0(\sigma)$ is equal to its orthogonal complement (with respect to $T_{\omega}$). Such a result was proved for number fields in \cite[Theorem 4.1.4]{tatethesis}, but  see for example \cite[Theorem 2.21]{dolce_phd} for the function field case. Therefore we conclude that $a_\sigma\in\mathbb \mathbb A_0(\sigma)$.

\paragraph{$\overline y$ horizontal.}  
We know that  $\mathbb A_{\overline y}=\mathbf A_{\overline y}((t_y))$, where  $t_y$ is a local parameter for $K_{\overline y}\subset A_{\widehat{01}} $  and therefore any
$a_{\overline y}$ has the following expression:
$$ \mathbb A_{\overline y}=\mathbf A_{\overline y}((t))\ni a_{\overline y}=\sum_{i\ge m} a_it^i_y\quad \text{with } a_i\in\mathbf A_{\overline y}\,.$$
We can also take $\omega=fdt_y$. 
\noindent
Then for any $r\in\mathbb Z$ and any $c\in k(y)$: 

\begin{equation}\label{for_every_r}
d_{\omega}(a_{\overline y},cf^{-1}t_y^r)=\prod_{x\in\overline y}\Cres_{x,\overline{y}}(a_{x,\overline y} t_y^rcdt_y)=\xi^{dt_y}(a_{\overline y} t_y^rc)=1\,.
\end{equation}
\noindent
Then $\xi^{dt_y}(a_{\overline y}t_y^rc)$ is a standard character of the one dimensional adeles $\mathbf A_{\overline y}$ calculated at $ca_{-r-1}$. Since equation (\ref{for_every_r}) holds for every $r\in\mathbb Z$,  and $k(y)$ is equal to $k(y)^\perp$ in $\mathbf A_{\overline y}$ (again \cite[Theorem 4.1.4]{tatethesis}), we can conclude that $a_i\in k(y)$ for every $i$. This means that $a_{\overline y}\in K_{\overline y}$.

\paragraph{$\overline y=y$ vertical.}
 Let $\overline t\in k(y)$ be a uniformizing parameter and  consider a lift  $t\in \mathscr O_{X,y}$. Put $\omega=fdt$ and let $L=K_p\{\!\{t\}\!\}$ a standard subfield of $K_{x,y}$. We know that there exists an integer $s$ such that $\tr_{K_p|K_b}(\mathfrak p^s_{K_p})\subseteq \mathfrak p_{K_b}$. Fix $r\in\mathbb Z$ such that $\tr_{K_{x,y}|L}(p^{r}a_{x,y})\in \mathfrak p^s_L$.  For any $m\in\mathbb Z$ we have $f^{-1}p^rt^m\in K_y\subset A_{\widehat{01}}$, so since  $a\in A^\perp_{\widehat{01}}$ we obtain that
$$d(a,f^{-1}p^rt^m)=\sum_{x\in y}\res_{x,y}(p^ra_{x,y}\cdot t^mdt)\in\mathcal O_b\,.$$
Since  $\tr_{K_{x,y}|L}(p^{r}a_{x,y})\in \mathfrak p^s_L$ and $\tr_{K_p|K_b}(\mathfrak p^s_{K_p})\subseteq \mathfrak p_{K_b}$:
$$\sum_{x\in y}\overline{\res_{x,y}(p^ra_{x,y}\cdot t^mdt)}=0$$
Now we apply \cite[Corollary 2.23]{MM3} to write
$$0=\overline{\sum_{x\in y}{\res_{x,y}(p^{r}a_{x,y}\cdot t^mdt)}}=\sum_{x\in y}\overline{\res_{x,y}(p^{r}a_{x,y}\cdot t^mdt)}=\sum_{x\in y}e_{x,y}\res^{(1)}_{x,y}(\overline{p^{r}a_{x,y}}\cdot \overline{t^m}d\overline{t})=$$
$$=\sum_{x\in y}\res^{(1)}_{x,y}(e_{x,y}\overline{p^{r}a_{x,y}}\cdot \overline{t}^md\overline {t})
$$
where:
\begin{itemize}
    \item [\sbt] $\res^{(1)}_{x,y}:E_{x,y}\to k(b)$ is the one dimensional residue on $E_{x,y}$.
    \item[\sbt] $e_{x,y}:=e(K_{x,y}|K_b)$ is the ramification degree.
\end{itemize}
The above  relation holds for any $m\in\mathbb Z$, and moreover we apply the same one-dimensional argument used in the case of the curves at infinity  to conclude  that $k(y)$ is equal to  $k(y)^\perp$ in $\mathbf A_y$. It follows that $(\overline{a_{x,y}})_{x\in y}\in k(y)$, therefore $a_y\in K_y$.
\endproof

Before proving that $A_{\widehat{02}}$ is self-orthogonal we need to study with more detail the structure of a neighborhood of a point $x\in X$ such that $\varphi(x)=b$. Let's denote with $\spec^1\mathcal O_x$ the set of prime ideals of height $1$ in $\mathcal O_x$, then a curve $y$ passing by $x$ corresponds to the set of local branches $y(x)\subset\spec^1 \mathcal O_x$. But there might be some elements $\mathfrak q\in\spec^1 \mathcal O_x $ which don't correspond to any curve passing by $x$, those are exactly those ideals:

$$\mathfrak T(x):=\{\mathfrak q\in\spec^1 \mathcal O_x\colon \mathfrak q\cap\mathscr O_{X,x}=(0)\}\,.$$
The elements of $\mathfrak T(x)$ are called \emph{transcendental curves (passing by $x$)}. 

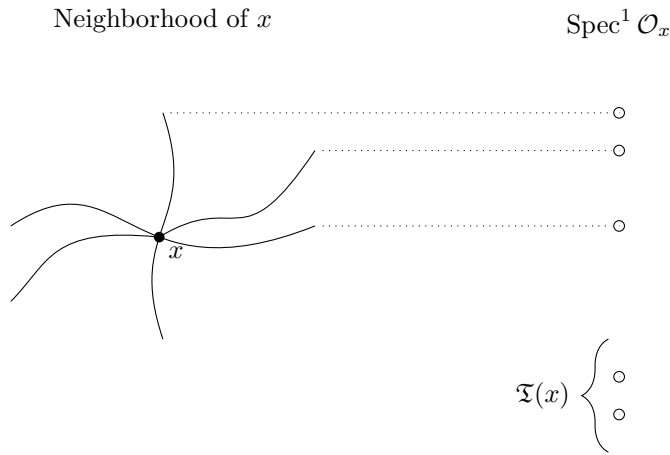
\begin{figure}[htp!]
\centering

\begin{tikzpicture}

\draw [name path=c1](2, -1.5)..controls (1.5,0) and (2.5,0).. (2,1.5);
\draw [name path=c2](0, 0)..controls (1.5,1) and (1.5,-1).. (4,0);
\fill [name intersections={
of=c1 and c2,  by={b}}]
(b) circle (2pt)node[anchor=north west]{$x$};

\draw[name path=c3] (0,-1)..controls (.5,-.5) and (.5,0)..(b);

\draw[name path=c4] (b)..controls (3,0.5) and (3,-.5)..(4,1);

\draw[dotted](2.1, 1.5)--(7.9,1.5);
\draw[dotted](4.1, 0)--(7.9,0);
\draw[dotted](4.1,1)--(7.9,1);

\draw (8,1.5) circle [radius=2pt];
\fill (8,1.5) circle; 

\draw (8,0) circle [radius=2pt];
\fill (8,0) circle; 

\draw (8,1) circle [radius=2pt];
\fill (8,1) circle; 

\draw (8,-2) circle [radius=2pt];
\fill (8,-2) circle;

\draw (8,-2.5) circle [radius=2pt];
\fill (8,-2.5) circle;

\draw (8,3) node[anchor=north]{$\spec^1 \mathcal O_x$};

\draw (2,3) node[anchor=north]{Neighborhood of $x$};

\draw [decorate,decoration={brace,amplitude=10pt, mirror,raise=4pt},yshift=0pt]
(8,-1.5) -- (8,-3) node [black,midway,xshift=-1cm] {$\mathfrak T(x)$};

\end{tikzpicture}

\caption{\footnotesize{A visual representation of the correspondence between prime ideals of $\mathcal O_x$ and curves passing by $x$. For simplicity we assumed that the curves are nonsingular at $x$, hence $y(x)$ is exactly a point in $\spec \mathcal O_x$ for any $y$.}}

\end{figure}

 Also for any $\mathfrak q\in \mathfrak T(x)$ it is possible to construct a $2$-dimensional local field $K_{x,\mathfrak q}$ and the residues $\res_{x,\mathfrak q}:\Omega_{x,\mathfrak q}\to K_b$,  $\Cres_{x,\mathfrak q}:\Omega_{x,\mathfrak q}\to \mathbb T$ in the usual way. But transcendental curves have the following pathological behaviour:
\begin{lemma}\label{transc}
Fix $\omega\in\Omega^1_{K(X)|K}$ and let $\mathfrak q\in\mathfrak T(x)$, then $\Cres_{x,\mathfrak q}(\omega)=1$. Moreover if $g\in K'_x$, then  $\Cres_{x,\mathfrak q}(g\omega)=1$.
\end{lemma}
\proof
The first claim follows immediately from the fact that $K(X)\subseteq (\mathcal O_x)_{\mathfrak q}$. For the second one it is enough to notice that $K_x=K(X)\mathcal O_x$ is dense in $K'_x$ and use the first part of the lemma.
\endproof

The presence of transcendental curves is a subtlety in the adelic theory. In fact, in general $K_x$ is obviously a proper subring of $K'_x$, but  the following result from commutative algebra ensures that  $K_x$ and $K'_x$ coincide if and only if there are no transcendental curves passing by $x$.

\begin{proposition}\label{equality_transc1}
Let $A$ be a Noetherian, regular, local domain  and let $\widehat A$ be the completion with respect to its maximal ideal. Then the product $\widehat A \fr{A}$ is a field if and only if for any nonzero prime $\mathfrak q\subset \widehat A$, $\mathfrak q \cap A\neq (0)$.
\end{proposition}
\proof
Since $A$ is regular and local, also $\widehat A$ is regular and local, which implies that $\widehat A$ and  $\widehat A \fr{A}$ are integral domains as well. It follows that $\widehat A\otimes_A \fr{A}\cong\widehat A\fr{A}$. Then it is well known (e.g. \cite[pag. 47]{matsu}) that we have an homeomorphism:
\begin{equation}\label{mats}
\spec(\widehat A\otimes_A \fr{A})\cong S:=\left\{\mathfrak q\in\spec\widehat A\colon \mathfrak q\cap A=(0)\right\}\,.
\end{equation}
$(\Rightarrow)$ $\widehat A\otimes_A \fr{A}$ contains only a prime ideal, the trivial one, so by the homeomorphism (\ref{mats}), $S$ contains only one element, which is $(0)$.
\noindent
$(\Leftarrow)$ If $S$ contains only $(0)$, then by the homeomorphism (\ref{mats}) the only prime ideal of   $\widehat A\otimes_A \fr{A}$ is $(0)$, which means that  $\widehat A\otimes_A \fr{A}$ is a field.
\endproof
\begin{corollary}\label{equality_transc2}
Fix a closed point $x\in X$, then $K_x=K'_x$ if and only if there are no transcendental curves passing by $x$.
\end{corollary}
\proof
By definition $K_x=K(X)\mathcal O_x$ is the smallest ring containing $K(X)$ and $\mathcal O_x$, so the claim follows from proposition \ref{equality_transc1}.
\endproof
Now let's put
$$\mathbf A_{X,x}:=\mathbf A_X\cap \prod_{y\ni x} K_{x,y}\,,$$
$$\mathbf A'_{X,x}:=\sideset{}{'}\prod_{\mathfrak q\in \spec^1\mathcal O_x} K_{x,\mathfrak q}\quad \text{with resp. to }\, \mathcal O_{x,\mathfrak q}\,,$$
and note that $\mathbf A'_{X,x}\supseteq \mathbf A_{X,x}$. Lemmas \ref{lemforA02} and \ref{red_tospecial_point} below will be used to show the inclusion $A_{\widehat{02}}^\perp\subseteq A_{\widehat{02}}$. The first one will be a  modified version of \cite[Lemma 3.3]{otherliu2}, so we present a proof. The second one will be just \cite[Lemma 3.4]{otherliu2} rewritten with our notation, so for its proof we remand the reader to the appropriate reference.

\begin{remark}
The paper \cite{otherliu2} shows only some local calculations regarding residues on the space $\mathbf A'_{X,x}$. Moreover the space denoted as $K_x$ in \cite{otherliu2} is exactly our $K'_x$.
\end{remark}


\begin{lemma}\label{fracpower}
Let $R$ be a ring, then 
$$\fr\left(R[[t]]\right)=F:=\left\{\sum_{i\ge m}a_it^i\in \fr(R)((t))\colon \exists r\in R \text{ such that } a_i\in R\left[1/r\right]\,,\forall i \right\}\,.$$
In particular, we deduce that in general $\fr\left(R[[t]]\right)$ is strictly contained in $\fr(R)((t))$.
\end{lemma}
\proof
Since $F$ is a field containing  $R[[t]]$, we have to show only the inclusion  $\fr\left(R[[t]]\right)\subseteq F$.  Let $\phi(t)=\frac {f(t)}{g(t)}\in\fr\left(R[[t]]\right) $  with $f(t), 0\neq g(t)\in R[[t]]$. Write $g(t)=t^k(r-t\gamma(t))=t^kr(1-\frac {t}{r}\gamma(t))$ with $k\geq 0$, $0\neq r\in R$  and $\gamma(t)\in R[[t]]$.  Then  $\frac {1}{g(t)}=t^{-k}\sum \frac {t^n}{r^n}(\gamma (t))^n$   and  $\phi(t)=\sum_{i=m }^{\infty } c_it^i$  where $m\in \mathbb Z$ depends on $\phi$ and each $c_i$ is of the form $c_i=\frac {\rho_i}{r^{\nu^i}}$ with $\rho_i\in R$ and $\nu_i\in \mathbb N$.  
 \endproof

\begin{lemma}\label{lemforA02}
Assume that $\mathcal O_x=\mathcal O_b[[t]]$. Fix a rational differential form $\omega\in\Omega^1_{K(X)|K}$ and let $a=(a_{x,\mathfrak q})\in\mathbf A'_{X,x}$ such that:
\begin{equation}\label{the_product}
\prod_{\mathfrak q\in\spec^1\mathcal O_x }\Cres_{x,\mathfrak q} (ga_{x,\mathfrak q}\omega)=1\quad\text{for any $g\in K'_x$\,,}
\end{equation}
then $a\in K'_x$.
\proof
There is a well known classification result for the elements  $\mathfrak q\in \spec^1 \mathcal O_b[[t]]=\spec^1 \mathcal O_x$ (see for example \cite[Lemma 5.3.7]{neucoh}):
\begin{itemize}
    \item[\sbt] $\mathfrak q=\mathfrak q_v:=\pi_b\mathcal O_x$, where $\pi_b$ is the uniformizing parameter of $\mathcal O_x$. This is the only prime ideal such that $K_{x,\mathfrak q_v}$ is of mixed characteristic.
    \item[\sbt] $\mathfrak q=h_{\mathfrak q}\mathcal O_x$ where $h_{\mathfrak q}\in\mathcal O_b[t]$ is an irreducible Weierstrass polynomial, i.e. $h_{\mathfrak q}=t^d+a_1t^{d-1}+\ldots+a_d$ with $a_i\in\mathfrak p_{K_b}$.
\end{itemize}
  Without loss of generality we may assume that $a_{x, \mathfrak q}\in \mathcal{O}_{x,\mathfrak q}$ for any $\mathfrak q\neq \mathfrak q_v$ since multiplying $(a_{x,\mathfrak q})$ by any nonzero element in  $K'_x$ amounts to an equivalent problem. Moreover, for the same reason  we can also assume for simplicity that $\omega=1dt$.

  For any $\mathfrak q\neq \mathfrak q_v$ and any uniformizing parameter $\pi_\mathfrak q$ for the $2$-dimensional local fields $K_{x,\mathfrak q}$, we can choose the following isomorphism:
  \begin{eqnarray*}
  K_{x,\mathfrak q} &\xrightarrow{\cong}& E_{x,\mathfrak q}((h_\mathfrak q))\\
  \pi_{\mathfrak q}&\mapsto& h_{\mathfrak q}(t)\,.
  \end{eqnarray*}
  In other words $t$ can be identified with a root of the polynomial equation  $h_{\mathfrak q}(t)=\pi_{\mathfrak q}$. By Hensel's lemma we deduce that such a root exists and it is integral, thus we can write:
  $$t=\sum_{i\ge0} c_i\pi^i_{\mathfrak q}\quad \text{with $c_i\in E_{x,\mathfrak q}$}\,.$$
 The following two easy results are fundamental:
  \begin{itemize}
      \item[$(i)$] $h_\mathfrak q\in\mathcal O^\times_{x,\mathfrak q'}$ for any $\mathfrak q'\neq \mathfrak q,\mathfrak q_v$. This is obvious from the definition of $\mathcal O^\times_{x,\mathfrak q'}$.
      \item[$(ii)$] $t\in\mathcal O_{x,\mathfrak q'}$ for any $\mathfrak q'\neq \mathfrak q $. Assume by contradiction that $t\notin\mathcal O_{x,\mathfrak q'}$  and let $h_{\mathfrak q}=t^d+a_1t^{d-1}+\ldots+a_d$, then  by $(i)$:
      $$0=v_{x,\mathfrak q'}(t^d+a_1t^{d-1}+\ldots+a_d)=\min\left\{v_{x,\mathfrak q'}(t^d), v_{x,\mathfrak q'}(a_1t^{d-1}),\ldots, v_{x,\mathfrak q'}(a_d)\right\}=$$
      $$=\min\left\{dv_{x,\mathfrak q'}(t), (d-1)v_{x,\mathfrak q'}(t),\ldots, 0\right\}=dv_{x,\mathfrak q'}(t)<0$$
      which cannot be true.
      \end{itemize}
 \noindent
If for any $\mathfrak q'\neq \mathfrak q_v$ we write
$$a_{x,\mathfrak q'}=\sum_{i\ge 0} a_{i,\mathfrak q'}\pi^i_{\mathfrak q'}\,, \quad  a_{i,\mathfrak q'}\in E_{x,\mathfrak q'}\,,$$
by $(i)-(ii)$ and equation (\ref{the_product}), for any $n\ge 0$ we have 
$$
\prod_{\mathfrak q'\in \spec^1\mathcal O_x} \Cres_{x,\mathfrak q'} (h_{\mathfrak q}^{-1}t^n a_{x,\mathfrak q'} \omega)=\Cres_{x,\mathfrak q} (h_{\mathfrak q}^{-1}t^n a_{x,\mathfrak q}\omega) \cdot
\Cres_{x,\mathfrak q_v} (h_{\mathfrak q}^{-1}t^n a_{x,\mathfrak q_v}\omega) =1\,.
$$
Therefore, we have the equality

\begin{equation}\label{eqforA02}
\Cres_{x,\mathfrak q} (h_{\mathfrak q}^{-1}t^n a_{x,\mathfrak q}\omega)=
\Cres_{x,\mathfrak q_v} (h_{\mathfrak q}^{-1}t^n a_{x,\mathfrak q_v}\omega)^{-1}\,,
\end{equation}
but by definition
\begin{equation}\label{eqforA02-}
\Cres_{x,\mathfrak q}(h_{\mathfrak q}^{-1}t^n a_{x,\mathfrak q}\omega)= \psi_b\left(\tr_{E_{x,\mathfrak q}|K_b}(c_0^n a_{0,\mathfrak q})\right)\,.
\end{equation}
Since we can take $1,c_0,...,c_0^{\deg h_{\mathfrak q}-1}$ as a basis of $E_{x,\mathfrak q}$ over $K_b$, equations (\ref{eqforA02}) and (\ref{eqforA02-}) imply that $\tr_{E_{x,\mathfrak q}|K_b}(\lambda a_{0,\mathfrak q})$ is determined by $a_{x,\mathfrak q_v}$ for any $\lambda \in E_{x,\mathfrak q}$. By using non-degeneracy of the trace pairing 
\begin{eqnarray*}
E_{x,\mathfrak q}\times E_{x,\mathfrak q}&\to &K_b\\
(u,s)&\mapsto&\tr_{E_{x,\mathfrak q}|K_b}(us)
\end{eqnarray*}
we conclude that the element $a_{0,\mathfrak q}$ is uniquely determined by $a_{x, \mathfrak q_v}$. We can conduct the same calculations  for $h_{\mathfrak q}^{-i-1}t^n a_{x,\mathfrak q}$, to  see that $a_{i,\mathfrak q}$ is determined by $a_{x,\mathfrak q_v}$ for any positive integer $i$.  It leads us to a conclusion that $a_{x,\mathfrak q}$ is uniquely determined by $a_{x,\mathfrak q_v}$ for any $\mathfrak q\neq \mathfrak q_v$.

So, we are reduced to show that  $a_{x,\mathfrak q_v}$ is in $K'_x$. Recall that $K_{x,\mathfrak q_v}\cong K_b\{\!\{t\}\!\}$, so we can write
$$
a_{x,\mathfrak q_v}=\sum_{i\in \mathbb{Z}}a_{i,\mathfrak q_v}t^i\,,\quad  a_{i,\mathfrak q_v}\in K_b\,.
$$
Now, by putting $\mathfrak p_0=t\mathcal O_x$ and reasoning similarly as above we obtain 
$$
 \Cres_{x,\mathfrak q_v}(t^{i-1}a_{x,\mathfrak q_v}\omega)^{-1}=\Cres_{x,\mathfrak p_0}(t^{i-1} a_{x,\mathfrak p_0}\omega)=1 \, , \quad \text{for all } i\geq 1\,.
$$
It means that $a^{-1}_{-i,\mathfrak q_v}\in\mathcal O_b$ for any $i\ge 1$. By definition of $K_b\{\!\{t\}\!\}$, we know that there exists $N>0$  such that $a_{-i,\mathfrak q_v}\in \mathcal O_b$ for $i\ge N$. In other words if $i\ge N$ and $a_{-i,\mathfrak q_v}\neq 0$, then  $a_{-i,\mathfrak q_v}\in \mathcal O_b^\times$. Since $\lim_{j\to-\infty}a_{j,\mathfrak q_v}=0$, we conclude that it has to exists $M>0$ such that $a_{-i,\mathfrak q_v}=0$ for $i\ge M$. This proves that  $a_{x,\mathfrak q_v}\in K_b((t))$.  Again thanks to the definition of $K_b\{\!\{t\}\!\}$, we know that there exists $m\in\mathbb Z$ such that $v_b(a_{i,\mathfrak q_v})\ge m$, which means that for any choice of  uniformizing parameter $s\in\mathcal O_b$, then  $a_{i,\mathfrak q_v}=s^{m+j} g$ with $j\ge 0$ and $g\in\mathcal O^\times_b$. We distinguish two cases:
\begin{itemize}
    \item[\sbt] If $m<0$, then $a_{i,\mathfrak q_v}=(1/s)^{-m}\cdot s^jg\in\mathcal O_b[1/s]$
    \item[\sbt] If $m\ge0$, then $a_{i,\mathfrak q_v}=(1/s)\cdot s^{m+j+1}g\in\mathcal O_b[1/s]$
\end{itemize}
Thanks to lemma \ref{fracpower} we conclude that $a_{x,\mathfrak q_v}\in K_x'$.
\endproof
\end{lemma}

Let $x\in X$ such that $\varphi(x)=b$; in general $\mathcal O_x$ is a finite ring extension of  $\mathcal O^\#_x:=\mathcal O_b[[t]]$. For any prime $u\in\spec^1\mathcal O^\#_x$ we have the  $2$-dimensional local field $K^\#_{x,u}$ obtained by the usual process of completion/localization. In general we can construct all local adelic objects relative to the flags $x\in u\in\spec^1\mathcal O^\#_x$. Such objects arising from the special ring  $\mathcal O^\#_x$ will be marked with the symbol $\#$ to distinguish them from the usual ones. Let $\mathfrak q\in\spec^1\mathcal O_x$ be a prime sitting over $u$, then we have a finite field extension $K_{x,\mathfrak q}|K^\#_{x,u}$ and a trace map $\tr_{K_{x,\mathfrak q}|K^\#_{x,u}}$ which extends directly at the level of differential forms:

\begin{eqnarray*}\label{trace_diff}
\tr_{K_{x,\mathfrak q}|K^\#_{x,u}}:\Omega^1_{x,\mathfrak q}& \to &\Omega^{1,\#}_{x,u}\\
fdt&\mapsto&  \tr_{K_{x,\mathfrak q}|K^\#_{x,u}}(f)dt
\end{eqnarray*} 
Such a map is exactly  the abstract trace map for differential forms defined in \cite{MM3} and mentioned in section \ref{res_sect}. We recall that in \cite{MM3} it is also proved that the residue is functorial with respect to the trace, which in our case means that $\res_{x,\mathfrak q}=\res^\#_{x,u}\circ\tr_{K_{x,\mathfrak q}|K^\#_{x,u}}\,.$
The local trace map defined above can be further generalized to an adelic trace:

\begin{eqnarray*}
\tr_x\colon\,\, \mathbf A'_{X,x} &\to& \left(\mathbf A'_{X,x}\right)^\#\\
(a_{x,\mathfrak q})_{\mathfrak q} &\mapsto &\left(\sum_{\mathfrak q|u}\tr_{K_{x,\mathfrak q}|K^\#_{x,u}}(a_{x,\mathfrak q})\right)_{u}
\end{eqnarray*}
where with the notation $\mathfrak q|u$ we denote all ideals $\mathfrak q\in\spec^1\mathcal O_x$ sitting over $u$.
\begin{lemma}\label{red_tospecial_point}
Let $f\in \mathbf A'_{X,x}$ such that $\tr_x(fg)\in (K'_x)^\#$ for any $g\in K'_x$, then $f\in K'_x$.
\end{lemma}
\proof
See \cite[Lemma 3.4]{otherliu2}.
\endproof


\begin{theorem}\label{self_ort_A02}
$A_{\widehat{02}}^\perp= A_{\widehat{02}}$\,.
\end{theorem}
\proof
First of all let's prove that $A_{\widehat{02}}\subseteq A_{\widehat{02}}^\perp$. We have to show that for any $\alpha,\beta\in A_{\widehat{02}} $, $d_{\omega}(\alpha,\beta)=\xi^{\omega}(\alpha\beta)=1$. Let $a=\alpha\beta$, then 

$$\xi^{\omega}(a)=\prod_{\substack {x\in \overline y,\\ \overline y\subset \widehat X}}\Cres_{x,\overline y}(\omega a_{x,\overline y})\prod_{\substack {p\in X_{\sigma},\\ \sigma\in B_\infty}}\Cres_{p,\sigma}(\omega a_{p,\sigma})=\prod_{\substack {x\in X,\\ \overline y\ni x}}\Cres_{x,\overline y}(\omega a_{x,\overline y})\prod_{\substack {p\in X_{\sigma},\\ \overline y\ni p,\\ \sigma\in B_\infty}}\Cres_{p,\sigma}(\omega a_{p,\sigma})\Cres_{p,\overline y}(\omega a_{p,\overline y})\,.$$
We can conclude $\xi^{\omega}(a)=1$ thanks to proposition \ref{improved_rec}(1) and from the explicit definition of $A_{\widehat{02}}$ at infinity.

Now we show the inclusion $A_{\widehat{02}}^\perp\subseteq A_{\widehat{02}}$. Fix $a=(a_{x,\overline y})\times(a_{p,\sigma})\in A_{\widehat{02}}^\perp $ and assume $\omega=fdt$, we consider two cases:

\paragraph{$x=p$ is a point on $X_\sigma$.} For any $g\in\mathbb C((t))$ we consider the element $(f^{-1}g,f^{-1}g)\in \Delta_{p,\sigma}$, then if $\overline y$ is the unique horizontal curve containing $p$ we obtain
$$\Cres_{p,\overline y}(a_{p,\overline y}f^{-1}g\omega)\cdot \Cres_{p,\sigma}(a_{p,\sigma}f^{-1}g\omega)=1\,.$$
This means
\begin{equation}\label{forA02}
\res_{p,\overline y}(a_{p,\overline y}gdt)- \res_{p,\sigma}(a_{p,\sigma}gdt)\in \ker\psi_\sigma= \frac{1}{2}\mathbb Z+\mathbb Ri\,.
\end{equation}
Since equation (\ref{forA02}) holds for any $g\in\mathbb C((t))$, it is clear that it must be $(a_{p,\overline y},a_{p,\sigma})\in \Delta_{p,\sigma}$.

\paragraph{$x$ is a point on $X$.} Recall that $\mathcal O_x$   is a finite ring extension of $\mathcal O_b[[t]]$. 

We first treat the case where there are transcendental curves passing by $x$; let's extend the element $(a_{x,\overline y})_{\overline y\ni x}$ to an element $(a'_{x,\mathfrak q})_{\mathfrak q}\in\mathbf A'_{X,x}$ in the following way: for a transcendental curve $\mathfrak q\in\mathfrak T(x)$ let's insert $a'_{x,\mathfrak q}\in K_x$; at all other primes nothing changes. Now let $g\in K_x$, then:

\begin{equation}\label{verybigprod}
\prod_{\mathfrak q\in\spec^1 \mathcal O_x}\Cres_{x,\mathfrak q}(a'_{x,\mathfrak q}g\omega)=\underbrace{\prod_{\overline y\ni x} \Cres_{x,\overline y}(a_{x,\overline y}g\omega)}_{(i)} \underbrace{\prod_{\mathfrak q\in\mathfrak T(x)}\Cres_{x,\mathfrak q}(a'_{x,\mathfrak q}g\omega)}_{(ii)}=1 
\end{equation}
where $(i)=1$ because  $(a_{x,\overline y})\times(a_{p,\sigma})\in A_{\widehat{02}}^\perp $ and $(ii)=1$ thanks to lemma \ref{transc}. Since $K_x$ is dense in $K'_x$, equation (\ref{verybigprod}) implies that for any $h\in K'_x$
\begin{equation}\label{verybigprod1}
\prod_{\mathfrak q\in\spec^1 \mathcal O_x}\Cres_{x,\mathfrak q}(a'_{x,\mathfrak q}h\omega)=1\,.
\end{equation}
Now we use equation (\ref{verybigprod1}) and the functoriality of the residue with respect to the trace map:

$$
\mathcal O_b\ni\sum_{\mathfrak q} \res_{x,\mathfrak q}(a'_{x,\mathfrak q}h\omega)=\sum_{u}\sum_{\mathfrak q|u}\res^\#_{x,u}\left(\tr_{K_{x,\mathfrak q}|K^\#_{x,u}}(a'_{x,\mathfrak q}h\omega) \right)=
$$
$$
=\sum_{u}\res^\#_{x,u}\left(\sum_{\mathfrak q|u}\tr_{K_{x,\mathfrak q}|K^\#_{x,u}}(a'_{x,\mathfrak q}h\omega) \right)= \sum_{u}\res^\#_{x,u}\left(\sum_{\mathfrak q|u}\tr_{K_{x,\mathfrak q}|K^\#_{x,u}}(a'_{x,\mathfrak q})h\omega \right)\,.
$$
By lemma \ref{lemforA02} we can conclude that $\tr_x\left(a'_{x,\mathfrak q}\right)_{\mathfrak q}\in (K')^\#_x$ diagonally.  By replacing $a_{x,\overline y}$ with $ca'_{x,\overline y}$ for any $c\in K'_x$ we can again conclude that  $\tr_x\left(ca'_{x,\mathfrak q}\right)_{\mathfrak q}\in (K')^\#_x$ diagonally.  At this point we can use lemma \ref{red_tospecial_point}  to conclude that $(a'_{x,\mathfrak q})_{\mathfrak q}\in K'_x$. It means that $(a_{x,\overline y})_{\overline y\ni x}\in K_x$ by the choice of  $(a'_{x,\mathfrak q})_{\mathfrak q}$.

If there are no transcendental curves passing by $x$, then $\mathbf A'_{X,x}=\mathbf A_{X,x}$ and $K_x=K'_x$ by corollary \ref{equality_transc2}. Then we can apply  a simplified  version of the argument used above to conclude the proof.
\endproof


\begin{remark} 
We were informed by I. Fesenko that there is an alternative proof of theorem \ref{self_ort_A02} which uses an arithmetic version of his argument in \cite{fe0}.
\end{remark}

\section{Idelic interpretation of Arakelov intersection theory}\label{subs_idelic_arakelov}
A prerequisite for this section is the whole appendix \ref{ar_g}. In \cite{adI}, it is described how to get a lift of the Deligne pairing (i.e. the schematic part of the Arakelov intersection number) at the level of ideles. Let's summarise the result: first of all we consider the idelic complex attached to the (uncompleted) surface $X$
\begin{equation}\label{id_coomplex}
\begin{tikzcd}
\mathcal A^\times_X: & A^\times_0\oplus A^\times_1\oplus A^\times_2 \arrow[r, "d_\times^0"]& A^\times_{01}\oplus A^\times_{02}\oplus A^\times_{12}\arrow[r, "d_\times^1"]& A^\times_{012}\\
& (a_0,a_1,a_2) \arrow[r, mapsto] & (a_0a^{-1}_1,a_2a^{-1}_0,a_1a^{-1}_2) &\\
&& (a_{01},a_{02},a_{12})\arrow[r, mapsto]& a_{01}a_{02}a_{12}\\
\end{tikzcd}
\end{equation}
and we note that we have a surjective map:
\begin{eqnarray*}
p:\;\; \ker(d^1_\times) &\to& \Div(X)\\
(\alpha,\beta,\alpha^{-1}\beta^{-1})&\mapsto & \sum_{y\subset X} v_y(\alpha_{x,y})[y]\,.
\end{eqnarray*}
Then by globalising the Kato's symbol, we define ad idelic Deligne pairing $\left<\,,\,\right>_i:\ker(d^1_\times)\times \ker(d^1_\times)\to\Pic(B)$ which descends to the Deligne pairing  $\left<\,,\,\right>:\Pic(X)\times\Pic(X)\to\Pic(B)$. In turn, the Deligne  pairing is strictly related  to intersection theory because for any two divisors $D,E\in \Div(X)$, the class in $\Pic(B)$ of the divisor 
$$\left<D,E\right>=\varphi_\ast i(D,E)=\sum_{x\in X} [k(x):k(\varphi(x))]\,i_{x}(D,E)\, [\varphi(x)]$$
is equal to $\left<\mathscr O_X(D), \mathscr O_X(D)\right>$. Note that we have used the brackets $\left<\,,\,\right>$ to denote two different (but strictly related) objects, but the clash of notations shouldn't confuse the reader.

The contribution at infinity to the Arakelov intersection pairing is given by the $\ast$-product between Green functions, so the next step in our theory is to find an idelic description of it. The part at infinity of the full adelic ring $\mathbf A_{X}\oplus\prod_{\sigma\in B_\sigma}(\mathbf A_{X_\sigma}\oplus \mathbf A_{X_\sigma})$ is given by $\mathbf A_{X_\sigma}\oplus \mathbf A_{X_\sigma}$ (for each $\sigma$), so we want to find a surjective map:
$$(\mathbf A^\times_{X_\sigma}\oplus \mathbf A^\times_{X_\sigma})\supseteq S\to \mathbb Z G(X_\sigma)$$
where $S$ is an adequate subset of $\mathbf A^\times_{X_\sigma}\oplus \mathbf A^\times_{X_\sigma}$ still to be determined and $\mathbb ZG(X_\sigma)$ is the vector space of Green functions on $X_\sigma$ with integer orders.
 \begin{remark}
First of all let's introduce a notation. For any $a=(a_x)\in\mathbf A_{X_\sigma}$,  with $a(x)$ we denote the projection of $a_x$ onto the residue field $\mathbb C$ (when it is well defined).
\end{remark}
Let $\mathcal F(X_\sigma,\mathbb R)'$ be the set of functions $f:U\subseteq X_\sigma\to\mathbb R$ whose domain $U$ is the whole $X_\sigma$ minus a finite set of points, then we have the following map:
\begin{eqnarray*}
\Theta:\mathbf A_{X_\sigma}^\times \times \mathbf A_{X_\sigma}^\times &\to& \mathcal F(X_\sigma,\mathbb R)'\\
(a,b) &\mapsto& -\log(ba\overline a):=[x\mapsto -\log\left(b(x)a(x)\overline{a(x)}\right)]
\end{eqnarray*}
where $\overline{a(x)}$ denotes the complex conjugate. Note that $\mathbb ZG(X_\sigma)\subset \mathcal F(X_\sigma,\mathbb R)'$, then put 

$$G(\mathbf A^\times_{X_\sigma}):=\{(a,b)\in\Theta^{-1}(\mathbb ZG(X_\sigma))\colon v_x(a_x)=\ord_x^G(\Theta(a,b)),\; \forall x\in X_\sigma\}\,.$$
We get the map:
$$\pi_\sigma:=\Theta|_{G(\mathbf A^\times_{X\sigma})}:G(\mathbf A^\times_{X\sigma})\to \mathbb ZG(X_\sigma)\,.$$
\begin{proposition}\label{idgrsurj}
The map $\pi_\sigma$ is surjective.
\end{proposition}
\proof
Let $g\in \mathbb ZG(X_\sigma)$, by proposition \ref{Zgreen}, there exist a $C^\infty$ hermitian invertible sheaf $(\mathscr L, h)$ on $X$ and a meromorphic section $s=\{(s_j, U_j)\}$ of $\mathscr L$ such that we can write:
$$g=-\log(h(s,s))\,.$$
We can choose $a\in\mathbf A^\times_{X_\sigma}$ such that $a(x)=s(x)$ (when $s(x)$ is well defined) and $v_x(a_x)=\ord_x(s)$ for any $x\in X_\sigma$. Now we can write
$$g(x)=-\log(h_{x}(a(x),a(x)))\,.$$
Since $z\mapsto h_{x}(z\overline z)$ is a complex absolute value, we have $h_{x}(z\overline z)=w_{x}z\overline z$ with $w_{x}\in\mathbb  C$. Let's choose $b=(b_{x})\in\mathbf A^\times_{X_\sigma}$ such that  $b(x)=w_{x}$, then 
$$g(x)=-\log\left(b(x)a(x)\overline{a(x)}\right)\,.$$
The fact that $v_x(a_x)=\ord_x^G(g)$ follows directly from the fact that for any hermitian metric $h$ and meromorphic section $s$ we have the equality:
$$\divi^G(-\log(h(s,s)))=\divi(s)\,.$$
(See proposition \ref{gf1}).
\endproof

So far, we have the idelic description of Green functions with integer orders thanks to the projection $\pi_\sigma$. Now let's fix a (normalised) K\"ahler fundamental form $\Omega_\sigma$ on $X_\sigma$ and  consider $G_0^{\Omega_\sigma}(\mathbf A^\times_{X_\sigma}):=\pi^{-1}_\sigma(\mathbb ZG^{\Omega_\sigma}_0(X_\sigma))$, $G^{\Omega_\sigma}(\mathbf A^\times_{X_\sigma}):=\pi^{-1}_\sigma(\mathbb ZG^{\Omega_\sigma}(X_\sigma))$. For pairs $(\alpha,\beta)\in G^{\Omega_\sigma}(\mathbf A^\times_{X_\sigma})\times  G^{\Omega_\sigma}(\mathbf A^\times_{X_\sigma})$ such that $\divi^G((\pi_\sigma(\alpha))$ and  $\divi^G(\pi_\sigma(\beta))$ have no common components we want to find a product $\alpha\ast_i\beta$ such that the following equality holds:
$$
\begin{tikzcd}
(\alpha,\beta)\arrow[d, mapsto]\arrow[dr,mapsto]& \\
(\pi_\sigma(\alpha),\pi_\sigma(\beta))\arrow[r, mapsto] & \alpha\ast_i\beta=\pi_\sigma(\alpha)\ast\pi_\sigma(\beta)
\end{tikzcd}
$$
As a consequence of the symmetry of the $\ast$-product we will get also the symmetry of $\ast_i$. For any $\alpha=(a,b)\in G^{\Omega_\sigma}(\mathbf A^\times_{X_\sigma})$ let's put:
$$\xi(\alpha):=e^{\int_{X_\sigma}\log(ba\overline a)\Omega_\sigma}$$
\begin{definition}
Let $\alpha=(a,b), \beta=(c,d)\in G^{\Omega_\sigma}(\mathbf A^\times_{X_\sigma}) $, then the idelic $\ast$-product is defined as:
$$\alpha\ast_i\beta:=-\sum_{x\in X_\sigma}v_x(c_x)\log\left(b(x)a(x)\overline{a(x)}\xi(\alpha)\right)+\log(\xi(\alpha))\ideg(c)+\log(\xi(\beta))\ideg(a)\,.$$
\end{definition}

\begin{proposition}
$(\alpha,\beta)\in G^{\Omega_\sigma}(\mathbf A^\times_{X_\sigma})\times  G^{\Omega_\sigma}(\mathbf A^\times_{X_\sigma})$ such that $\divi^G(\pi_\sigma(\alpha))$ and  $\divi^G(\pi_\sigma(\beta))$ have no common component; then $\alpha\ast_i\beta=\pi_\sigma(\alpha) \ast\pi_\sigma(\beta)$.
\end{proposition}
\proof
Put $g_1=\pi_\sigma(\alpha)$   and $g_2=\pi_\sigma(\beta)$, then by proposition \ref{gff4} we can write $g_1=g_{1,0}+c_1$ and $g_2=g_{2,0}+c_2$ for, $g_{1,0},g_{2,0}\in G^{\Omega_\sigma}_0(X_\sigma)$, $c_1=\log(\xi(\alpha))$ and  $c_2=\log(\xi(\beta))$. An easy calculation shows that:
$$g_1\ast g_2=\sum_{x\in X_\sigma}\ord^G_x(g_{2,0})g_{1,0}(x)+c_1\sum_{x\in X_\sigma}\ord^G_x(g_{2,0})+c_2\sum_{x\in X_\sigma}\ord^G_x(g_{1,0})\,.$$
Then it is enough to note the following equalities:
$$
\ord_x^G(g_{1,0})=\ord_x^G(g_1)=v_x(a_x)\,,
$$
$$
\ord_x^G(g_{2,0})=\ord_x^G(g_2)=v_x(c_x)\,,
$$
$$
g_{1,0}(x)=g_1(x)-\log(\xi(\alpha))=-\log(b(x)a(x)\overline{a(x)})-\log(\xi(\alpha))\,.
$$
\endproof

Let's write an element $\alpha\in\mathbf A^\times_{\widehat X}=\mathbf A^\times_{X}\oplus\prod_{\sigma\in B_\infty}(\mathbf A^\times_{X_\sigma}\oplus \mathbf A^\times_{X_\sigma})$ in the following way:
$$\alpha=\alpha_X\times(\alpha_\sigma)_\sigma$$
with $\alpha_X\in\mathbf A^\times_X$ and $\alpha_\sigma\in \mathbf A^\times_{X_\sigma}\oplus \mathbf A^\times_{X_\sigma}$, then we have a surjective map:

\begin{eqnarray*}
\widehat{p}:\ker (d^1_\times)\oplus\prod_{\sigma}G(\mathbf A^\times_{X_\sigma})&\to&\Div(X)\oplus\bigoplus_{\sigma} G(X_\sigma)\\
\alpha=\alpha_X\times(\alpha_\sigma)_\sigma &\mapsto& \left(p(\alpha_X),\sum_\sigma \pi_\sigma(\alpha_\sigma)X_\sigma\right)
\end{eqnarray*}
where $p:\ker (d^1_\times)\to \Div(X)$ is the usual projection on usual divisors and $\pi_\sigma:  G(\mathbf A^\times_{X_\sigma})\to G(X_\sigma)$ is the projection on Green functions.
\begin{definition}\label{id_ar_pa}
Let's put  
$$\Div\left(\mathbf{A}^\times_{\widehat X}\right):={\widehat p}^{\;-1}\left(\Div_{\ar}(X,\Omega)\right)\,,$$
and let $\alpha,\beta\in\Div\left(\mathbf{A}^\times_{\widehat X}\right)$ such that $(\widehat p(\alpha), \widehat p(\beta))\in\Upsilon_{\ar}$
then the \emph{idelic Arakelov intersection pairing}
is given by:
$$\alpha.\beta:=\deg\left(\left<\alpha_X,\beta_X\right>_i\right)+\frac{1}{2}\sum_{\sigma}\varepsilon_\sigma\,\alpha_{\sigma}\ast_i\beta_\sigma$$
where $\deg$ is the usual degree of line bundles, $\left<\,,\,\right>_i$ is the idelic Deligne pairing and  $\alpha_{\sigma}\ast_i\beta_\sigma$ is the idelic $\ast$-product.
\end{definition}

We have to check that definition \ref{id_ar_pa} gives the correct extension of the Arakelov pairing.
\begin{theorem}
Let $\alpha,\beta\in\Div\left(\mathbf{A}^\times_{\widehat X}\right)$ such that $\widehat p(\alpha)=\widehat D$ and $\widehat p(\beta)=\widehat E$, with $(\widehat D,\widehat E)\in\Upsilon_{\ar}$, then $\alpha.\beta=\widehat D.\widehat E$. In other words the idelic Arakelov intersection pairing extends to a pairing:
\begin{eqnarray*}
\Div\left(\mathbf{A}^\times_{\widehat X}\right)\times \Div\left(\mathbf{A}^\times_{\widehat X}\right)&\to&\mathbb R
\end{eqnarray*}
and the following diagram is commutative:
$$
\begin{tikzcd}
\Div\left(\mathbf{A}^\times_{\widehat X}\right)\times \Div\left(\mathbf{A}^\times_{\widehat X}\right)\arrow[d,"{\hat p\times\hat p}"]\arrow[dr]& \\
\Div_{\ara}(X)\times \Div_{\ara}(X)\arrow[r]& \mathbb R
\end{tikzcd}
$$
\end{theorem}
\proof
It follows easily from the definitions.
\endproof

\begin{appendices}
\section{Semi-topological algebraic structures}\label{ap_ST}
\subsection{Basic notions}\label{ap_ST1}
\begin{definition}
A topological abelian group $(G,\tau)$  is  \emph{linearly topologised} (or has a \emph{linear topology}) if there is a local basis at $0$ made of subgroups. A morphism between linearly topologised groups is a continuous homomorphism. The category of linearly topologised group is denoted by $\catname{LTAb}$.
\end{definition} 

\begin{proposition}\label{lintop}
Let  $G$ be an abelian group and fix a non-empty collection of subgroups  $\mathcal F=\{U_i\}_{i\in I}$. If  $G$ is endowed with the topology $\tau$ generated by $\{x+U_i\}_{i\in I, x\in G}$,  then it becomes a linearly topologised group. 
\end{proposition}
\proof
First we show that $G$ is a topological group: we want the inversion $\iota: G\to G$ and the sum $\sigma: G\times G\to G$ to be continuous. We check this for the subbase $\{x+U_i\}_{i\in I, x\in G}$. Obviously $\iota^{-1}(U_i+x)=U_i-x\in\tau$. Then we prove that the following equality holds:
$$\sigma^{-1}(U_i+x)=\bigcup_{y\in G}(U_i+y)\times (U_i+x-y)\,.$$
The inclusion $\supseteq$ is evident, so let $(z,z')\in \sigma^{-1}(U_i+x)$, then $z=u+(x-z')$ for $u\in U_i$. If we write $z'=0+x-(x-z')$ and we put $y=x-z'$ we finally get $(z,z')=(u+y,0+x-y)\in (U_i+y)\times (U_i+x-y)$.\\
For the last statement consider the family 
$$\mathcal B:=\{U\in\tau\colon \text{ U is finite intersection of elements of } \mathcal F\}\,.$$
Then $\mathcal B$ is a local basis at $0$ made of subgroups.
\endproof
\begin{definition}
The linear topology on an abelian group $G$ obtained from a family of subgroups $\{U_i\}_{i\in I}$, as it is described in proposition \ref{lintop}, is called \emph{the linear topology generated by} $\{U_i\}_{i\in I}$.
\end{definition}
In this setting, concepts like initial and final topologies are well defined. Let $G$ be an abelian group and consider some homomorphisms of groups $\left\{\varphi_\alpha: G\to H_\alpha\right\}_\alpha$ and $\left\{\psi_\beta: H_\beta\to G\right\}_\beta$, where the $H_\alpha$ and $H_\beta$ are all linearly topologised. The \emph{initial linear topology} on $G$ with respect to $\left\{\varphi_\alpha\right \}_\alpha$  is the linear topology generated by 
$$\left\{\varphi^{-1}_\alpha(V_\alpha)\colon V_\alpha\subseteq H_\alpha \text{ is an open subgroup}\right\}_\alpha\,.$$
This is the coarsest linear topology which makes all the $\varphi_\alpha$ continuous. The \emph{final linear topology} on $G$ with respect to $\left\{\psi_\beta\right\}_\beta$  is the linear topology generated by 
$$\left\{U\subseteq G\colon U \text{ is a subgroup and } \psi_\beta^{-1}(U)  \text{ is open for any } \beta \right\}\,.$$
This is the finest linear topology which makes all the $\psi_\beta$ continuous. 
\begin{proposition}
$\catname{LTAb}$ is an additive category and moreover it admits inverse and direct limits.
\end{proposition}
\proof
The nontrivial statements are those involving the categorical limits. In particular  ${\varprojlim_ i} G_i$ and ${\varinjlim_j} G_j$ are the usual  limits in the category of groups, endowed respectively with the inital and final linear topology.
\begin{remark}
By commodity, in the category of linearly topologised groups, we call the limits ${\varprojlim_ i} G_i$ and ${\varinjlim_j} G_j$  respectively \emph{linear inverse limit} and \emph{linear direct limit}.
\end{remark}

\begin{definition}
 A \emph{ST ring}  (ST stands for semi-topological) is a ring $A$ endowed with a topology  satisfying the following two properties:
\begin{enumerate}
\item[\sbt] $(A,+)$ is a linearly topologised abelian group.
\item[\sbt] For any $a\in A$ the map $\lambda_a :A\to A,$ such that $\lambda_a(x)=ax$, is continuous.
\end{enumerate}
A morphism of ST rings is a continuous homomorphisms of rings. The category of ST rings is denoted as $\catname{STRing}$. Moreover  $B$ is a ST $A$-algebra if there is a morphism of ST rings $\varphi: A\to B$. The category of ST A-algebras is $A$-$\catname{STAlg}$.
\end{definition}

\begin{proposition}\label{stinvlim}
$\catname{STRing}$ and $A$-$\catname{STAlg}$ admit inverse and direct limits.
\end{proposition}
\proof
We show it only for rings. Let $A=\varprojlim_i A_i$ be the usual inverse limit in the category of rings and topologise its additive structure by taking the linear inverse limit topology. Thus we have the coarsest linear topology on $(A,+)$ such that the projections $\pi_j:A\to A_j$ are continuous. Assume that $\Lambda_{(a_i)}$ is the multiplication by $(\ldots,a_i,a_{i+1},\ldots)$ in $A$ and consider the composition: $A \xrightarrow{\Lambda_{(a_i)}} A\xrightarrow{\pi_j} A_j$, given by 
$$
x=(\ldots x_i,x_{i+1},\ldots) \mapsto  (\ldots a_ix_i,a_{i+1}x_{i+1},\ldots) \mapsto a_jx_j\,.
$$
Since $\pi_j\circ\Lambda_{(a_i)}(x)= \lambda_{a_j}\circ \pi_j(x)$, we can conclude that $\pi_j\circ\Lambda_{(a_j)}$ is continuous. Finally if $\pi_j^{-1}(V_j)\subset A$ is an element in the subbase of $A$, then  $\Lambda_{(a_i)}^{-1}(\pi_j^{-1}(V_j))$ is open in $A$.\\

Let $A=\varinjlim_i A_i$ be the usual direct limit in the category of rings and topologise its additive structure by taking the linear direct limit topology. Thus we have the finest linear topology on $(A,+)$ such that the maps $\phi_i:A_i\to A$ are continuous. Let's denote with $\mu_{ij}: A_i\to A_j$ the continuous homomorphims in the directed set $\{A_i\}_i$; moreover $\Lambda_{[(j,a)]}$ is the multiplication in $A=\left(\sqcup_i A_i\right)/_{\!\sim}$ for the fixed element $[(j,a)]$ where $
a\in A_j$. Note that the composition: $A_i \xrightarrow{\phi_i} A\xrightarrow{\Lambda_{[j,a]}} A$, given by 
$$
x \mapsto [(i,x)] \mapsto [k, \mu_{jk}(a)\mu_{ik}(x)]\,.
$$
is continuous. Thus if $U\subset A$ is open, then $\phi_1^{-1}\left(\Lambda_{[j,a]}^{-1}(U)\right)$ is open and by definition of final linear topology we can conclude that $\Lambda_{[j,a]}^{-1}(U)$ is open in $A$.
\endproof
\begin{definition}
 Let $A$ be a ST ring. A \emph{$ST$ $A$-module} is an $A$-module satisfying the following properties:
\begin{enumerate}
\item[\sbt] $M$ is a linearly topologised abelian group.
\item[\sbt] For any $a\in A$ and any $m\in M$ the maps $\lambda^M_a:M\to M$ and $\rho_m:A\to M$ such that $\lambda_a(x)=ax$ and $\rho_m(x)=xm$ are continuous. 
\end{enumerate}
A morphism of ST modules is a continuous homomorphism of $A$-modules. If  $A$ is a ST field then $M$ is called a \emph{ST vector space}.
\end{definition}
Given a ST $A$-module $M$, the subset $\overline{\{0\}}$ is  a submodule because of the continuity of $\lambda_a$, therefore we define
$$M^{\sep}:=M/\overline{\{0\}}$$
which is again a ST $A$-module if endowed with the quotient topology.
\begin{proposition}\label{finetop}
Let $A$ be a ST ring, and $M$ an $A$-module. If $M$ is endowed with the final linear topology with respect to the group homomorphisms $\rho_m: A\to M$, then $M$ is a ST $A$-module
\end{proposition}
\proof
See \cite[pag. 17]{yek1}.
\endproof
\begin{definition}
The topology on $M$ described in proposition \ref{finetop} is called the \emph{fine $A$-module topology}.
\end{definition}

\subsection{Ind/pro topologies}\label{indpro}
Now we present the crucial part of this very general theory. Given  a ST ring $A$, we describe two procedures called $(C)$ and $(L)$  that give canonical topologies of ST rings respectively on $\varprojlim_r A/\mathfrak p^r$ and $A_{\mathfrak p}$ for any prime ideal $\mathfrak p\subset A$. We need the following lemma:
\begin{lemma}\label{yeklem}
Let $\varphi:A\to B$  be a ring homomorphism where $A$ is a ST ring. Consider $B$ as an $A$-module endowed with the fine $A$-module topology, then $B$ is a ST ring.
\end{lemma}
\proof
\cite[Proposition 1.2.9.(b)]{yek1}.
\endproof
\begin{itemize}
\item[$(C)$] For any $r>0$ we put on $ A/\mathfrak p^r$ the fine $A$-module topology, so by lemma \ref{yeklem} $ A/\mathfrak p^r$ is a ST ring. By proposition \ref{stinvlim} we can endow $\varprojlim_r A/\mathfrak p^r$ with a structure of  ST ring. 

\item[$(L)$]  $A_\mathfrak p$ is naturally  an $A$-module, so we endow it with the fine $A$-module topology. Again by lemma \ref{yeklem} we conclude that $A_\mathfrak p$ is a ST ring.
\end{itemize}

Let $R$ be a ST ring ring and put on $A=R[t]$ the fine $R$-module topology. Consider the ring of formal Laurent power series $R((t))$, then as linear projective limit we have:  
$$R[[t]]=\varprojlim_r \frac{R[t]}{t^rR[t]}\,.$$

Therefore we consider on $R((t))$ the topology induced in the following way:
\begin{equation}
\begin{tikzcd}
A=R[t]\arrow[r,squiggly,"(C)"]& R[[t]] \arrow[r,squiggly,"(L)"]& R((t))\,.
\end{tikzcd}
\end{equation}
This is called the ind/pro-topology. We have an isomorphism of ST $R$-modules
$$R((t))\cong\left(\bigoplus_{n\in\mathbb N} R\right)\oplus\prod_{n\in\mathbb N} R$$
and each subgroup of the form $t^rR[[t]]$, for $r\in\mathbb Z$, is closed in $R((t))$ . 
\begin{remark}
If we start with a discrete field  $K=R$, then the ind/pro-topology on $R((t))$ is the discrete valuation topology.
\end{remark}
Let $\xi\in \widehat{R((t))}$ be a nontrivial character. The \emph{conductor of $\xi$} is $$c_{\xi}:=\min\left\{i\in\mathbb Z\colon \xi\in\left(t^iR[[t]]\right)^\perp\right\}\,.$$

\section{Arakelov geometry}\label{ar_g}
This section is just a collection of basic results about Arakelov geometry for arithmetic surfaces. We will maintain the same notations used so far for the arithmetic surface $\varphi:X\to B$. Moreover we assume the reader to be familiar with complex analytic geometry for Riemann surfaces.
\subsection{Green functions and \texorpdfstring{$\ast$}{}-product}
Let's fix a connected Riemann surface $C$.
\begin{definition}\label{defG}
A \emph{Green function} on $C$ is a map $g:U\subseteq C\to \mathbb R$ satisfying the following properties:
\begin{enumerate}
\item[$(1)$] $U=C\setminus\{x_1,\ldots,x_r\}$ for $r\in\mathbb N$.
\item[$(2)$] $g$ is a $C^\infty$ function on $U$.
\item[$(3)$] For any point $x\in C$ there exist a real number $a\in\mathbb R$ and a $C^\infty$ function $u$ on an open neighborhood of $x$ such that the equality:
$$g=a\log|z|^2+u$$
holds in an open  neighborhood of $x$ contained in a holomorphic chart $(V,z)$ centred in $x$.
\end{enumerate} 
\end{definition}
The number $a\in\mathbb R$ arising in condition $(3)$ of definition \ref{defG} depends only on the point $x$ and it is uniquely defined.
\begin{definition}
Let $g$ be a Green function on $C$ such that around a point $x\in C$ it can be written as 
$g=a\log|z|^2+u$. Then we put $\ord_x^G(g):=-a$ and we call it  \emph{the Green order of $g$ at $x$}.
\end{definition}
Clearly $\ord_x^G(g)\neq 0$ if and only if $x$ is a point out from the domain of $g$, i.e. only at a finite number of points. The Green functions on $C$ form a real vector space $G(C)$, and for any $g,g'\in G(C)$
$$\ord^G_x(\lambda g)=\lambda\ord^G_x(g)\quad\text{for any $\lambda\in\mathbb R$}\,, $$
$$\ord^G_x(g+g')=\ord^G_x(g)+\ord^G_x(g')\,.$$
Let's denote with  $\Div(C)_{\mathbb R}:=\Div(C)\otimes_{\mathbb Z}\mathbb R$ the vector space of $\mathbb R$ divisors on $C$, then we have a $\mathbb R$-linear map:
\begin{eqnarray*}
\divi^G: G(C)&\to& \Div(C)_{\mathbb R}\\
g &\mapsto &\sum_{x\in C}\ord^G_x(g)[x]
\end{eqnarray*}

\begin{proposition}\label{gf1}
Let $(\mathscr L,h)$ be a $C^{\infty}$ hermitian invertible sheaf on $C$, and let $s$ be a nonzero meromorphic section of $\mathscr L$, then the map $-\log (h(s,s))$ is a Green function on $C$ such that $\divi^G(-\log (h(s,s)))=\divi(s)$.
\end{proposition}
\proof
See \cite[lemma 4.8]{Mor}.
\endproof
The following result is an immediate consequence of proposition \ref{gf1}:
\begin{proposition}\label{gf2}
The map $\divi^G: G(C)\to \Div(C)_{\mathbb R}$ is surjective.
\end{proposition}

Let's define a very important subspace of $\mathbb ZG(C)$:

\begin{definition}
The vector space of Green functions with integer orders on $C$ is:
$$\mathbb ZG(C):=\left\{g\in G(C)\colon \ord^G_x(g)\in\mathbb Z\,\;\forall x\in C\right\}$$

\end{definition}
The next result shows that any Green function which induces a divisor on $C$ is actually of the form $-\log (h(s,s))$ for some meromorphic section $s$ of a $C^\infty$ hermitian invertible sheaf $(\mathscr L, h)$.
\begin{proposition}\label{Zgreen}
Let $g\in \mathbb ZG(C)$ , then there exist a $C^\infty$ hermitian invertible sheaf $(\mathscr L, h)$ on $C$ and a meromorphic section $s$ of $\mathscr L$ such that $g=-\log (h(s,s))$.
\end{proposition}
\proof
See again \cite[lemma 4.8]{Mor}.
\endproof

From now on, in this subsection we fix a K\"ahler fundamental form $\Omega$ on $C$ such that $\int_C\Omega=1$.  Let's define  some  subsets of $G(C)$:
$$G^\Omega(C):=\{g\in G(C)\colon \Delta_{\bar\partial}(g) \;\text{is constant}\}\,,$$
$$G^\Omega_0(C):=\{g\in G^\Omega(C)\colon \int_C g\Omega=0 \}\,,$$
$$\mathbb ZG^\Omega(C):=\mathbb ZG(C)\cap G^\Omega(C)\,,$$
$$\mathbb ZG^\Omega_0(C):=\mathbb ZG(C)\cap G^\Omega_0(C)\,.$$
\begin{theorem}\label{gf3}
The map $\divi^G|_{G^\Omega_0(C)}: G^\Omega_0(C)\to \Div(C)_{\mathbb R}$ is an isomorphism.
\end{theorem}
\proof
See \cite[theorem 4.10]{Mor}.
\endproof

\begin{proposition}\label{gff4}
For any $g\in G^\Omega(C)$ there exists a unique decomposition $g=g_0+c$  for  $g_0\in G^\Omega_0(C)$  and $c\in\mathbb R$.
\end{proposition}
\proof
See again \cite[theorem 4.10]{Mor}.
\endproof
\begin{definition}
 The inverse map of $\divi^G|_{G^\Omega_0(X)}$ is denoted as:
\begin{eqnarray*} 
 \mathcal G^\Omega:\Div(X)_\mathbb R &\to& G^\Omega_0(X)\\
 D &\mapsto& \mathcal G^{\Omega}(D)
 \end{eqnarray*}
 and we can define the following function:
 \begin{eqnarray*}
 g^{\Omega}:(X\times X)\setminus \Delta_{X\times X}&\to &\mathbb R\\
 (p,q)&\mapsto & g^{\Omega}(p,q):=\mathcal G^{\Omega}([p])(q)
 \end{eqnarray*}
 where $\Delta_{X\times X}$ denotes the diagonal subset of $X\times X$. 
\end{definition}
By construction $g^{\Omega}$ is $C^\infty$ in the variable $q$, but, as we will see soon (corollary \ref{symm}),  $g^{\Omega}$ turns out to be symmetric, therefore it is $C^\infty$. Since $g^\Omega(p,\cdot)\in G^\Omega_0(X)\subset  G^\Omega(X) $, then $dd^c(g^\Omega(p,\cdot))=\alpha\Omega$ for a constant $\alpha\in\mathbb C$, but
 $$1=\deg^{G}(g^\Omega(p,\cdot))=\int_X dd^c(g^\Omega(p,\cdot))=\int_X \alpha\Omega=\alpha\,.$$
Hence $\alpha=1$ and 
\begin{equation}\label{poi_eq}
dd^c(g^\Omega(p,\cdot))=\Omega.
\end{equation}
 Thus, amongst all Green functions, those of the form $g^\Omega(p,\cdot)$ satisfy the  Poisson differential equation (\ref{poi_eq}). this feature will be very useful for intersection theory.

 Another important property is that for any fixed $p\in X$:
\begin{equation}\label{boundarycond}
\int_{X}g^{\Omega}(p,\cdot)\Omega=\int_{X} \mathcal G^{\Omega}([p])\Omega=0
\end{equation}
because $\mathcal G^{\Omega}([p])\in G^\Omega_0(X)$.
\begin{remark}
$g^\Omega$ can be defined as the \emph{unique} function on $(X\times X)\setminus\Delta_{X\times X}$ with values in $\mathbb R$ satisfying the following properties:
\begin{itemize}
\item[$(1)$] Around any point $p\in X$ we can write $g^{\Omega}(p,\cdot)=-\log|z|^2+u$, where $z$ is a chart centred in $p$ and $u$ is $C^\infty$.
\item[$(2)$] $dd^c(g^\Omega(p,\cdot))=\Omega$.
\item[$(3)$] $\int_{X}g^{\Omega}(p,\cdot)\Omega=0$.
\end{itemize}
This is how Arakelov defined $g^\Omega$ in \cite{Ar} and \cite{Ar1}. In the literature $g^\Omega$ is usually called \emph{the Green function of $X$ (with respect to $\Omega$)}\footnote{Actually the conditions which uniquely define $g^\Omega$ in \cite{Ar} and \cite{Ar1} are slightly different from the ones listed here, and moreover they may vary in other references. For instance it is common to find different constants for  the differential Poisson equation, or the Green function might be defined as $G=\exp(g^\Omega)$. Of course these discrepancies are fixed when the Green function is applied for intersection theory.}. Here we used a different approach (and notations), indeed  $g^\Omega$ was constructed directly by using the isomorphism $\Div(X)_{\mathbb R}\cong G^{\Omega}_0(X)$.
\end{remark}
 
 \begin{definition}
Let $g_1,g_2\in G(C)$ such that $\divi^G(g_1)$ and $\divi^G(g_2)$ have no common components then the \emph{$\ast$-product} between $g_1$ and $g_2$ is the real  number:
$$g_1\ast g_2:=\widetilde {g_1}(\divi^G(g_2))+\int_{C}dd^c(g_1)g_2\,,$$
where $dd^c=\frac{1}{2\pi}\partial\overline \partial$.
\end{definition}
\begin{remark}
It is necessary to assume that $\divi^G(g_1)$ and $\divi^G(g_2)$ have no common components otherwise $\widetilde {g_1}(\divi(g_2))$ wouldn't be well defined.
\end{remark}
\begin{theorem}\label{astprod}
Let $g_1,g_2\in G(C)$ such that $\divi^G(g_1)$ and $\divi^G(g_2)$ have no common components, then $g_1\ast g_2=g_2\ast g_1$.
\end{theorem}
\proof
See \cite[proposition 4.12]{Mor}.
\endproof
\begin{corollary}\label{symm}
$g^{\Omega}(p,q)=g^{\Omega}(q,p)$ for any $p\neq q$.
\end{corollary}
\proof
By using the properties of the elements in $G^{\Omega}_0$, it is easy to verify that 
$$\mathcal G^{\Omega}([p])\ast \mathcal G^{\Omega}([q])= g^{\Omega}(p,q);\quad \mathcal G^{\Omega}([q])\ast \mathcal G^{\Omega}([p])= g^{\Omega}(q,p)\,.$$
Hence the conclusion follows immediately from theorem \ref{astprod}.
\endproof
Note that for any three different points $p,q,t\in X$ and  coefficients $a,b\in\mathbb R$ we have that:
$$\mathcal G^{\Omega}(a[p]+b[q])\ast \mathcal G^{\Omega}([t])=a\mathcal G^{\Omega}([p])\ast \mathcal G^{\Omega}([t])+b\mathcal G^{\Omega}([q])\ast \mathcal G^{\Omega}([t])\,.$$
Therefore if $D=\sum_{p\in X} a_p[p]$ and $E=\sum_{q\in X} b_q[q]$ are two real divisors of $X$ with no common components, then it is customary to define:
\begin{equation}\label{ArGondivi}
g^{\Omega}(D,E):=\sum_{p\neq q} a_pb_q g^{\Omega}(p,q)\,.
\end{equation}
\begin{remark}
The important point to emphasize here is that for Green functions $g_1,g_2\in G^{\Omega}_0$, i.e. coming from some real divisors on $X$, the integral appearing in $g_1\ast g_2$ vanishes. This means that for such kind of Green functions, the nature of the $\ast$-product is ``less analytic'', indeed it depends only on the value of $g_1$ or $g_2$ at a finite set of points. 
\end{remark}

\subsection{Arakelov intersection pairing}
On each Riemann surface $X_\sigma$ we fix a K\"ahler form $\Omega_\sigma$ such that $\int_{X_\sigma}\Omega_\sigma=1$, and we put $\Omega:=\{\Omega_\sigma\}_{\sigma\in B_\infty}$. For any divisor $D\in\Div(X)$,  $D_\sigma:=\varphi_\sigma^\ast D\in\Div(X_\sigma)$ denotes its pullback through $\varphi_\sigma$. Consider the additive group $\mathbf G(X):=\oplus_{\sigma\in B_\infty} G(X_\sigma)$ and its subgroup, depending on $\Omega$, $\mathbf G(X,\Omega):=\bigoplus_{\sigma\in  B_{\infty}} G^{\Omega_\sigma}(X_\sigma)$. By commodity 
we  write any element of $\mathbf G(X)$ (or of $\mathbf G(X,\Omega)$) as a finite formal linear combination $\sum_{\sigma} g_\sigma X_\sigma$ for $g_\sigma\in\mathbf G(X) $ (or $g_\sigma\in\mathbf G(X,\Omega)$). 
  
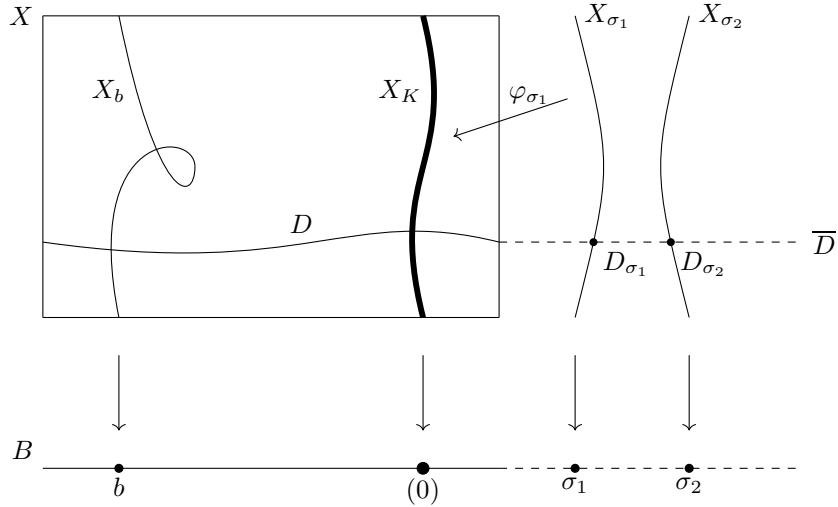
\begin{figure}[htp]
\centering
\begin{tikzpicture}

\draw  (-2,-2)  .. controls (-2.5,0.5) and (-1,0.5) .. (-1,0);
\draw (-1,0) .. controls (-1,-0.5) and (-1.5,-0.5) .. (-2,2); 
\draw (-3,2)node[anchor=east]{$X$}--(3,2);
\draw (-3,2)--(-3,-2);
\draw (-3,-2)--(3,-2);
\draw (3,2)--(3,-2);

\draw (-3,-1)..controls (0.5,-1.5) and (1,-0.5)..(3,-1);
\draw (0.4,-1)node[anchor=south]{$D$};

\draw[name path=c2,dashed](3,-1)--(7,-1)node[anchor=west]{$\overline D$};

\draw [line width=2.2pt](2, -2)..controls (1.5,0) and (2.5,0).. (2,2);
\draw (2.1,1)node[anchor=east]{$X_K$};
\draw (-1.8,1)node[anchor=east]{$X_b$};

\draw [name path=c1](4, -2)..controls (4.5,0) and (4.5,0).. (4,2)node[anchor=west]{$X_{\sigma_1}$};

\draw [name path=c3](5.5, -2)..controls (5,0) and (5,0).. (5.5,2)node[anchor=west]{$X_{\sigma_2}$};

\fill [name intersections={
of=c1 and c2,  by={b}}]
(b) circle (1.5pt)node[anchor=north west]{$D_{\sigma_1}$};

\fill [name intersections={
of=c3 and c2,  by={b}}]
(b) circle (1.5pt)node[anchor=north west]{$D_{\sigma_2}$};

\draw (-3,-4)node[anchor= south east]{$B$}--(3,-4);
\draw[dashed] (3,-4)--(7,-4);
\draw (2,-4) circle [radius=2.2pt];
\fill (2,-4) circle [radius=2.2pt]node[anchor=north]{$(0)$};

\draw [->] (-2,-2.5)--(-2,-3.5);
\draw [->] (2,-2.5)--(2,-3.5);
\draw [->] (4,-2.5)--(4,-3.5);
\draw [->] (5.5,-2.5)--(5.5,-3.5);

\draw [<-] (2.4,0.4)--(3.9,0.9);
\draw (3.4,0.7)node[anchor=south]{$\varphi_{\sigma_1}$};

\draw (-2,-4) circle [radius=1.5pt];
\fill (-2,-4) circle [radius=1.5pt]node[anchor=north]{$b$};

\draw (4,-4) circle [radius=1.5pt];
\fill (4,-4) circle [radius=1.5pt]node[anchor=north]{$\sigma_1$};

\draw (5.5,-4) circle [radius=1.5pt];
\fill (5.5,-4) circle [radius=1.5pt]node[anchor=north]{$\sigma_2$};

\end{tikzpicture}
\caption{footnotesize {A schematization of an arithmetic surfaces $\varphi:X\to B$ such that $B_\infty=\{\sigma_1,\sigma_2\}$ (for instance $B=\spec \mathbb Z[i]$). $X_b$ is a vertical divisor over the closed point $b$, $D$ is a horizontal divisor such that $D_{\sigma_1}$ and $D_{\sigma_2}$ are prime divisors respectively on $X_{\sigma_1}$ and $X_{\sigma_2}$.}}
\label{figu1}
\end{figure}

\begin{definition}\label{defardiv}
The group of \emph{Arakelov divisors} on $\widehat{X}$ is:
$$\Div_{\ar}(X,\Omega):=\left\{\left (D,\sum_{\sigma}g_\sigma X_\sigma\right)\in \Div(X)\times \mathbf G(X,\Omega)\colon \divi^G(g_\sigma)=D_\sigma\right\}\,.$$
We often denote the element $(0,X_\sigma)\in\Div_{\ar}(X,\Omega) $ simply with the symbol $X_\sigma$.
\end{definition}
It is important to understand the  geometry lying behind the above apparently mysterious definition. Fix an Arakelov divisor $(D,\sum_\sigma g_\sigma X_\sigma)$, by theorem \ref{gf3} and proposition \ref{gff4} we can write  
\begin{equation}\label{ar_dec}
g_\sigma=\mathcal G^{\Omega_\sigma}(D_\sigma)+\alpha_\sigma
\end{equation}
 where $\alpha_\sigma\in\mathbb R$ is uniquely determined. Figure \ref{figu1} highlights the fact that $D_\sigma$, which is a finite set of points on $X_\sigma$, can been interpreted as the ``prolongation'' of $D$ on the curve $X_\sigma$; thus, it makes sense to define the Arakelov divisor
$$\overline D:=\left(D,\sum_\sigma \mathcal G^{\Omega_\sigma}(D_\sigma)X_\sigma\right)\in \Div_{\ar}(X,\Omega)$$ 
which will be called \emph{completion} of $D$ in $\widehat{X}$ (this is consistent with the notion of completed horizontal curve given before). By equation (\ref{ar_dec}) we have the following unique decomposition  of $(D,\sum_\sigma g_\sigma X_\sigma)$ in $\Div_{\ar}(X)$:
\begin{equation}\label{ar_dec1}
\left(D,\sum_\sigma g_\sigma X_\sigma\right)=\overline D+\sum_\sigma\alpha_\sigma X_\sigma
\end{equation}
where the linear combination $\sum_\sigma\alpha_\sigma X_\sigma$ can be evidently read as a ``real divisor'' on $\widehat{X}$ with support made of  curves at infinity. In perfect analogy with the usual notion of divisor, equation (\ref{ar_dec1}) tells us that an Arakelov divisor can be interpreted as a formal linear combination of ``curves'' in $\widehat X$, such that the coefficients of the curves at infinity are in $\mathbb R$. The presence of this real coefficients  underlines once again the fact that the curves at infinity have an analytic nature.
From the above  discussion we recover the original definition of the group of Arakelov divisors given in \cite{Ar} and \cite{Ar1}:
\begin{proposition}
There is an isomorphism of groups:
$${\Div_{\ar}(X,\Omega)}\cong \Div(X)\oplus \mathbb R^{(B_\infty)}$$
\end{proposition}
\proof
Thanks to  equation (\ref{ar_dec1}) we can define the isomorphism:
$$\left(D,\sum_\sigma g_\sigma X_\sigma\right)\mapsto D+\sum_\sigma\alpha_\sigma[\sigma]\,.$$
\endproof

Now we want to introduce the concept of principal Arakelov divisor, in other words we want to define an Arakelov divisor associated to an element of $K(X)$. Recall that $K(X)$ is also the function field of $X_K$, so the morphism $\varphi_\sigma:X_\sigma\to X_K$ induces a field embedding 
$$\varphi^\#_{\sigma}:K(X)\hookrightarrow \mathbb C(X_\sigma)\,.$$ 
For any rational function $f\in K(X)$ we put by simplicity  $f_\sigma:=\varphi^\#_\sigma(f)$. Moreover let $\mathscr O_\sigma$ be the sheaf of regular functions on $X_\sigma$, then as usual $f_\sigma$  can  be identified with a holomorphic map $X_\sigma\to \mathbb C $ at all but finitely many points:
$$p\mapsto f_{\sigma,p}\mapsto \overline f_{\sigma,p}\in k(p)\cong\mathbb C$$
Then it is easy to see that  $-\log|f_\sigma|^2$ is a Green function on $X_\sigma$ such that $\partial\bar\partial(-\log|f_\sigma|^2)=0$, therefore $-\log|f_\sigma|^2\in G^{\Omega_\sigma}(X_\sigma)$. 

\begin{proposition}
Let $f\in K(X)^\times$, then $\divi^G(-\log|f_\sigma|^2)=(f)_\sigma$, where $(f)_\sigma$ is the pullback of the principal divisor $(f)$.
\end{proposition}
\proof
Fix a point $p\in X_\sigma$, let $x=\varphi_\sigma(p)$ and consider $f$ as a rational function on $X_K$. If $\varpi_\sigma$ is a local parameter in  $\mathscr O_{\sigma,p}$ and $\varpi$ is a local parameter in $\mathscr O_{X_K,x}$, then 

$$f_\sigma=\varpi_\sigma^{v_p(\varphi^\#_\sigma(\varpi)){v_x(f)}}u\quad \text{for } u\in\mathscr O_{\sigma,p}\,.$$
This implies that $\ord^G_p(-\log|f_\sigma|^2)=v_p(\varphi^\#_\sigma(\varpi)){v_x(f)}$, but $v_p(\varphi^\#_\sigma(\varpi))$ is precisely the ramification index $e_{\varphi_\sigma,p}$, hence $\ord^G_p(-\log|f_\sigma|^2)=e_{\varphi_\sigma,p}v_x(f)$. So, we finally have:

$$\divi^G(-\log|f_\sigma|^2)=\sum_{p\in X_\sigma} e_{\varphi_\sigma,p}v_{\varphi_\sigma(p)}(f)[p]=(f)_\sigma\,.$$
\endproof
Now the following definition makes sense:
\begin{definition}
Let $f\in K(X)^\times$ be a rational function. It induces an Arakelov divisor in the following way:
$$\widehat{(f)}:=\left((f),\sum_\sigma-\log|f_\sigma|^2X_\sigma\right)\in \Div_{\ar}(X,\Omega)\,.$$
The group
$$\Princ_{\ar}(X,\Omega):=\left\{\widehat{(f)}\,: f\in K(X)\right\}$$
is called the group of \emph{principal Arakelov divisor} and $\CH^1_{\ar}(X,\Omega):=\frac{\Div_{\ar}(X,\Omega)}{\Princ_{\ar}(X,\Omega)}$ is the \emph{Arakelov Chow group}. Two Arakelov divisor are said \emph{linearly equivalent} if they are contained in the same class in $\CH^1_{\ar}(X,\Omega)$.
\end{definition}

Moreover for any principal Arakelov divisor $\widehat{(f)}$ we get the following decomposition:
$$\widehat{(f)}=\overline{(f)}+\sum_\sigma\left(\int_{X_\sigma}- \log|f_\sigma|^2\Omega_\sigma\right) X_\sigma\,.$$

\begin{proposition}\label{nocommon}
Let $D,E$ be two finite divisors on $X$ with no common components, then for any $\sigma\in B_\infty$ the divisors $D_\sigma$ and $E_\sigma$  on $X_\sigma$ have no common components.
\end{proposition}
\proof
Omitted
\endproof
Let's denote as $\Upsilon_{\ar}\subset \Div_{\ar}(X,\Omega)\times \Div_{\ar}(X,\Omega)$ the set of couples of Arakelov divisors with no common components on $X$, then we can define the Arakelov intersection pairing on $\Upsilon_{\ar}$:
\begin{definition}
Let $\widehat D:=\left(D,\sum_\sigma g_\sigma X_\sigma\right),\widehat E:=(E,\sum_\sigma l_\sigma X_\sigma)$ be two Arakelov divisors such that $(\widehat{D},\widehat{E})\in\Upsilon_{\ar}$. Thanks to proposition \ref{nocommon} we can define  an Arakelov divisor on $B$:\,\footnote{Note that we assume $D$ and $E$ to have no common components in order to ensure that the $\ast$-product between green functions is well defined for any $\sigma\in B_\infty$.}
$$\left<\widehat{D},\widehat{E}\right>_{\ar}:=\left<D,E\right>+\sum_\sigma  g_\sigma\ast l_\sigma\,[\sigma]\;\;\in \Div_{\ar}(B)$$
where 

$$\left <D, E\right>:=\varphi_\ast i(D,E)=\sum_{x\in X} [k(x):k(\varphi(x))]\,i_{x}(D,E)\, [\varphi(x)]$$ 
and $\ast$ is the product between Green functions. If $\widehat d=\sum_{\mathfrak p\in B}n_{\mathfrak p}[\mathfrak p]+\sum_{\sigma\in B_{\infty}}\alpha_\sigma[\sigma]$ is an Arakelov divisor on the base $B$, its degree  is defined as:
$$
\deg_{\ar}(\widehat d):=\sum_{\mathfrak p\in B}n_{\mathfrak p}\log \mathfrak N(\mathfrak p)+\frac{1}{2}\sum_{\sigma\in B_{\infty}}\epsilon_\sigma\alpha_\sigma\,.
$$
In particular we use the notation $D.E:= \deg_{\ar}(\left<D,E\right>)$, and the \emph{Arakelov intersection number} of $\widehat D$ and $\widehat E$ is:
$$\widehat D.\widehat E:=\deg_{\ar}\left(\left<\widehat{D},\widehat{E}\right>_{\ar}\right)=D.E+\frac{1}{2}\sum_\sigma\epsilon_\sigma\,g_\sigma\ast l_\sigma\,\,\in \mathbb R\,,$$
\end{definition}
The following proposition summarizes some properties of the Arakelov intersection pairing:
\begin{proposition}\label{propar}
Let $(\widehat{D},\widehat{E}), (\widehat{D_j},\widehat{E_j})\in\Upsilon_{\ar}$ with $j=1,2$, then
\begin{itemize}
\item[$(1)$] $\widehat D.\widehat E=\widehat E.\widehat D$ (symmetry).
\item[$(2)$] $(\widehat{D_1}+\widehat{D_2}).(\widehat{E_1}+\widehat{E_2})=\sum_{j,k=1}^2 \widehat{D_j}.\widehat{E_k}$ ($\mathbb Z$-bilinearity).
\item[$(3)$] If $\widehat{D}=(D,\sum_\sigma g_\sigma X_\sigma)$ and $f\in K(X)^\times$ such that $(D,(f))\in \Upsilon$, then 
$$\left<\widehat{D},\widehat{(f)}\right>_{\ar}=\widehat{\left(N_D(f)\right)}\;\in\Princ_{\ar}(B)\,.$$
In particular $\widehat{D}.\widehat{(f)}=0$.
\end{itemize}
\end{proposition}
\proof
See \cite[section 4.4]{Mor}.
\endproof
The Arakelov intersection number can be extended to an intersection pairing on the whole $\Div_{\ar}(X,\Omega)$ and induces a natural intersection pairing on $\CH^1_{\ar}(X,\Omega)$.
\begin{proposition}
 The Arakelov intersection number extends to any two Arakelov divisors in $\Div_{\ar}(X,\Omega)\times\Div_{\ar}(X,\Omega) $ and moreover descends naturally to pairing  on $\CH^1_{\ar}(X,\Omega)\times \CH^1_{\ar}(X,\Omega)$.
\end{proposition}
\proof
See \cite[section 4.4]{Mor}.
\endproof
Now we interpret the Arakelov intersection pairing in a more geometric way by using the decomposition  given in equation (\ref{ar_dec1}). Fix two Arakelov divisors $\widehat D,\widehat E\in \Upsilon_{\ar}$, then we can write 
$$\widehat{D}=\overline{D}+\sum_{\sigma}\alpha_\sigma X_\sigma=\left(D, \sum_\sigma \mathcal G^{\Omega_\sigma}(D_\sigma)X_\sigma\right)+\left(0,\sum_\sigma \alpha_\sigma X_\sigma\right)\,,$$

$$\widehat{E}=\overline{E}+\sum_{\sigma}\beta_\sigma X_\sigma=\left(E, \sum_\sigma \mathcal G^{\Omega_\sigma}(E_\sigma)X_\sigma\right)+\left(0,\sum_\sigma \beta_\sigma X_\sigma\right)\,.$$
In order to find explicitly $\widehat D.\widehat E$, by bilinearity and symmetry of the intersection pairing it is enough to understand how calculate the following three terms:
\begin{enumerate}
\item[$(i)$] $\overline D.\overline E$; namely the intersection of two completed divisors.
\item[$(ii)$] $\overline D.(0,\sum_\sigma\beta_\sigma X_\sigma)$; namely the intersection between a completed divisor and a divisor at infinity. Clearly $(0,\sum_\sigma\alpha_\sigma X_\sigma).\overline E$ is obtained in the same way.
\item[$(iii)$] $(0,\sum_\sigma\alpha_\sigma X_\sigma).(0,\sum_\sigma\beta_\sigma X_\sigma)$; that is the intersection of divisors composed only by curves at infinity.
\end{enumerate}
For $(i)$ let's evaluate $\mathcal G^{\Omega_\sigma}(D_\sigma)\ast\mathcal G^{\Omega_\sigma}(E_\sigma)$. By the bare definition of the $\ast$-product and $g^{\Omega_\sigma}$:
$$\mathcal G^{\Omega_\sigma}(D_\sigma)\ast\mathcal G^{\Omega_\sigma}(E_\sigma)=g^{\Omega_\sigma}(D_\sigma,E_\sigma)+\int_{X_\sigma} dd^c\left( \mathcal G^{\Omega_\sigma}(D_\sigma)\right)\mathcal G^{\Omega_\sigma}(E_\sigma)\,,$$
but since $\mathcal G^{\Omega_\sigma}(D_\sigma),\mathcal G^{\Omega_\sigma}(E_\sigma)\in G^{\Omega_\sigma}_0(X_\sigma)$, it is straightforward  to verify that the integral on the right hand side is $0$. Therefore we get:
\begin{equation}\label{inter1}
\overline D.\overline E= D.E+\frac{1}{2}\sum_\sigma \epsilon_\sigma g^{\Omega_\sigma}(D_\sigma,E_\sigma)\,.
\end{equation}
In order to calculate $(ii)$ we need $\mathcal G^{\Omega_\sigma}(D_\sigma)\ast\beta_\sigma$:
$$\mathcal G^{\Omega_\sigma}(D_\sigma)\ast\beta_\sigma=\beta_\sigma\ast \mathcal G^{\Omega_\sigma}(D_\sigma)=\beta_\sigma\deg(D_\sigma)+\int_{X_\sigma}dd^c(\beta_\sigma)\mathcal G^{\Omega_\sigma}(D_\sigma)=\beta_\sigma\deg(D_\sigma)\,,$$
thus we obtain
\begin{equation}\label{inter2}
\overline D.(0,\sum_\sigma\beta_\sigma X_\sigma)=\frac{1}{2}\sum_\sigma \epsilon_\sigma\beta_\sigma \deg(D_\sigma)\,.
\end{equation}
Finally $(iii)$ is trivial since $\alpha_\sigma\ast\beta_\sigma=0$ and we have:
\begin{equation}\label{inter3}
(0,\sum_\sigma\alpha_\sigma X_\sigma).(0,\sum_\sigma\beta_\sigma X_\sigma)=0\,.
\end{equation}

\end{appendices}

\bibliographystyle{plain}
\bibliography{adII.bib}

\Addresses

\end{document}